\documentclass[12pt,a4paper]{amsart}
\usepackage{mystyle,mymacros}
\usepackage[all]{xy} 		
\usepackage{longtable}
\usepackage{tikz}
\usepackage{amsxtra,amsmath,amsfonts,amscd, amssymb, mathrsfs}

\makeatletter
\@namedef{subjclassname@2010}{%
 \textup{2010} Mathematics Subject Classification}
\makeatother

\title[Diophantine geometry and uniformity]{Diophantine and tropical geometry, and uniformity of rational points on curves}

\author{Eric Katz}
\address{Eric Katz, Department of Mathematics, The Ohio State University University, 231 W. 18th Ave., Columbus, OH 43210}
\email{katz.60@osu.edu}

\author{Joseph Rabinoff}
\address{Joseph Rabinoff, School of Mathematics, Georgia Institute of Technology, Atlanta, GA 30332-0160, USA}
\email{jrabinoff@math.gatech.edu}

\author{David Zureick-Brown}
\address{David Zureick-Brown, Dept. of Math and CS, Emory University, 400 Dowman
 Dr., W401, Atlanta, GA 30322, USA}
\email{dzb@mathcs.emory.edu}

\newcommand{\defi}[1]{\textbf{#1}} 				
 
 \newcommand{\calX}{\mathcal{X}}
 
\def\presuper#1#2%
  {\mathop{}%
   \mathopen{\vphantom{#2}}^{#1}%
   \kern-\scriptspace%
   #2}
\newcommand{\Ab}{{\operatorname{Ab}}}
\newcommand{\BC}{{\operatorname{BC}}}
\newcommand\BCint{\presuper\BC\int}
\newcommand\Abint{\presuper\Ab\int}
\newcommand{\Log}{\operatorname{Log}}
\newcommand{\GSp}{\operatorname{GSp}}

\begin{document}

\begin{abstract}

We describe recent work connecting combinatorics and tropical/non-Archimedean geometry to Diophantine geometry, particularly the uniformity conjectures for rational points on curves and for torsion packets of curves. The method of Chabauty--Coleman lies at the heart of this connection, and we emphasize the clarification that tropical geometry affords throughout the theory of $p$-adic integration, especially to the comparison of analytic continuations of $p$-adic integrals and to the analysis of zeros of integrals on domains admitting monodromy.
\end{abstract}


\maketitle

\section{Introduction}
\label{sec:introduction}

\defi{Diophantine geometry} studies solutions to systems of equations over arithmetically interesting fields (e.g., $\Q, \F_p$); cornerstones of the subject include:
\begin{enumerate}
\item \defi{Faltings' theorem}: finiteness of the set $X(\Q)$ of solutions to equations defining an algebraic curve $X$ of genus $g \geq 2$ (e.g., $y^2 = f(x)$ with $f$ squarefree and of degree at least 5);
\item the \defi{Weil conjectures}: the number of $\F_{p^n}$-valued solutions to a system of equations $X$ are governed by strict and surprising formulas (which depend on the topology 
of $X$);
\item \defi{Mazur's theorem}: the solutions to $E\colon y^2 = f(x)$ with $f$ a squarefree cubic and $x,y \in \Q$ form an abelian group $E(\Q)$, and the possibilities for the torsion subgroup of $E(\Q)$ have been completely classified.
\end{enumerate}

A central outstanding conjecture in Diophantine geometry is the following.
\begin{conjecture}[Uniformity]
 There exists a constant $B(g,\Q)$ such that every smooth curve $X$ over $\Q$ of genus $g\geq 2$ has at most $B(g,\Q)$ rational points. 
\end{conjecture}
The uniformity conjecture famously follows \cite[Theorem 1.1]{CaporasoHM:uniformity} from the widely believed conjecture of Lang--Vojta, which is the higher dimension analogue of Faltings' theorem. Initially, the uniformity conjecture was considered somewhat outrageous; this implication was initially taken as possible evidence \emph{against} the Lang--Vojta conjecture and led to a frenzied hunt for examples (or better: families) of curves with many points.



\subsection{Bounds on $\#X(\Q)$}

The proofs of Mordell due to Faltings, Vojta, and Bombieri \cite{Faltings:bookArithmeticGeometry, Vojta:siegels-theorem-in-the-compact-case, Bombieri:mordell-conjecture-revisited} give various astronomical upper bounds on the quantity $\#X(\Q)$. (See \cite{szpiro:peu} for a discussion of these bounds which involve a lattice theoretic estimate applied to $\Jac_X(\Q)\otimes \R$ and incorporates constants involving the isogeny class of $\Jac_X$ and various height functions on $X$.)
  An earlier partial proof due to Chabauty \cite{Chabauty1941} was later refined by Coleman to give a more modest upper bound. See~\S\ref{sec:chabauty} for a longer discussion of this method.

\begin{theorem}[Coleman, \cite{Coleman:effectiveChabauty}]
\label{T:Coleman-bound}
 Let $X$ be a curve of genus $g$ and let $r = \rank_{\Z} \Jac_X(\Q)$. Suppose $p > 2g$ is a prime of {\bf good reduction}. Suppose $r < g$. Then
 $$\#X(\Q) \leq \#X(\F_{p}) + 2g - 2.$$
\end{theorem}

In some cases, this bound allows one to compute the set $X(\Q)$ (rather than its cardinality) exactly. For example, Gordon and Grant \cite{GordonG:1993-computing-the-mordell-weil} showed that the Jacobian of the genus 2 hyperelliptic curve
        \begin{equation}\label{eq:eg.curve}
	X\colon y^2 = x(x-1)(x-2)(x-5)(x-6)
      \end{equation}
has rank 1. A quick computer search finds that
\[
X(\Q) \supset \{ (0, 0), (1, 0), (2, 0), (5, 0), (6, 0), (3, \pm 6), (10, \pm 120), \infty
\}.
\]
One can check that $X$ has good reduction at $p = 7$. Applying Coleman's theorem
gives an upper bound $X(\Q) \leq 10$, and in particular rigorously proves that
the known points of $X(\Q)$ constitutes the full collection of rational points.
In other words, we have completely solved the degree-$5$ Diophantine
equation~\eqref{eq:eg.curve}, which, as always in arithmetic geometry, is a
demonstration of the power of the Chabauty--Coleman method.

This example is very special -- the bound $\#X(\F_{p}) + 2g - 2$ will be usually much larger than the true value of $\#X(\Q)$. Nonetheless, refinements to this method frequently allow one to provably compute $X(\Q)$ (with computational assistance from e.g., MAGMA); see \cite{BruinS:2010-mordell-weil-sieve} as a starting point.
\\

Gradual improvements to Chabauty and Coleman's method by various authors have refined Coleman's bound and removed some hypotheses from Coleman's theorem. 

Lorenzini and Tucker removed the hypothesis that $X$ have good reduction at the prime $p$.

 \begin{theorem}[Lorenzini, Tucker, \cite{LoreniniT:thue}, Corollary 1.11] 
 \label{T:PMC}
 Suppose $p > 2g$ and let $\cX$ be a proper regular model of $X$ over $\Z_p$. Suppose $r < g$. Then
 $$\#X(\Q) \leq \#\cX_{\F_p}^{\sm}(\F_{p}) + 2g - 2$$
 where $\cX_{\F_p}^{\sm}$ is the smooth locus of the special fiber $\cX_{\F_p}^{\sm}$.
 \end{theorem}
\textsc{}
The ``$2g-2$'' term of Coleman's theorem arises from application of Riemann--Roch on $X_{\F_p}$ and is the source of his ``good reduction'' hypothesis. Lorenzini and Tucker's generalization recovers the $2g-2$ term via Riemann--Roch on $X_{\Q_p}$ and a more involved $p$-adic analytic argument. A later, alternative proof \cite[Theorem A.5]{McCallumP:chabautySurvey} instead recovers the $2g-2$ term via arithmetic intersection theory on $\cX$ and adjunction.

In \cite{LoreniniT:thue}, the authors ask if one can refine Coleman's bound when the rank is small (i.e., $r \leq g-2$), which was subsequently answered by Stoll.

 \begin{theorem} [Stoll, \cite{Stoll:2006-independence-of-rational-points}, Corollary 6.7]
 \label{T:stollthm}
 With the hypothesis of Theorem \ref{T:Coleman-bound},
 \[
 \#X(\Q) \leq \#X(\F_{p}) + 2r.
 \]
 \end{theorem}
Stoll's idea was, instead of using a single integral, to tailor his choice of integral to each residue class. The geometric input that coordinates the information between residue classes and allows the refinement from $2g-2$ to $2r$ is \textbf{Clifford's theorem} for the reduction $X_{\F_p}$ and is the source of the good reduction hypothesis. 

Unlike Lorenzini and Tucker's generalization of Coleman's theorem, where they
replace Coleman's use of Riemann--Roch on $X_{\F_p}$ with Riemann--Roch on
$X_{\Q_p}$, it does not seem possible to replace Stoll's use of Clifford's
theorem on $X_{\F_p}$ with Clifford's theorem on $X_{\Q_p}$. Matt Baker
suggested that it might be possible to generalize Stoll's theorem to curves with
bad, \emph{totally degenerate} reduction (i.e., $X_{\F_p}$ is a union of
rational curves meeting transversely) using ideas from tropical geometry (see
the recent survey \cite{bakerJ:degeneration-survey} on tropical geometry and
applications), in particular the notion of ``chip firing'', Baker's combinatorial definition of rank, and Baker--Norine's \cite{BakerN:RR} combinatorial Riemann--Roch and Clifford theorems. Baker was correct, and in fact an enrichment of his theory led to the following common generalization of Stoll's and Lorenzini and Tucker's theorems.

\begin{theorem}[Katz, Zureick-Brown, \cite{katzZB:tropicalChabauty}]
\label{T:katzB}
 Let $X/\Q$ be a curve of genus $g$ and let $r = \rank \Jac_X(\Q)$. Suppose $p > 2r+2$ is a prime, that $r < g$, and let $\cX$ be a proper regular model of $X$ over $\Z_{p}$. Then
 \[
	 \#X(\Q) \leq \#\cX_{\F_p}^{\text{sm}}(\F_{p}) + 2r. 
 \]
\end{theorem}
For example, applying this theorem (with $p = 5$) to the rank 1, genus 3 hyperelliptic curve 
	\[
	X\colon -2 \cdot 11 \cdot 19 \cdot 173 \cdot y^2 = (x-50)(x-9)(x-3)(x+13)(x^3 + 2x^2 + 3x + 4)
	\]
gives that 
\[
X(\Q) = \{ \infty, (50,0), (9,0),(3,0),(-13,0), (25 , 20247920 ), (25 , -20247920 )
\}.
\]
Theorem \ref{T:katzB} is proved via a combination of \defi{combinatorics}, \defi{$p$-adic analysis}, and \defi{tropical geometry}. One can associate to a singular curve with transverse crossings its dual graph as in Figure~\ref{fig:skeleton}: component curves become nodes, and intersections correspond to edges.

Baker's detailed study \cite{Baker:specialization} of the relationship between
\emph{linear systems} on curves and on finite graphs includes a semicontinuity
theorem for ranks of linear systems (as one passes from the curve to its dual
graph), and purely graph theoretic analogues of Riemann--Roch and Clifford's
theorem. Baker's theory works best with totally degenerate curves (i.e., Mumford
curves, where each component of the special fiber is a $\P^1$), and the heart of Theorem \ref{T:katzB} is an enrichment of Baker's theory which accounts for the additional geometry when the components have higher genus.
\\

Another success of tropical tools is the following. Instead of rational points, one can also apply a variant of Chabauty and Coleman's method to symmetric powers of curves and study points of $X$ defined over \emph{some} number field $K \supset \Q$ of degree at most $d$. Building on earlier partial bounds \cite{Klassen:1993-algebraic-points-low-degree} and explicit approaches \cite{Siksek:2009-symmetric-power-chabauty}, Jennifer Park recently proved the following.

\begin{theorem}[Park, \cite{Park-symmetric-chabauty-arxiv}]
 Let $d \geq 1$, $p$ a prime, and $g \geq 2$. Then there exists a number $N(p,g,d)$ such that for curve $X/\Q$ with good reduction at $p$, with $\rank_{\Z} \Jac_X(\Q) \leq g-d$, and satisfying an additional technical condition (related to excess intersection), 
\[
\#\big\{P \in X(K) \mid \deg K \leq d \text{ and } P \text{ does not belong
  to the ``special set'' of } X\big\}
\]
is at most $N(p,d,g)$.
\end{theorem}

The number $N(p,g,d)$ can be effectively computed; for example, if one restricts
to odd degree hyperelliptic curves, then $N(2,3,3) \leq 1539$. The significant new idea in her work is to apply \textbf{tropical intersection theory} to the multivariate power series which arise in symmetric power Chabauty.

\subsection{Bhargavology}

There is an expectation that ``most'' of the time a curve will have as few points as possible. Similarly, one expects that the Mordell--Weil rank of the Jacobian of a curve is usually 0 or 1. In practice one needs a family of curves over a rational base to even make this precise, and thus often restricts to families of hyperelliptic (or sometimes low genus plane) curves.

Recent years have witnessed rapid progress on this flavor of problem; results on
ranks include \cite{bhargava2013average} (elliptic curves),
\cite{bhargava2012average} (Jacobians of hyperelliptic curves), and
\cite{Thorne:E6-arithmetic} (certain families of plane quartics). Combining
these rank results with Chabauty's method and other techniques, several recent
results
\cite{PoonenS:mostOddDegree,bhargava2013most,shankar2013average,bhargava2013pencils}
prove that a strong version of the uniformity conjecture (in that there are no ``non-obvious'' points) holds for ``random curves'' in such families; see \cite{ho:howManySurvey} for a recent survey.

\subsection{Uniform bounds for rational points on curves}

While \cite{katzZB:tropicalChabauty} was notable as one of the first substantial
applications of tropical methods to a Diophantine problem, it is far from a
proof of uniformity -- it requires $r < g$ (which is expected to hold for a
random curve, but not all curves); more crucially, $\#\cX_{\F_p}^{\text{sm}}(\F_{p})$ may be arbitrarily large. In fact, there are even hyperelliptic curves whose regular models have many smooth $\F_{p}$-points: a hyperelliptic curve with some ramification points sufficiently $p$-adically close together will possess a regular model whose special fiber will have long chains of $\P^1$'s.
\\

A recent breakthrough of M. Stoll \cite{stoll:uniform} removed, for hyperelliptic curves, the dependence from Chabauty--Coleman on a \emph{regular} model, and derived (for $r \leq g-3$) a \emph{uniform} bound on $\#X(\Q)$.
\begin{theorem}[Stoll, \cite{stoll:uniform}]
\label{T:stoll-uniform-hyperelliptic}
 Let $X$ be a smooth hyperelliptic curve of genus $g$ and let $r = \rank_{\Z} \Jac_X(\Q)$.
 Suppose $r\, \leq \, g-3$. Then 
 $$\#X(\Q) \leq 8(r+4)(g-1) + \max\{1,4r\}\cdot g.$$
\end{theorem}

Exploiting the full catalogue of tropical and non-Archimedean analytic tools, we overcome Stoll's hyperelliptic restriction and prove uniform bounds for \emph{arbitrary} curves of small rank. 

\begin{theorem}[Katz--Rabinoff--Zureick-Brown, \cite{KatzRZB-uniform-bounds}]
\label{T:KRZB}
 Let $X$ be {\textbf{any}} smooth curve of genus $g$ and let $r = \rank_{\Z} \Jac_X(\Q)$.
 Suppose $r\, \leq \, g-3$. Then
 $$\#X(\Q) \leq 84g^2 - 98g + 28.$$
\end{theorem}

We discuss an overview of the proofs of Theorems
\ref{T:stoll-uniform-hyperelliptic} and \ref{T:KRZB}
in~\S\ref{ss:stoll-discussion} and~\S\ref{ss:KRZB-prelude};
\S\ref{sec:berk-curv-skel} and~\S\ref{sec:integration} are dedicated to an
exposition of the tropical and non-Archimedean analytic details underlying the proof of Theorem \ref{T:KRZB}, and \S\ref{sec:uniformity-results} explains how these tools tie together to give a proof of Theorem \ref{T:KRZB}.

\subsection{Effective Uniform Manin-Mumford}
\label{ss:effective-manin-mumford}

The following pair of results in the direction of ``uniform Manin-Mumford''
also make essential use of the tropical toolkit. \

Declare two points $P,Q$ to be equivalent if $mP$ is linearly equivalent to $mQ$
on $X$ for some integer $m\geq 1$, and define a \defi{torsion packet} to be an
equivalence class under this relation. Equivalently, a torsion packet is the
inverse image of the torsion subgroup of the Jacobian of $X$, under an
Abel--Jacobi map.
Raynaud \cite{raynaud:maninMumford} famously proved that every torsion packet of a curve is finite. Many additional proofs, with an assortment of techniques and generalizations, were given later by \cite{buium:p-adic-jets, coleman:ramified_torsion_curves, Hindry1988,
 Ullmo1998,Pila2008} and others. Several of these proofs rely on $p$-adic methods, with Coleman's approach via $p$-adic integration especially related to ours.

A \emph{uniform} bound on the size of the torsion packets of a curve of genus
$g\geq 2$ is expected but still conjectural.  It is known to follow from the Zilber--Pink conjecture (see \cite[Theorem 2.4]{stoll:uniform}), which is a generalization of the Andre--Oort conjecture. The first, unconditional result concerns \emph{rational} torsion packets.

\begin{theorem}[\cite{KatzRZB-uniform-bounds}]
\label{T:uniformity-K-torsion}
Let $g\geq 3$. Then for any smooth, proper, geometrically connected genus-$g$ curve $X/\Q$, and any
Abel--Jacobi embedding $\iota\colon X \hookrightarrow J$ into its Jacobian, we have
\[
\#\iota^{-1}\big(J(\Q)_{\tors}\big) \leq 84g^2 - 98g + 28.
\]
\end{theorem}

 No such uniformity result was previously known. The second result concerns \emph{geometric} torsion packets and requires the following restriction on the reduction type. Let $X$ be a smooth, proper, geometrically connected curve of genus $g\geq 2$ over $\Q$. Let $p$ be a prime and let $\calX$ be the stable model of $X$ over $\overline{\Z}_p$. For each irreducible component $C$ of the mod $p$ reduction $\calX_s$ of $\calX$ let $g(C)$ denote its geometric genus and let $n_C$ denote the number points of the normalization of $C$ mapping to nodal points of $\calX_s$. We say that $X$ \defi{satisfies condition} ($\dagger$) \defi{at $p$} provided that
\begin{equation*}
 g > 2g(C) + n_C
\end{equation*}
for each component $C$ of $\calX_s$.

\begin{theorem}
\label{T:uniformity-torsion}
Let $g\geq 4$. Then for any smooth, proper, geometrically connected genus $g$ curve $X/\Q$ which satisfies condition \textup{($\dagger$)} at some prime $p$ of $\Q$, and for any Abel--Jacobi embedding $\iota\colon X_{\bar \Q}\hookrightarrow J_{\bar \Q}$ of $X_{\bar \Q}$ into its Jacobian, we have
\[
\#\iota^{-1}\big(J(\bar \Q)_{\tors}\big) \leq (16g^2-12g)\,N_p\big((4\cdot 7^{2g^2+g+1})^{-1},\,2g-2\big)
\] 
where 
\begin{equation*}\label{eq:explicit.geom.bound}
 N_p(s,N_0) = \min\big\{N\in\Z_{\geq 1}~:~s(n-N_0)>\lfloor \log_p(n) \rfloor ~\forall
 n\geq N\big\}. 
\end{equation*}
\end{theorem}

The condition ($\dagger$) is satisfied at $p$, for instance, when $X$ has
totally degenerate trivalent stable reduction.

A uniform bound as in Theorem~\ref{T:uniformity-torsion} for the size of
geometric torsion packets was previously known~\cite{buium:p-adic-jets} for
curves of \defi{good reduction} at a fixed prime $p$. Theorem
\ref{T:uniformity-torsion} uses Coleman's
observation~\cite{coleman:ramified_torsion_curves} that the $p$-adic integrals
considered in Chabauty's method vanish, unconditionally, on torsion
packets. Coleman additionally deduces uniform bounds in many situations, still
in the good reduction case: for instance, if $X/\Q$ has ordinary good reduction
at $p$ and its Jacobian $J$ has potential CM, then
$\#\iota^{-1}\big(J(\bar\Q)_{\tors}\big) \leq
gp$. Theorem~\ref{T:uniformity-torsion}, on the other hand, applies to curves
with highly \emph{degenerate} reduction, hence approaches the uniform
Manin--Mumford conjecture from the other extreme. 

The full power of the general machinery developed in \cite{KatzRZB-uniform-bounds} is needed for the proof of Theorem~\ref{T:uniformity-torsion}, which is striking in that it uses $p$-adic integration to bound the number of \emph{geometric} (i.e., over $\C$) torsion points. 

\section{The method of Chabauty--Coleman}
\label{sec:chabauty}

Let $X$ be a smooth, proper, geometrically connected curve of genus $g\geq 2$,
and let $J = \Jac_X$ be its Jacobian.
Building on an idea of Skolem \cite{skolem-1934-verfahren}) (who studied
products of the multiplicative group via $p$-adic methods), Chabauty
\cite{Chabauty1941} gave the first substantial progress toward the Mordell
conjecture, conditional on the hypothesis that the rank $r$ of the group
$J(F)$ of $F$-points of the Jacobian of $X$ is strictly less than the genus
$g$, where $F$ is a number field. See \cite{McCallumP:chabautySurvey} for a marvelous survey of the
method. We briefly exposit Chabauty's ideas here and highlight various details
of the proofs of Theorems \ref{T:Coleman-bound}, \ref{T:PMC}, \ref{T:stollthm},
\ref{T:katzB}, \ref{T:stoll-uniform-hyperelliptic}, and \ref{T:KRZB}

Chabauty's theorem considers the effect of $p$-adic integration on rational points. For brevity we consider only $F = \Q$.
\begin{theorem}[Chabauty, \cite{Chabauty1941}]
If $r<g$, then $X(\Q)$ is a finite set.
\end{theorem}
The finiteness argument goes as follows. One very naturally considers the commutative diagram
\[
\xymatrix{
X(\Q) \ar[r] \ar[d] & J(\Q)\ar[d]\\
X(\Q_p) \ar[r] & J(\Q_p).
}
\]
The $p$-adic closure $\overline{J(\Q)}$ of $J(\Q)$ inside the $p$-adic manifold $J(\Q_p)$ (a manifold in the very down-to-earth sense of Bourbaki \cite{Bourbaki-lie-groups-1-3}) satisfies the dimension bound
\begin{equation}
 \label{eq:dimension-bound}
 \dim \overline{J(\Q)} \leq r
\end{equation}
\cite[Lemma 4.2]{McCallumP:chabautySurvey}. In particular, if $r < g$, then the $p$-adic closure $\overline{J(\Q)}$ is a submanifold of smaller dimension; the intersection $X(\Q_p) \cap \overline{J(\Q)}$ can then be shown to be zero dimensional, hence discrete and compact, and thus finite. 

To deduce finiteness of the intersection, Chabauty constructed a locally analytic function $f_\omega$ (actually, a $p$-adic integral) near $P\in X(\Q)$ such that if $P,Q\in X(\Q)$ specialize to the same $\F_p$-point in a model of $X$ over $\Z_p$, then $f_\omega(P)=f_\omega(Q)$.  Because an analytic function cannot take a value infinitely often, there can only be finitely many rational points specializing to a given $\F_p$-point, hence there are only finitely many rational points. 

The function $f_\omega$ arises in the following fashion.  The proof of the
dimension bound~\eqref{eq:dimension-bound}  passes to the Lie algebra 
\[
 \Lie J_{\Q_p} \cong H^0(J_{X,\Q_p},\Omega^1)^{\vee} \cong \Q_p^g,
\]
 whence
 \begin{align*}
 \dim \overline{J(\Q)} = \dim \overline{\log J(\Q)} = \rank _{\Z_p} (\Z_p \log J (\Q)) \leq & \rank _\Z \log J (\Q) \\
\leq & \rank _\Z J (\Q) = r.   
 \end{align*}
The essential facts packaged in the first two equalities are that 
\begin{enumerate}
\item $\log$ is a local diffeomorphism, and 
\item the closure of a subgroup of a $\Q_p$-vector space is its $\Z_p$-span.
\end{enumerate}
Now, under the hypothesis $r < g$, the codimension of $\overline{\log J(\Q)}$ is at least $g-r$, and thus there exists a subspace $V$ of linear functionals on $\Lie J(\Q_p)$ of dimension at least $g-r$ which vanish on $\overline{J(\Q)}.$    

\subsection{Integration.} 
\label{ss:integration-warmup}
The functions $f_{\omega}$ arise as integrals as follows. A $1$-form $\omega\in V$ can be expanded in a power series in a uniformizer at a smooth point $P$ of $X(\Q)$.   The term-by-term antiderivative of $\omega$ is an analytic function $f_\omega$ near $P$ that obeys $f_\omega(P)=f_\omega(Q)$ for all $Q\in X(\Q)$ near $P$. 

More concretely, when $Q,Q' \in X(\Q_p)$ reduce to the same point of $X(\F_p)$, the integrals $\int_Q^{Q'} \omega$ are explicit and rather non-exotic objects (e.g., one can work quite directly with them in MAGMA and SAGE). The key observations leading to their definition are the following.
\begin{enumerate}
\item (\textbf{Local structure}) Given $P \in X^{\sm}(\F_p)$, the \textbf{tube} $D_P \subset X(\Q_p)$ of points reducing to $P$ is a $p$-adic disc -- given any $Q \in D_P$, any uniformizer $t$ at $Q$ defines an analytic isomorphism $D_P \cong p\Z_p$.
\item (\textbf{de Rham}) A $p$-adic disc has trivial de Rham cohomology, and so the restriction of any $\omega$ to $D_P$ admits a power series expansion
\[
\omega|_{D_p} = \sum_{i= 0}^{\infty} a_it^idt \in \Z_p\ps t dt.
\]
\item (\textbf{First fundamental theorem of calculus}) For $Q_1,Q_2 \in D_P$, the integral is determined \emph{formally} as 
\[
\int_{Q_1}^{Q_2} \omega := \int_{t(Q_1)}^{t(Q_2)} \sum_{i= 0}^{\infty}  a_it^idt = \left.\left( \sum_{i= 0}^{\infty}  \frac{a_i}{i+1}t^{i+1}\right)\right|^{t(Q_2)}_{t(Q_1)}.
\]

\end{enumerate}

\begin{eg}[From the survey \cite{McCallumP:chabautySurvey}]
 Consider the genus 2 curve
 \[
 X\colon y^2 = x^6 + 8x^5 + 22x^4 + 22x^3 + 5x^2 + 6x + 1.
 \]
A quick search with a computer gives 
\begin{enumerate}
 \item[] $X(\Q) \supset \{\infty^+, \infty^-, (0,\pm 1), (-3, \pm 1) \}$, and 
 \item[] $X(\F_3) = \{\infty^+, \infty^-, (0, \pm 1) \}$.
 \end{enumerate}
 Points $(x,y)$ reducing to $\widetilde{Q} = (0,1)$ are given by 
\[
\begin{array}{rl}
 x = &3\cdot t, \text{ where } t \in \Z_p \\
 y = & \sqrt{x^6 + 8x^5 + 22x^4 + 22x^3 + 5x^2 + 6x + 1} = 1 + x^2 + \cdots
 \end{array}
\]
and integration of the form $\omega = xdx/y$ is given by 
\[
\displaystyle{\int_{(0,1)}^{t(P)}\frac{xdx}{y} = \int_{0}^{t} (x - x^3 + \cdots)dx}.
\]

\end{eg}

This explains how to integrate between pairs of points $Q, Q'$ in the same tube $D_P$ (which is all that is needed for Chabauty's original proof, and for Coleman's original theorem); Coleman was able to define $f_\omega$ on all of $X$ (i.e., determine $\int_Q^{Q'} \omega$ for $Q, Q'$ in different tubes, a sort of ``$p$-adic analytic continuation'') via more involved techniques.

\subsection{Coleman's explicit Chabauty}

In a seminal paper of the 1980's, 
Coleman revisited Chabauty's method and realized one could carve out an explicit upper bound. 
\begin{theorem}[\cite{Coleman:effectiveChabauty}]
\label{T:coleman-effective-chabauty}
 Let $X$ be a curve of genus $g$ and let $r = \rank_{\Z} J(\Q)$.
 Suppose $p > 2g$ is a prime of {\bf good reduction}. Suppose $r <
 g$. Then
 $$\#X(\Q) \leq \#X(\F_{p}) + 2g - 2.$$
\end{theorem}

A modified statement holds for $p \leq 2g$ or for $F \neq \Q$. Also note: {\bf this does not prove uniformity} (since the first prime $p$ of good reduction, and hence $\#X(\F_{p})$, might be large).

Here, as in Chabauty's original proof, one makes use of the inclusion of $X(\Q)$ into the more tractable stand-in $X(\Q_p)\cap \overline{J(\Q)}$.  This is a general theme of extensions that go under the name {\em Chabauty method} where one bounds $X(\Q)$ by bounding an enlargement usually arising as a subset of $X(\Q_p)$. 

Coleman's major insight was that the zeros of the integrals $f_{\omega}$ are amenable to a rather explicit $p$-adic local analysis (via Newton polygons). The analysis in Coleman's proof of Theorem \ref{T:coleman-effective-chabauty} roughly breaks into the following steps. 
 \begin{enumerate}
 \item \textbf{Local Bounds (``$p$-adic Rolle's'')}:
the number of zeros of $f_\omega$ in a tube $D_P$ is at most $1 + n_P$, where $n_P = \#\left(\div \omega \, \cap D_P\right)$. 
\vspace{2pt}
 \item \textbf{Global geometric step}: by 
 degree considerations, $\sum n_P = 2g-2$.
\vspace{2pt}
 \item \textbf{Coleman's bound}: $\sum_{P \in X(\F_p)} (1 + n_P) = \#X(\F_p) + 2g-2$.
 \end{enumerate}
 
 We will discuss the local bounds below, but first we say some words about (2).  Let $\cX$ be a smooth proper model of $X$ over $\Z_p$.
 The sheaf of regular $1$-forms has a canonical extension $\Omega^1_{\cX/\Z_p}$, and any rational differential $\omega \in H^0(X_{\Q_p},\Omega^1)$ extends (after multiplication by an appropriate power of $p$) to a regular nonzero section of $\Omega^1_{\cX/\Z_p}$ whose restriction $\omega|_{\cX_{\F_p}}$ to the special fiber $\cX_{\F_p}$ is a regular non-zero $1$-form.  As we discuss below, we have
 \[\#\div(\omega|_{\cX_{\F_p}})=\sum_{P\in \cX_{\F_p}(\overline{\F}_p)}  \#\left(\div \omega \, \cap D_P\right).\]
 From this, we obtain (2).  Coleman's bound (3), of course, follows directly from (1) and (2).
 
 \subsection{Local bounds}
\label{ss:local-bounds}

 The local bound is an exercise in Newton polygons.  As above, for a point $Q$ of $X$ specializing to $P\in \cX(\overline{\F}_p)$, a choice $t$ of uniformizer at $Q$ induces a $p$-adic analytic isomorphism from the tube $D_P$ to the $p$-adic disc $p \overline{\Z}_p = \{ x : x \in \overline{\Q}_p \mid |x|<1\}$.  Note also that all elements of $X(\Q_p)$ in $D_P$ lie within the smaller disc $D'=\{ x : x \in \overline{\Q}_p \mid |x| \leq 1/p\}$.
 
 The $1$-form $\omega$ can be expanded near $Q$ as
 \[\omega=\sum a_n t^n dt.\]
Rescaling, we may suppose that $\min(\val(a_n))=0$.  
Let $n_0(\omega)$ be the smallest value of $n$ such that $\val(a_n)=0$.  By a standard Newton polygon argument, $\omega$ has $n_0(\omega)$ zeros in $D_P$.  Therefore, by the degree computation $n_0(\omega)\leq 2g-2$.  The reduction of $\omega$ mod $p$ has the form $t^{n_0}g(t)$ where $g(0)\neq 0$.  Consequently, $\omega|_{\cX_{\F_p}}$ has a zero of order $n_0$ at $P$.
 
 Now, integrating $\omega$ term-by-term gives
 \[f_\omega=\sum \frac{a_n}{n+1} t^{n+1} dt+c\]
 for some constant $c$.  To get some control over $c$, we will concede one zero at $Q$ and suppose that $f_\omega(Q)=c=0$.  If $f_\omega$ has no zeros in $D_P$, then our bound is still true.  Now, the coefficients of $f_\omega$ satisfy $\val\left(\frac{a_n}{n+1}\right)=\val(a_n)-\val(n+1)$, and thus 
 the Newton polygon of $f_\omega$ is very nearly the Newton polygon of $\omega$ shifted one unit to the right, so we would expect $f_\omega$ to have at most $n_0(\omega)+1$ zeros in $D'$.  If $p>n_0(\omega)-2$, this is indeed the case (see \cite[Lemma 5.1]{McCallumP:chabautySurvey}), and otherwise there is a small error term coming from the appearance of powers of $p$ in the denominators of the coefficients of $f_\omega$ (see \cite[Section 6]{Stoll:2006-independence-of-rational-points} for a precise error term).  

 \subsection{The case of bad reduction}

One can extend the above to the case of bad reduction to prove the following theorem \cite{LoreniniT:thue}.
\begin{theorem} Let $X$ be a curve of genus $g$ and let $r = \rank_{\Z} J(\Q)$.
Let $p>2g$ be a prime.  Let $\cX$ be a proper regular model for $X$ over $\Z_p$.  Then 
 \[
 \#X(\Q) \leq \#\cX_{\F_p}^{\sm}(\F_{p}) + 2g - 2,
 \]
where $\cX_{\F_p}^{\sm}$ denotes the smooth locus of the special fiber $\cX_{\F_p}$. 
 \end{theorem}

We will outline the argument given in \cite{McCallumP:chabautySurvey}.  By standard arguments involving regular models, every rational point of $X$ specializes to a smooth $\F_p$-point of $\cX_{\F_p}$ so we only have to consider tubes around such points.
 
 The sheaf of regular $1$-forms has an extension to the invertible canonical sheaf $\omega_{\cX/\Z_p}$ to which $\omega$ extends as a rational section.  
Now, given any reduced component $C$ of $\cX_{\F_p}$, there is an integer $m(C)$ such that $p^{m(C)}\omega$ restricts to a non-zero rational differential along $C$.  For $P$ a smooth point of $\cX_{\F_p}$ lying on $C$, $n_P$ is equal to the order of vanishing of $(p^{m(C)})\omega|_C$ at $P$.  By a computation in the intersection theory of arithmetic surfaces (see \cite[Appendix A]{McCallumP:chabautySurvey}), one can show
 \[\sum_{P\in \cX^{\sm}_0(\F_p)} n_P \leq 2g-2.\]
 The proof of the local bound proceeds as before.  Note that we only need to compute $p$-adic integrals between points in the same tubes.  Such a $p$-adic integral is just evaluating a primitive of $\omega$ and there are no difficulties with its definition on curves of bad reduction.
 
\subsection{Stoll bounds} 

 The Chabauty bounds can be tightened further if $r\leq g-2$ using an argument of Stoll \cite{Stoll:2006-independence-of-rational-points}; this was extended to the bad reduction case by the first and third authors \cite{katzZB:tropicalChabauty}.
 \begin{theorem} 
Let $X$ be a curve of genus $g$ and let $r = \rank_{\Z} J(\Q)$.
Let $p>2g$ be a prime.  Let $\cX$ be a proper regular model for $X$ over $\Z_p$.  Then
 \[
 \#X(\Q) \leq \#\cX_{\F_p}^{\sm}(\F_{p}) + 2r.
 \]
 \end{theorem}
 
 The proof involves optimizing the choice of $1$-form on each tube.  Recall that the space $V$ of differentials such that 
 \begin{equation}
   \label{eq:chabautys-theorem}
\int_Q^{Q'} \omega = 0 \text{ for } Q, Q' \in X(\Q) \text{ and } \omega \in V
\end{equation}
has dimension at least $g-r$. For each $P \in X(\F_p)$, set
\[
\widetilde{n_P} := \min_{\omega \in V \setminus \{0\}} \#\left(\div \omega \, \cap D_P\right).
\]
Then Coleman's local analysis still gives that $\#\left(D_P \cap X(\Q)\right) \leq 1 + \widetilde{n_P}$. Stoll proved that 
\begin{equation}
\label{eq:stoll}
\sum \widetilde{n_P} \leq 2r,
\end{equation}
whence 
\[
\#X(\Q) \leq \sum_{P \in X(\F_p)} (1 + \widetilde{n_P}) \leq \#X(\F_p) + 2r.
\]
Stoll's proof of the inequality \ref{eq:stoll} exploits another ``global geometric tool" from the theory of algebraic curves: Clifford's theorem. Consider the divisor $D = \sum \widetilde{n_P} P$. Then, essentially by construction,
$V \subset H^0(X_{\F_p}, \Omega^1(-D))$. Together with Clifford's theorem and Chabauty's theorem this gives
\[
g-r \leq \dim H^0(X_{\F_p}, \Omega^1(-D)) \leq \frac{1}{2}\deg \Omega^1(-D) + 1,
\]
and simplifying gives Stoll's inequality.

In the bad reduction case, one first reduces to the case of semistable models and then uses a \emph{combinatorial} version of Clifford's theorem adapted to that setting.  Such a result is phrased in terms of abelian rank functions in \cite{katzZB:tropicalChabauty} which extend and make use of the theory of linear systems on graphs due to Baker and Norine \cite{BakerN:RR}.  Such arguments were systematized into the theory of linear systems on metrized complexes of curves in the work of Amini and Baker  \cite{Amini2012}.

\subsection{Uniformity for hyperelliptic curves} 
\label{ss:stoll-discussion}
A recent breakthrough of M. Stoll \cite{stoll:uniform} removed (for hyperelliptic curves) the dependence from Chabauty--Coleman on a \emph{regular} model, and derived (for $r \leq g-3$) a \emph{uniform} bound on $\#X(\Q)$.
\begin{theorem}[Stoll, \cite{stoll:uniform}]
\label{T:stoll-uniform-hyperelliptic-2}
 Let $X$ be a smooth hyperelliptic curve of genus $g$ and let $r = \rank_{\Z} J(\Q)$.
 Suppose $r\, \leq \, g-3$. Then 
 $$\#X(\Q) \leq 8(r+4)(g-1) + \max\{1,4r\}\cdot g.$$
\end{theorem}

 The key elements of Stoll's proof are:
\begin{enumerate}
\item analysis of $p$-adic integration on \emph{annuli}, rather than tubes, 
\item a comparison of different analytic continuations of $p$-adic integrals, and
\item a $p$-adic Rolle's theorem for hyperelliptic integrals.
\end{enumerate}

To elaborate on (1): the utility of a \emph{regular} model $\cX$ is that the reduction of a rational point $P \in \cX(\Q)$ is a smooth point of $P' \in \cX(\F_p)$, and the tube of $P'$ is a disk, which greatly facilitates the analysis in any Chabauty-type argument. The limitation is that $\#\cX^{\sm}(\F_p)$ (and consequently any upper bound on $\#\cX(\Q)$ arising from the analysis described above) can be arbitrarily large.

Stoll's insight was to forgo regularity and (morally) work instead with a \emph{stable} model $\cX$. The advantage is that it is quite straightforward to bound $\#\cX(\F_p)$ in terms of $p$ and $g$. The cost is that the tube over a node is an annulus, instead of a disk, and one must execute the Chabauty--Coleman analysis over annuli.

This creates several significant technical difficulties and leads to a discussion of  (2). So far this survey has only discussed integration between points in the same residue disk. When $\cX$ is smooth, Coleman (surprisingly) proved that there is a unique way to extend  the ``naive'' integrals to all of $\cX$ (i.e., a unique way to integrate between points in different residue disks). In general, though, there are multiple ways to ``analytically continue'' these integrals and various issues that arise. 

First, restricting a differential $\omega$ to an annulus gives a Laurent
expansion. One can no longer naively integrate the $dt/t$ term, and in fact
making sense of $\int dt/t$ more or less amounts to picking a branch of the $p$-adic logarithm, of which there are many. But beyond this, there are still multiple natural ways to extend integration across disks: the \emph{abelian} integral  $\Abint$ (arising from punting the problem to the Lie algebra of the Jacobian), and a rather different integral $\BCint$ due to Berkovich and Coleman. Each has relative merits:
\begin{enumerate}
\item $\Abint$ is the integral used in Chabauty's methods (i.e., so that Equation \ref{eq:chabautys-theorem} holds), but
\item $\BCint$ can be computed via primitives (as in Section \ref{ss:integration-warmup}, (3)) on \emph{annuli}  (while $\Abint$ cannot), which is necessary for any type of analysis leading to an explicit upper bound.
\end{enumerate}
See Section \ref{sec:integration} for a much longer discussion of these technicalities. 

This tension is resolved via the stronger hypothesis $r\leq g-3$. Let $V$ be as in Equation \ref{eq:chabautys-theorem}. Stoll proved that on an annulus, $\Abint$ and $\BCint$ differ by a linear form in $\omega$, and thus, if $\dim V \geq 2$, then there is some non-zero $\omega \in V$ such that $\Abint \omega = \BCint \omega$. If, additionally, $\dim V \geq 3$, then we may also choose $\omega$ to have trivial residue (i.e., the $dt/t$ term of $\omega$ is zero). Under these hypotheses, one can thus work directly with the Laurent series expansion of $\omega$ on any given annulus and attempt a Newton polygon analysis as in Subsection \ref{ss:local-bounds}.

Finally, the hyperelliptic restriction of Stoll's theorem arises from his proof of (3). A crucial step is to prove a ``$p$-adic Rolle's theorem'' for integrals on annuli, i.e., to relate the zeros of an integral $\int \omega$ to the zeros of $\omega$ itself (as in Subsection \ref{ss:local-bounds} above), and for Stoll's very direct approach to $p$-adic Rolle's it is essential that differentials on a hyperelliptic curve have an explicit description $\omega = f(x)dx/y$. 

\subsection{Uniformity} 
\label{ss:KRZB-prelude}

This discussion culminates with the following theorem, which uses the full
catalogue of tropical and non-Archimedean analytic tools to overcome the hyperelliptic restriction and prove uniform bounds for \emph{arbitrary} curves of small rank.

\begin{theorem}[Katz--Rabinoff--Zureick-Brown, \cite{KatzRZB-uniform-bounds}]
\label{T:KRZB-2}
 Let $X$ be {\textbf{any}} smooth curve of genus $g$ and let $r = \rank_{\Z} J(\Q)$.
 Suppose $r\, \leq \, g-3$. Then
 $$\#X(\Q) \leq 84g^2 - 98g + 28.$$
\end{theorem}

\begin{remark}
A few aspects of the proof deserve immediate elaboration.

\begin{enumerate}
\item Tropical and combinatorial tools afford a proof of a very general $p$-adic Rolle's theorem, free of the necessity of explicit equations. This is the main additional input needed to promote Stoll's result from hyperelliptic to general curves.
\item Chabauty's method (of bounding $X(\Q)$ via $p$-adic integrals) is systematically redeveloped using Berkovich spaces instead of rigid spaces.
\item Consequently, Stoll's \emph{comparison of various analytic continuations}
  of the integrals is greatly clarified by the use of Berkovich machinery -- the
  uniformization (i.e., the universal cover, in the sense of algebraic topology)
  of the Berkovich curve plays an essential role, as does Baker--Rabinoff's
  \cite{bakerR:skeleton} comparison of the tropical and algebraic Abel--Jacobi maps.
\item Traditionally, one restricts an integral $\int \omega$ to a covering of the curve by balls, where the integral becomes a power series and can be analyzed directly. Stoll's achievement was to perform this analysis on annuli, where instead the formal object is a Laurent series. 

  The techniques developed in \cite{KatzRZB-uniform-bounds} allow one to analyze these integrals on much larger domains, namely, basic wide opens (in the terminology of Coleman); this is the key input to the uniform Manin--Mumford results discussed in Subsection \ref{ss:effective-manin-mumford}.
\end{enumerate}

\end{remark}

Most of the rest of this survey is dedicated to expositing these tools
and their application to Theorem \ref{T:KRZB-2}.

\section{Berkovich curves and skeletons}
\label{sec:berk-curv-skel}

In this section we explain the basics of the structure theory of Berkovich
curves and their skeletons, with a brief foray into potential theory.  This
non-Archimedean analytic language provides a robust and convenient framework for
producing the combinatorial data of tropical meromorphic functions and divisors
on metric graphs, from the algebro-geometric data of sections of line bundles on
algebraic curves.  Our treatment of Berkovich curves is quite threadbare; for
more detail, see the book of Baker--Rumely~\cite{Baker2010}, the paper of
Baker--Payne--Rabinoff~\cite{Baker2013}, and Berkovich's
book~\cite{Berkovich:spectralTheory}.

Let $K$ be an algebraically closed field which is complete with respect to a
nontrivial, non-Archimedean valuation $v\colon K\to\R\cup\{\infty\}$.  We will
take $K$ to be $\C_p$ in the sequel; the hypothesis that $K$ be algebraically
closed is mainly for simplicity of exposition, and because all of our arguments
involving Berkovich spaces are geometric in nature.  Let $R$ be the valuation
ring of $K$, let $\fm_R$ be its maximal ideal, and let $k = R/\fm_R$ be its
residue field.  Let $|\scdot| = \exp(-v(\scdot))\colon K\to\R_{\geq 0}$ be an
absolute value associated to $v$.

\subsection{Analytification}\label{sec:analytification}
To a variety (or locally finite-type scheme) $X$ over $K$, one can functorially
associate a non-Archimedean analytification $X^\an$, in the sense of Berkovich.
The theory of $K$-analytic spaces is meant to mirror the classical theory of
complex-analytic spaces as closely as possible.  The topological space
underlying $X^\an$ is locally compact and locally contractible, and reflects
many of the algebro-geometric properties of $X$: for instance, $X$ is separated
if and only if $X^\an$ is Hausdorff; $X$ is proper if and only if $X^\an$ is
compact and Hausdorff; $X$ is connected if and only if $X^\an$ is
path-connected.  (These facts hold true even if the valuation on $K$ is trivial,
an observation which provides curious equivalent valuation-theoretic definitions
of these purely scheme-theoretic properties.)  By functoriality, the
analytification contains $X(K)$ in a natural way, but unlike for $\C$-analytic
spaces, the set of $K$-points of $X$ should be regarded as a relatively
unimportant subset of $X^\an$ which lies ``at infinity''.  Indeed, the subspace
topology on $X(K)\subset X^\an$ coincides with the natural totally disconnected
topology on $X(K)$, defined by local embeddings into $K^n$, yet $X^\an$ is
path-connected when $X$ is connected!

The topological space $X^\an$ underlying an affine variety $X = \Spec(A)$ is
defined to be the space of all multiplicative seminorms
$\|\scdot\|\colon A\to\R_{\geq 0}$ that restrict to $|\scdot|$ on $K$, equipped
with the topology of pointwise convergence.  By \defi{multiplicative seminorm}
we mean a homomorphism of multiplicative monoids which satisfies the strong
triangle inequality.  For a general variety $X$, the topological space $X^\an$
is constructed affine-locally by gluing.

\begin{egsub}\label{eg:evaluation.norm}
  Let $X = \Spec(A)$ be an affine $K$-variety and let $P\in X(K)$.  Then
  $\|f\|_P \coloneq |f(P)|$ defines a multiplicative seminorm.  This is the point
  of $X^\an$ corresponding to $P$.  Note that the kernel of $\|\scdot\|_P$, i.e.,
  the set $\{f\mid \|f\|_P=0\}$, is the maximal ideal of $A$ corresponding to $P$.
\end{egsub}

\begin{egsub}\label{eg:supremum.norm}
  Let $r\in\R$ and let $\rho = \exp(-r)$.  The rule
  \[ \bigg\|\sum a_it^i\bigg\|_r\coloneq\max\big\{|a_i|\rho^i\big\}
  \sptxt{i.e.}
  -\log\bigg\|\sum a_it^i\bigg\|_r = \min\big\{v(a_i) + ir\big\} \]
  defines a multiplicative seminorm $\|\scdot\|_r\colon K[t]\to\R_{\geq0}$,
  hence a point of $\bA^{1,\an}$.  The proof that $\|\scdot\|_r$ is
  multiplicative is essentially Gauss' lemma (on the content of a product of
  polynomials); for this reason, these seminorms are often called Gauss points
  (at least when $r\in v(K^\times)$).
  Note that the kernel of $\|\scdot\|_r$ is trivial, so that
  $\|\scdot\|_r$ is in fact a norm on $K[t]$.  It turns out that
  \[ \|f\|_r = \sup\big\{|f(x)|\mid x\in K, |x|\leq\rho\big\}, \]
  which is the supremum norm over the $K$-ball of radius $\rho$
  around $0$.  Of course, composing $\|\scdot\|_r$ with any translation
  $t\mapsto t-x$ yields another multiplicative seminorm on $K[t]$, namely, the
  supremum norm over the ball of radius $\rho$ around $x\in K$.  The limit as
  $r\to\infty$ of $\|f\|_r$ is the seminorm $f\mapsto|f(0)|$ of
  Example~\ref{eg:evaluation.norm}; this explains the heuristic that the
  $K$-points of $\bA^{1,\an}$ lie ``at infinity''.

\end{egsub}

Berkovich classified points of $\bA^{1,\an}$ into four types:
\begin{itemize}
  \item Type~1: These are the points $f\mapsto|f(x)|$ for $x\in K$.
  \item Type~2: These are the points $f\mapsto\|f(t-x)\|_r$ for $x\in K$ and
    $r\in v(K^\times)$, the value group of $K$.
  \item Type~3: These are the points $f\mapsto\|f(t-x)\|_r$ for $x\in K$ and
    $r\notin v(K^\times)$.
  \item Type~4: These mysterious points are constructed from infinite descending
    sequences of closed balls in $K$ with empty intersection.  (Fields for which
    no such sequence exists are called \defi{spherically complete}.  An example
    is the field of Hahn series $\C\ls{t^\R}$.  The field $\C_p$ is very far from
    being spherically complete although there is a spherical complete field
    $\Omega_p$ containing it.)  These points are necessary for local compactness
    of $\bA^{1,\an}$, but will not play any role in this paper.
\end{itemize}
In particular, nearly all%
\footnote{In terms of categorization, not cardinality} of the points of
$\bA^{1,\an}$ are described in Example~\ref{eg:supremum.norm}.
See~\cite[\S1.2]{Baker2010}.

\begin{notnsub}
  Let $X$ a be $K$-variety, let $x\in X^\an$, and let $f$ be a regular function
  defined on an affine neighborhood $U = \Spec(A)\subset X$ containing $x$ in its
  analytification.  We write
  \[ |f(x)| \coloneq \|f\|_x, \]
  where $\|\scdot\|_x$ is the multiplicative seminorm on $A$ which corresponds
  to $x$.
\end{notnsub}

So far we have only discussed the topological space underlying  the
analytification of a variety; we have not discussed analytic functions.  We will
do so briefly through some examples.

\begin{egsub}\label{eg:ball}
  Fix a parameter $t$ on $\bA^1$, so $\bA^1 = \Spec(K[t])$, and let $\rho>0$.
  The \defi{closed ball of radius $\rho$} around zero is
  \[ \bB(\rho) \coloneq \big\{x\in\bA^{1,\an}\mid |t(x)| \leq \rho\big\}. \]
  The analytic functions on $\bB(\rho)$ consist of those power series
  $f = \sum a_it^i\in K\ps t$ such that $|a_i|\rho^i\to 0$ as $i\to\infty$.
  Such a power series $f$, called \textit{strictly convergent}, is a uniform
  limit on $\bB(\rho)$ of polynomials.  In particular, it makes sense to
  evaluate $f$ on points of $\bB(\rho)$.  This is true in general of analytic
  functions. 

  The \defi{open ball of radius $\rho$} around zero is
  \[ \bB(\rho)_+ \coloneq \big\{x\in\bA^{1,\an}\mid |t(x)| < \rho\big\}. \]
  An analytic function on $\bB(\rho)_+$ is a power series $f = \sum a_it^i$
  which is analytic on $\bB(\tau)$ for all $\tau<\rho$, i.e., such that
  $|a_i|\tau^i\to 0$ as $i\to\infty$ for all $\tau<\rho$.

  See Figure~\ref{fig:ball} for an illustration, and~\cite[Chapter~1]{Baker2010},
  \cite[\S2]{Baker2013}, and~\cite[Chapter~4]{Berkovich:spectralTheory} for details.

  \begin{figure}[ht]
    \centering
    \includegraphics{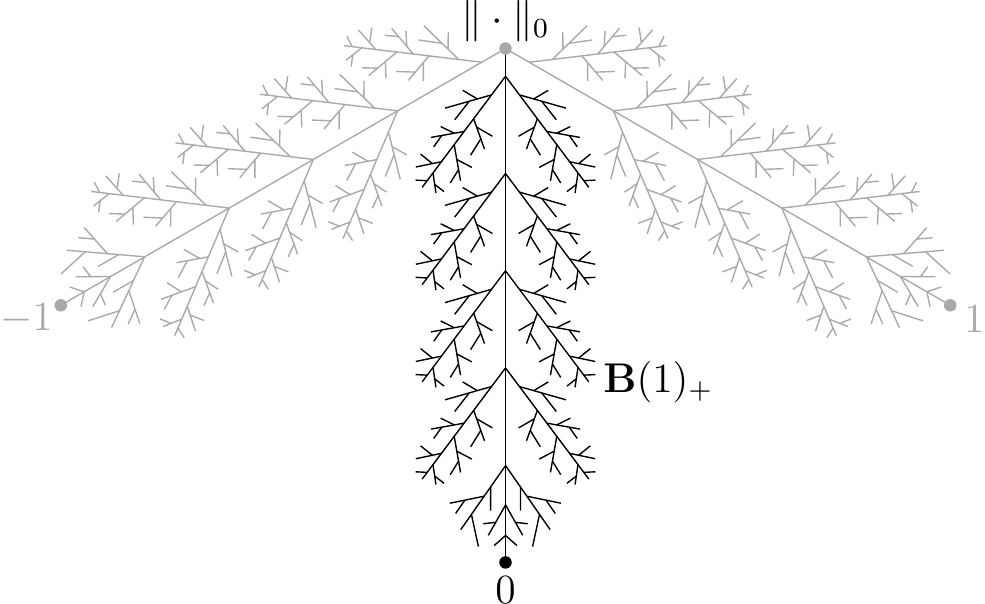}
    
    \caption{The closed ball $\bB(1)$, and its open sub-ball $\bB(1)_+$.  The
      trees depicted in grey represent open balls around points with distinct
      images in $k = R/\fm_R$ (so $\chr(k)\neq 2$ in the picture);
      there are infinitely many of these.
      Each open ball is an $\R$-tree: it is a metric tree with infinite
      ramification along infinitely many points along each edge.  The set of points
      of types~1--3 in $\bB(1)$ consist of the lines from the Gauss point
      $\|\scdot\|_0$ to the classical points $x\in R$, i.e.\ the lines
      \[ \big\{f(t)\mapsto\|f(t-x)\|_r\mid t\in\R_{\geq0}\cup\{\infty\}
      \big\} \] for $x\in R$.  (The type-4 points are also leaves.)}
    \label{fig:ball}
  \end{figure}

\end{egsub}

\begin{egsub}\label{eg:annulus}
  Let $\rho\in(0,1)$.  The \defi{closed annulus of modulus $\rho$} is
  \[ \bS(\rho) \coloneq \big\{x\in\bA^{1,\an}\mid|t(x)|\in[\rho,1]\big\}. \]
  The analytic functions on $\bS(\rho)$ consist of those infinite-tailed Laurent
  series $f = \sum_{i\in\Z}a_i t^i\in K\ps{t,t\inv}$ such that 
  $|a_i|\to 0$ as $i\to+\infty$ and $|a_i|\rho^i\to 0$ as $i\to-\infty$.  The
  set of such Laurent series forms a ring (although $K\ps{t,t\inv}$ does not).
  As in Example~\ref{eg:ball}, it makes sense to evaluate such Laurent series on
  points of $\bS(\rho)$.

  The \defi{open annulus of modulus $\rho$} is
  \[ \bS(\rho)_+ \coloneq \big\{x\in\bA^{1,\an}\mid|t(x)|\in(\rho,1)\big\}. \]
  An analytic function on $\bS(\rho)_+$ is an infinite-tailed Laurent series
  $f = \sum_{i\in\Z} a_it^i$ such that $|a_i|\tau^i\to 0$ as $i\to\pm\infty$ for
  all $\tau\in(\rho,1)$.

  The modulus $\rho$ is an isomorphism invariant of an open or closed annulus,
  as in complex analysis.

  See Figure~\ref{fig:annulus_skeleton} for an illustration,
  and~\cite[\S2]{Baker2013} for details.
\end{egsub}

\subsubsection{The skeleton of an annulus}\label{sec:skeleton-an-annulus}
Let $r=-\log(\rho) > 0$.  The mapping $\sigma\colon(0,r)\to\bS(\rho)_+$ defined
by $\sigma(s)=\|\scdot\|_s$ is a continuous embedding; its image is called the
\defi{skeleton} of $\bS(\rho)_+$ and is denoted $\Sigma(\bS(\rho)_+)$.  The
mapping $\trop\colon\bS(\rho)_+\to(0,r)$ defined by $\trop(x) = -\log|t(x)|$ is
a left inverse to $s$.  Identifying $(0,r)$ with its image
$\Sigma(\bS(\rho)_+)$, Berkovich
showed~\cite[Proposition~4.1.6]{Berkovich:spectralTheory} that
$\trop\colon\bS(\rho)_+\to\Sigma(\bS(\rho)_+)$ is the image of a strong
deformation retraction.  Moreover~\cite[Proposition~2.2.5]{Thuillier2005}, the
skeleton $\Sigma(\bS(\rho)_+)$ consists of all points of $\bS(\rho)_+$ that do
not admit a neighborhood isomorphic to $\bB(1)_+$, so that the skeleton is
well-defined as a set.  In fact it is canonically identified with the interval
$(0,r)$, up to the flip $x\mapsto r-x$; this is due to the fact that the modulus
$\rho$ is an isomorphism invariant of $\bS(\rho)_+$.

\begin{figure}[ht]
  \centering
  \includegraphics{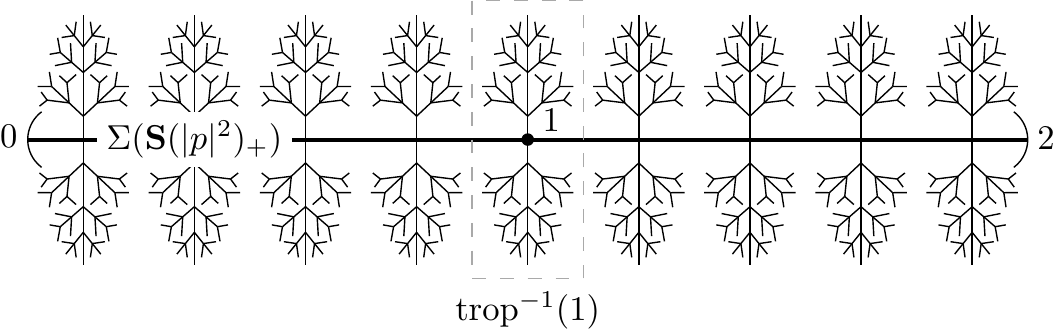}
  
  \caption{The open annulus $\bS(|p|^2)_+$ and its skeleton.  Here the base
    field is $\C_p$, with $v(p) = 1$.  The skeleton in the center is identified
    with the open interval $(0,2)$.  The rest of $\bS(|p|^2)$ is obtained by
    attaching a copy of $\bB(1)\setminus\bB(1)_+$ (the closed annulus of modulus
    $1$) to each $\Q$-point of $(0,2)$.  The retraction map takes each tree
    to its root.}
  \label{fig:annulus_skeleton}
\end{figure}

See Figure~\ref{fig:annulus_skeleton} and also~\cite[\S2]{Baker2013}, 

\subsection{Reduction}\label{sec:reduction}

Let $\cX$ be a proper and flat scheme over the valuation ring $R$ of $K$, and
let $X = \cX_K$.  Let $x\in X^\an$ be any point.  The kernel of an associated
seminorm $\|\scdot\|_x$ on some affine neighborhood in $X$ is a prime ideal,
hence determines a point $y\in X$, called the \defi{center} of $\|\scdot\|_x$.
The seminorm $\|\scdot\|_x$ extends to the local ring $\sO_{X,y}$, and defines a
norm on the residue field $\kappa(y) = \sO_{X,y}/\fm_{X,y}$, which by definition
extends the absolute value $|\scdot|$ on $K$.  We therefore obtain a morphism
$\Spec(\kappa(y))\to X$ over a morphism $\Spec(\kappa(y)^\circ)\to\Spec(R)$,
where $\kappa(y)^\circ$ is the valuation ring with respect to $\|\scdot\|_x$.
By the valuative criterion of properness, there exists a unique morphism
$\Spec(\kappa(y)^\circ)\to\cX$ sending the closed point of
$\Spec(\kappa(y)^\circ)$ to the special fiber $\cX_k$ of $\cX$.  Define
$\red(x)\in\cX_k$ to be the image of the closed point.  The resulting map
$\red\colon X^\an\to\cX_k$ is called the \defi{reduction} or
\defi{specialization} map.  It is surjective, and anti-continuous in the sense
that the inverse image of an open set is closed (and vice-versa).
See~\cite[\S.4]{Berkovich:spectralTheory}.

\begin{egsub}\label{eg:reduction.closed}
  Let $\cX = \bP^1_R$, and let $x\in R$.  This defines a point of $\bP^{1,\an}$
  given on $\bA^1 = \Spec(K[t])\subset\bP^1_K$ by the seminorm
  $f\mapsto|f(x)|\colon K[t]\to\R_{\geq 0}$.  The kernel of this seminorm is the
  prime ideal corresponding to $x$, so we have $y=x$ in the above notation, and
  the norm on $\kappa(x) = K$ defined by $\|\scdot\|_x$ is the usual absolute
  value $|\scdot|$.  As $x\in R$, the map $f\mapsto f(x)\colon K[t]\to K$
  restricts to a map $R[t]\to R$; its special fiber is the map $k[t]\to k$
  corresponding to the reduction $x\mod\fm_R\in k$.  So in this case, the
  reduction map $\red$ coincides with the quotient map $R\to k$.
\end{egsub}

\begin{egsub}\label{eg:reduction.gauss}
  Continuing with Example~\ref{eg:reduction.closed}, let $r\geq 0$ and let
  $x = \|\scdot\|_r\in\bA^{1,\an}\subset\bP^{1,\an}$ be the seminorm defined in
  Example~\ref{eg:supremum.norm}.  This norm has trivial kernel, so its center
  $y$ is the generic point of $\bP^1$.  In other words, $\|\scdot\|_r$ extends
  to a norm $\|\scdot\|\colon K(t)\to\R_{\geq 0}$ by
  $\|f/g\|_r = \|f\|_r/\|g\|_r$.  Let $\rho = \exp(-r)$, and recall that for
  $f = \sum a_it^i\in K[t]$, we have $\|f\|_r = \max\{|a_i|\rho^i\}$.  As
  $r\geq 0$, we have $\rho\leq 1$, so the polynomial ring $R[t]$ is contained in
  the valuation ring $K(t)^\circ\subset K(t)$ with respect to $\|\scdot\|_r$.
  This inclusion $R[t]\subset K(t)^\circ$ defines the extension morphism
  $\Spec(K(t)^\circ)\to\bA^1_R\subset\bP^1_R$.  A polynomial $f\in R[t]$ maps to
  the maximal ideal $K(t)^{\circ\circ}\subset K(t)^\circ$ if and only if
  $|a_i|\rho^i < 1$ for all $i\geq 0$.  When $r=0$ and $\rho=1$, this is the
  case if and only if $f\in\fm_R R[t]$, so the retraction of $K(t)^{\circ\circ}$
  is $\fm_R R[t]$, and thus $\red(x)$ is the generic point of $\bP^1_k$.  When
  $r > 0$ and $\rho < 1$, we have $|a_i|\rho^i<1$ for all $i > 0$, so the
  retraction of $K(t)^{\circ\circ}$ is the ideal $(\fm_R,t)\subset R[t]$, and
  therefore $\red(x)$ is the closed point $0\in\bA^1(k)$.

  In this way it is not hard to see that for a closed point
  $\td x\in \bA^1(k) = k$, the inverse image $\red\inv(\td x)$ is the open ball
  of radius $1$ around any point $x\in R$ reducing to $\td x$, and the inverse
  image of $\infty\in\bP^1(k)$ is the ``open ball of radius $1$ around
  $\infty\in\bP^1(K)$'', i.e., the image of $\bB(1)_+$ under the involution
  $t\mapsto t\inv\colon\bP^1\to\bP^1$.  The inverse image of the generic point
  of $\bP^1_k$ is the singleton set $\{\|\scdot\|_0\}$.
\end{egsub}

\begin{egsub}\label{eg:reduction.annulus}
  Fix $\varpi\in K$ with $0<|\varpi|<1$, and let $\cX\subset\bP^2_R$ be the
  variety defined by the homogeneous equation $uv = \varpi w^2$.  Let
  \[ \cU = \Spec(R[u,v]/(uv-\varpi))\subset \cX \]
  be the affine open defined by the nonvanishing of $w$.  As usual we set
  $X = \cX_K$ and $U = \cU_K$.  Note that $(u,v)\mapsto u$ defines an
  isomorphism $U\isom\bG_{m,K}$ with inverse $t\mapsto(t,\varpi/t)$; we use these
  isomorphisms to identify $U$ with an open subscheme of $\bA^1_K$.

  Let $\td x\in\cU(k)\subset\cX_k(k)$ be the point with $u=v=0=\varpi\mod\fm_R$.
  This is the point of $\cU$ defined by the ideal
  $(u,v,\fm_R)\subset R[u,v]/(uv-\varpi)$.  Choose $x\in U^\an$ such that
  $|u(x)|\leq 1$ and $|v(x)|\leq 1$.  Such $x$ correspond to the points of the
  closed annulus $\bS(|\varpi|)\subset\bA^{1,\an}$.  Let $y\in U$ be the center
  of $x$.  The image of $R[u,v]/(uv-\varpi)$ is contained in the valuation ring
  $\kappa(y)^\circ$ of $\kappa(y)$ with respect to $\|\scdot\|_x$, and it is
  clear that the maximal ideal $\kappa(y)^{\circ\circ}$ contracts to
  $(u,v,\fm_R)$ if and only if $|u(x)| < 1$ and $|v(x)| < 1$.  This is the case
  if and only if $x\in\bS(|\varpi|)_+\subset\bS(|\varpi|)$, so in summary,
  $\red\inv(\td x)\cong\bS(|\varpi|)_+$.  Of note here is the fact that the
  ``thickness'' of the singularity $\td x\in\td X$ is measured by the
  \emph{modulus} of the annulus $\red\inv(\td x)$.
\end{egsub}

Berkovich~\cite[Proposition~2.4.4]{Berkovich:spectralTheory} and
Bosch--L\"utkebohmert~\cite[Propositions~2.2 and~2.3]{BoschL:stableReductionAnd}
showed that Examples~\ref{eg:reduction.gauss} and~\ref{eg:reduction.annulus} are
representative of fairly general phenomena, in the following way.  By a
\defi{semistable $R$-curve} we mean a proper, connected, flat $R$-scheme $\cX$
of relative dimension one, with smooth generic fiber $X$ and reduced special
fiber with nodal singularities.  We also call $\cX$ a \defi{semistable model} of
$X$.

\begin{thmsub}[Berkovich; Bosch--L\"utkebohmert]\label{thm:formal.fibers}
  Let $\cX$ be a semistable $R$-curve, let $X = \cX_K$, and let $\td x\in\cX_k$.
  Then:
  \begin{enumerate}
  \item $\td x$ is a generic point if and only if $\red\inv(\td x)$ is a
    singleton set;
  \item $\td x$ is a smooth point if and only if $\red\inv(\td x)\cong\bB(1)_+$; and
  \item $\td x$ is a nodal point if and only if
    $\red\inv(\td x)\cong\bS(\rho)_+$ for some $\rho\in|K^\times|$.
  \end{enumerate}
\end{thmsub}

Case~(3) can be made more precise: if $\td x\in\cX_k$ is a node, then there is
an \'etale neighborhood of $\td x$ which is an \'etale cover of
$\Spec(R[u,v]/(uv-\varpi))$ for some $\varpi\in K^\times$ with $|\varpi|<1$; then
$\rho = |\varpi|$.

With the semistable reduction theorem, Theorem~\ref{thm:formal.fibers} easily
implies the following important structure theorem for analytic curves.

\begin{corsub}\label{cor:semistable.decomposition}
  Let $X$ be a smooth, proper, connected $K$-curve.  Then there exists a finite
  set $V\subset X^\an\setminus X(K)$ such that $X^\an\setminus V$ is isomorphic
  to a disjoint union of infinitely many open balls and finitely many open
  annuli.
\end{corsub}

Indeed, if $\cX$ is a semistable model for $X$, one can take the set $V$ to be
the inverse image of the set of generic points of $\cX_k$; what remains is the
disjoint union of the inverse images of the smooth points of $\cX_k$, which are
open balls, and of the nodal points of $\cX_k$, which are open annuli.  For this
reason, we call a decomposition as in
Corollary~\ref{cor:semistable.decomposition} a \defi{semistable decomposition}.
It turns out that, conversely, a semistable decomposition also gives rise to a
semistable model: see~\cite[Theorem~4.11]{Baker2013}.

\subsection{Skeletons}

\begin{figure}[ht]
  \centering
  \includegraphics{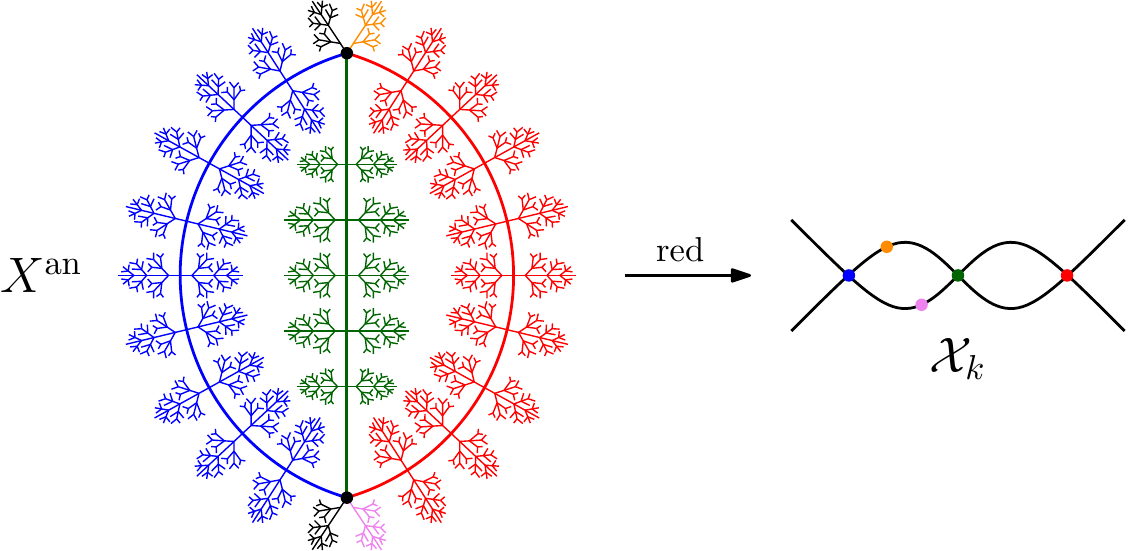}
  
  \caption{\textit{Left:} the analytification of a curve $X$ of genus at least $2$.
    \textit{Right:} the special fiber of a semistable model $\cX$ of this curve.  The
    thicker lines in $X^\an$ denote the skeleton $\Gamma_\cX$.  The annuli in
    red, green, and blue map to the indicated nodes in $\cX_k$ under $\red$.
    The vertices of $\Gamma_\cX$ map to the generic points of $\cX_k$.  The open
    balls in orange and purple map to the indicated smooth points of $\cX_k$.
    There are infinitely many open balls attached to (i.e.\ retracting to) each
    vertex; these correspond exactly to the smooth points of the corresponding
    component.}
  \label{fig:skeleton}
\end{figure}

Let $X$ be a smooth, proper, connected $K$-curve, and let
\begin{equation}\label{eq:ss.decomp}
X^\an = V \djunion\bS(\rho_1)_+\djunion\cdots\djunion\bS(\rho_n)_+
\djunion\Djunion\bB(1)_+
\end{equation}
be a semistable decomposition.  The \defi{skeleton} of $X$ associated to this
decomposition is the union of $V$ and the skeletons of the embedded open annuli:
\[ \Sigma(X,V) \coloneq
V\cup\Sigma(\bS(\rho_1)_+)\cup\cdots\cup\Sigma(\bS(\rho_n)_+). \]
An elementary argument in point-set topology implies that the closure of an open
annulus (resp.\ open ball) $U$ in the decomposition~\eqref{eq:ss.decomp}
consists of $U$ and one or two (resp.\ exactly one) point(s) of $V$.  From this
one can show that $\Sigma(X,V)$ is a closed subset of $X^\an$, which is
homeomorphic to a graph with vertex set $V$ and (open) edge set
$\{\Sigma(\bS(\rho_i)_+)\mid i=1,\ldots,n\}$.  In fact $\Sigma(X,V)$ is a
\emph{metric} graph: the length of the edge $\Sigma(\bS(\rho_i)_+)$ is by
definition $-\log\rho_i$, the logarithmic modulus of the open annulus
$\bS(\rho_i)_+$.  To summarize, a semistable decomposition of $X^\an$ gives rise
to a skeleton, which is a naturally embedded metric graph.  As a set,
$\Sigma(X,V)$ is the collection of all points of $X^\an$ which do not admit a
neighborhood disjoint from $V$ and isomorphic to an open unit ball.  The
skeleton \emph{$\Sigma(X,V)$ does not contain any $K$-points of $X$}: indeed,
the center of any point $x\in\Sigma(X,V)$ is the generic point of $X$.  See
Figure~\ref{fig:skeleton}, and refer to~\cite[\S3]{Baker2013} for details.

\subsubsection{The skeleton associated to a semistable model}
Suppose now that the semistable decomposition~\eqref{eq:ss.decomp} is the set of
fibers of the reduction map with respect to a semistable model $\cX$, as in
Theorem~\ref{thm:formal.fibers}.  In this case we will generally denote the
skeleton $\Sigma(X,V)$ by $\Gamma_\cX$.  The vertices in $V$ correspond to
the irreducible components of $\cX_k$: to be precise, $x\in V$ corresponds to
the closure of the generic point $\red(x)\in\cX_k$.  The edges in $\Gamma_\cX$
correspond, again via the reduction map, to the nodal points in $\cX_k$.  One
shows using the anti-continuity of the reduction map that, for a nodal point
$\td x\in\cX_k$, the skeleton of its inverse image
$\red\inv(\td x)\cong\bS(\rho)_+$ is the edge connecting the vertices
corresponding to the component(s) of $\cX_k$ containing $\td x$.  In other
words, in this case $\Gamma_\cX$ is naturally identified with the
\emph{incidence graph} of the components of the special fiber $\cX_k$.  In
particular, since $X$ admits infinitely many semistable models, it also has
infinitely many skeletons.  Note that the loop edges in $\Gamma_\cX$ correspond
to self-intersections of irreducible components of $\cX_k$.
See Figure~\ref{fig:skeleton} and~\cite[\S3]{Baker2013}.

\subsubsection{Retraction to the skeleton}\label{sec:retraction-skeleton}
Define a map $\tau\colon X^\an\to\Sigma(X,V)$ in the following way.  For
$x\in V$ we set $\tau(x) = x$.  For $x\in\bS(\rho_i)_+$ we set
$\tau(x) = \trop(x)\in\Sigma(\bS(\rho_i)_+)$, the retraction map
of~\S\ref{sec:skeleton-an-annulus}.  If $x$ is contained in an open ball $\bB(1)_+$ in
the semistable decomposition~\eqref{eq:ss.decomp}, we set $\tau(x)$ to be the
unique vertex $y\in V$ contained in the closure of $\bB(1)_+$.  The resulting
function $\tau\colon X^\an\to\Sigma(X,V)$ is a continuous retraction mapping.
Berkovich showed that $\tau$ is in fact the image of a strong deformation
retraction, so that $X^\an$ has the homotopy type of its skeleton.  See~\cite[Chapter~4]{Berkovich:spectralTheory}.

\begin{egsub}\label{eg:inverse.image.retraction}
  Suppose again that the semistable decomposition~\eqref{eq:ss.decomp} comes
  from a semistable model $\cX$ of $X$.  Let $\td x\in\cX_k$ be a nodal point,
  so $\red\inv(\td x)\cong\bS(\rho)_+$ is an open annulus.  By construction,
  \[ \tau\inv(\Sigma(\bS(\rho)_+)) = \red\inv(\td x) = \bS(\rho)_+, \]
  so the inverse image of an open edge under retraction is an open annulus.  Now
  let $\td x\in\cX_k$ be a generic point, let $x\in V$ be the unique point
  reducing to $\td x$, and let $C\subset\cX_k$ be the closure of $\td x$.  Using
  the anti-continuity of the reduction map, one shows in this case that
  \[ \tau\inv(x) = \red\inv(C^\sm), \]
  the inverse image of the set of smooth points of $C$ under reduction.
  Assuming $\cX$ is not itself smooth, so that $\Sigma(X,V)\neq\{x\}$, then
  $C^\sm$ is affine and $\tau\inv(x)$ is an \defi{affinoid domain} of $X^\an$.
\end{egsub}

\begin{rem*}
  The retraction $\tau\colon X^\an\to\Sigma(X,V)$ is very much analogous to the
  canonical deformation retraction of a once-punctured Riemann surface onto an
  embedded metric graph, important in the study of Teichm\"uller space.  Beware
  however that the first Betti number of $\Sigma(X,V)$ (as a topological
  space) is at most $g$, the genus of $X$, whereas in the setting of Riemann
  surfaces, the first Betti number of the skeleton is exactly $2g$.
\end{rem*}

\subsubsection{Decomposition into wide open subdomains}\label{sec:decomp-into-wide}
\begin{figure}[ht]
  \centering
  \includegraphics{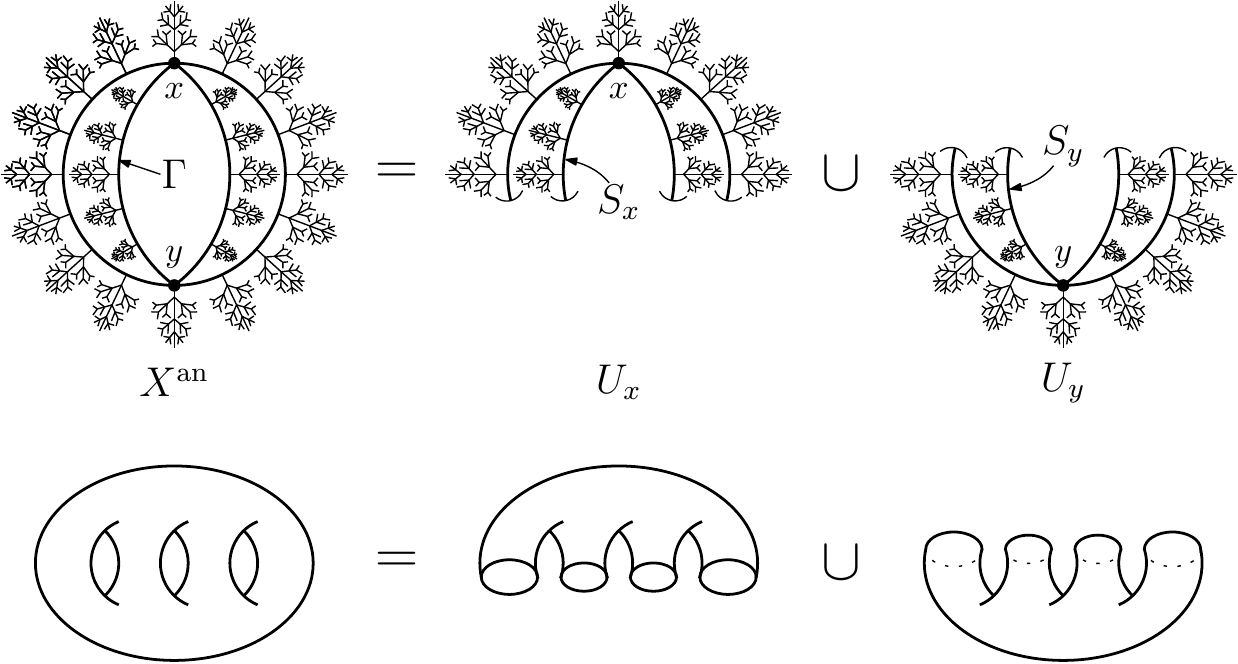}
  
  \caption{\textit{Above:} a ``pair-of-pants'' decomposition of a
    non-Archimedean curve $X^\an$ into basic wide open subdomains.  The thicker
    lines represent the skeleton $\Gamma = S_x\cup S_y$.
    \textit{Below:} the analogous picture for Riemann surfaces.  In both
    pictures, the intersection of the two open subsets is a disjoint union of
    four open annuli (cf.\ Figure~\ref{fig:annulus_skeleton}).}
  \label{fig:pants}
\end{figure}

The retraction to the skeleton allows us to make the following kind of
decomposition of $X^\an$.  For simplicity, we assume that
$\Gamma\coloneq\Sigma(X,V)$ does not contain any loop edges, i.e., that it
corresponds to a semistable model with smooth components (no
self-intersections).  For each vertex $x\in V$, let $S_x$ denote a
\defi{star neighborhood} around $x$ in $\Gamma$: this is the union of $x$ with a
connected open neighborhood of $x$ in each of the open edges adjacent to $x$.
We assume that the star neighborhoods $\{S_x\mid x\in V\}$ cover $\Gamma$.  Let
$U_x = \tau\inv(S_x)\subset X^\an$.  This open subspace of $X^\an$ is called a
\defi{basic wide open subdomain} in Coleman's
terminology~\cite[\S3]{ColemanReciprocity} (assuming the edge lengths are
contained in $v(K^\times)$).  See~\cite[Remark~2.20]{KatzRZB-uniform-bounds}.
If $x,y\in V$ share exactly one edge $e$, then $U_x\cap U_y = \tau\inv(e')$ for
the open interval $e' = S_x\cap S_y\subset e$.  Hence $U_x\cap U_y$ is
isomorphic to an open sub-annulus of the open annulus $\tau\inv(e)$.  In
general, $U_x\cap U_y$ will be isomorphic to a disjoint union of open annuli.
For this reason, the collection $\{U_x\mid x\in V\}$ is analogous to a
pair-of-pants decomposition of a Riemann surface: it is a collection of
connected open subsets that intersect along disjoint open annuli.  However, in
the non-Archimedean world, we cannot limit our pants to having only three
``legs'' (or two legs and a waist loop) --- indeed, the number of ``legs'' of
$U_x$ is the valency of the vertex $x$ in $\Gamma$, which cannot in general be
decreased.  See Figure~\ref{fig:pants}.

\subsection{Potential theory on Berkovich curves}\label{sec:potent-theory-berk}
In this subsection we fix a smooth, proper, connected $K$-curve $X$ with a
semistable $R$-model $\cX$ and associated skeleton $\Gamma\subset X^\an$.  Let
$f\in K(X)^\times$ be a nonzero rational function on $X$.  Any point of $\Gamma$
is centered at the generic point of $X$, so $|f(x)| = \|f\|_x$ is well-defined
and positive for all $x\in\Gamma$.  Let $F = -\log|f|\colon\Gamma\to\R$.

\begin{lemsub}\label{lem:log.PL}
  With the above notations, $F$ is piecewise affine-linear with integer slopes
  on $\Gamma$.
\end{lemsub}

By ``piecewise affine-linear with integer slopes,'' we mean that the restriction
of $F$ to each edge $e\subset\Gamma$ is differentiable at all but finitely many
points, and that the slopes of $F$ are integers, \emph{with respect to} either
of the two identifications of $e$ with an interval $(0,r)\subset\R$
(see~\S\ref{sec:skeleton-an-annulus}).  The proof of Lemma~\ref{lem:log.PL} is a
simple Newton polygon argument as applied to each open annulus in the semistable
decomposition of $X^\an$, but is perhaps best understood by example.
See also~\cite[Proposition~2.10]{Baker2013}.

\begin{egsub}\label{eg:log.PL}
  Let $f = p + t\in \C_p[t]$, regarded as a meromorphic function on
  $\bP^1_{\C_p}$.  Consider the embedded open annulus
  $A = \bS(|p|^{-2})_+\subset\bA^{1,\an}$ with skeleton $\Sigma(A) = (0,2)$.
  Recall (Example~\ref{eg:supremum.norm}) that a point $r\in(0,2)$ is identified
  with the seminorm $\|\scdot\|_r\in\Sigma(A)$ defined by
  $\|\sum a_it^i\|_r = \max\{|a_i|\rho^i\}$, where $\rho = \exp(-r)$.  Letting
  $F = -\log|f|\colon\Sigma(A) = (0,2)\to\R$, then, we have
  \[ F(r) = -\log\big(\max\{p\inv, \rho\}\big) = \min\{1, r\} =
  \begin{cases}
    r & \text{ if } r \leq 1 \\
    1 & \text{ if } r \geq 1,
  \end{cases} \]
  since in this case $a_0 = p$ and $a_1 = 1$.  Hence $F$ has one point of
  non-differentiability at $r=1$, with integer slopes $1$ and $0$.  See
  Figure~\ref{fig:graph}.

  \begin{figure}[ht]
    \centering
    \includegraphics{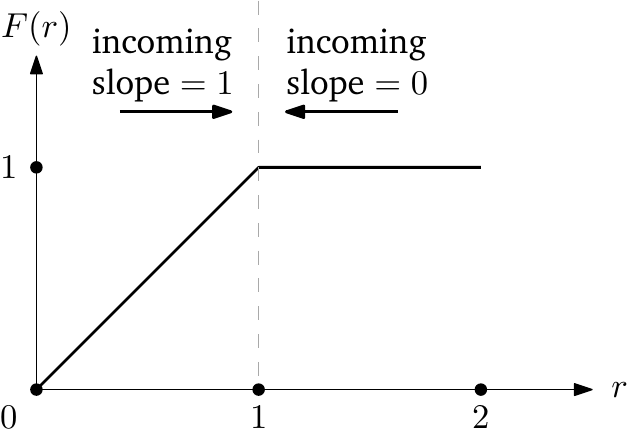}
    
    \caption{The function $F$ of Example~\ref{eg:log.PL}.  The divisor of $F$ is
    $(1)$: at that point, the sum of the incoming slopes is $1$, and at any
    other point $r\in(0,2)$, the incoming slope on one side equals the outgoing
    slope on the other side.}
    \label{fig:graph}
  \end{figure}
\end{egsub}

We call a piecewise affine-linear function with integer slopes
$F\colon\Gamma\to\R$ a \defi{tropical meromorphic function} on the metric graph
$\Gamma$.  If we declare that $-\log|f|$ is the \defi{tropicalization} of a
nonzero meromorphic function $f\in K(X)^\times$, then the tropicalization of a
meromorphic function is a tropical meromorphic function.

Let $F\colon\Gamma\to\R$ be a tropical meromorphic
function.  For $x\in\Gamma$ we define $\ord_x(F)\in\Z$ to be the sum of the
\emph{incoming} slopes of $F$ at $x$.  In other words, there are a number of
``tangent directions'' $\bv$ at $x$ in $\Gamma$, and
$\ord_x(F) = -\sum_\bv d_\bv F$, where $d_\bv F\in\Z$ is the derivative of $F$
in the direction $\bv$ (always with respect to the edge lengths).  The
\defi{divisor} of $F$ is the formal sum
\[ \div(F) \coloneq \sum_{x\in\Gamma}\ord_x(F)\cdot(x). \]
At almost all points $x$ in the interior of an edge, $F$ will be differentiable
at $x$, so that the incoming slope of $F$ in one direction equals the outgoing
slope in the other direction; it follows that $\ord_x(F) = 0$ at such a point,
so that $\div(F)$ is a finite sum.

\begin{egsub}\label{eg:log.PL.2}
  Continuing with Example~\ref{eg:log.PL}, on $\Sigma(A)$ we have
  $\ord_r(F) = 0$ for $r\neq 1$: for instance, for $r < 1$ the incoming slope of
  $F$ from the negative direction is $1$ and the incoming slope in the positive
  direction is $-1$ (since $F$ increases as $r$ increases).  We have
  $\ord_1(F) = 1$, since at that point the incoming slope from the positive
  direction is $0$ and the incoming slope from the negative direction is $1$.
  It follows that $\div(F) = (1)$.  See Figure~\ref{fig:graph}.
\end{egsub}

We denote the group of divisors on $\Gamma$, i.e., the free abelian group on the
points of $\Gamma$, by $\Div(\Gamma)$.  The retraction map
$\tau\colon X(K)\inject X^\an\to\Gamma$ of~\S\ref{sec:retraction-skeleton}
extends by linearity to a map 
\[ \tau_*\colon\Div(X)\to\Div(\Gamma) .\]
The following deep result is a restatement of Thuillier's Poincar\'e--Lelong
formula in non-Archimedean harmonic analysis, translated into this tropical
language.  We call it the \defi{Slope Formula for meromorphic functions}.
See~\cite[Proposition~3.3.15]{Thuillier2005}
and~\cite[Theorem~5.15 and Remark~5.16]{Baker2013}.

\begin{thmsub}[Slope Formula]\label{thm:poincare-lelong}
  Let $X$ be a smooth, proper, connected $K$-curve, and let
  $\Gamma\subset X^\an$ be a skeleton.  Let $f\in K(X)^\times$ and let
  $F = -\log|f|\colon\Gamma\to\R$.  Then
  \[ \tau_*(\div(f)) = \div(F). \]
\end{thmsub}

\begin{egsub}\label{eg:log.PL.3}
  Continuing with Examples~\ref{eg:log.PL} and~\ref{eg:log.PL.2}, we have
  $\div(f) = (p)$ and $\div(F) = (1)$.  By definition
  (see~\S\ref{sec:skeleton-an-annulus}),
  \[ \tau(p) = \trop(p) = -\log|t(p)| = -\log|p| = v(p) = 1, \]
  so Theorem~\ref{thm:poincare-lelong} recovers the fact that $f$ has a simple
  zero with valuation $1$.
\end{egsub}

In fact, for any point $x$ in the interior of an edge $e\subset\Gamma$, the
equality $\ord_x(F) = \sum_{\tau(y) = x}\ord_y(f)$ is more or less equivalent to
the theorem of the Newton polygon (see, e.g.,~\cite[\S6.4]{gouvea:pAdicNumbers-an-introduction}) as applied to the restriction of $f$ to the
open annulus $\tau\inv(e)$; the real content of
Theorem~\ref{thm:poincare-lelong} is that the formula also holds true at the
\emph{vertices} of $\Gamma$.  
The following corollary is then a purely combinatorial consequence of
Theorem~\ref{thm:poincare-lelong}.

\begin{corsub}\label{cor:zeros.from.slopes}
  With the notation in Theorem~\ref{thm:poincare-lelong}, let $U_x\subset X^\an$
  be a basic wide open subdomain for a vertex $x\in\Gamma$.  Suppose that $f$
  has no poles on $U_x(\C_p)$.  Let $e_1,\ldots,e_n$ be the edges of $\Gamma$ in
  $U_x$ adjacent to $x$, and let $\bv_i$ be the tangent direction at the other
  vertex of $e_i$, in the direction of $x$.  Then the number of zeros of $f$ on
  $U_x$ (counted with multiplicity) is equal to $\sum_{i=1}^n d_{\bv_i}F$.
\end{corsub}

We will apply Corollary~\ref{cor:zeros.from.slopes} to the antiderivative of an
exact $1$-form, in order to bound the number of rational points or torsion
points on $U_x$, in the style of Chabauty--Coleman.

\begin{egsub}\label{eg:log.PL.4}
  Continuing with Examples~\ref{eg:log.PL}, \ref{eg:log.PL.2},
  and~\ref{eg:log.PL.3}, we note that the annulus $A$ is a basic wide open
  subdomain with respect to any point $x$ on its skeleton, which we identify
  with the interval $(0,2)$.  The slope of $F$ at $0$ in the direction of $x$ is
  equal to $1$, and the slope of $F$ at $2$ in the direction of $x$ is $0$.
  Hence Corollary~\ref{cor:zeros.from.slopes} asserts that $f$ has a single zero
  on $A$.
\end{egsub}

Again, the result of Corollary~\ref{cor:zeros.from.slopes} is not hard to see
using Newton polygons when $U_x$ is an annulus; the reader might find this a
helpful first exercise before proving the general case.

\subsubsection{Model metrics on line bundles}
We will need a variant of Theorem~\ref{thm:poincare-lelong} which applies to
meromorphic sections of line bundles.  Let $\cL$ be an invertible sheaf on our
semistable model $\cX$, with generic fiber $L$.  Let $x\in X^\an$, let $y\in X$
be its center, and let $\iota_y\colon\Spec(\kappa(y)^\circ)\to\cX$ be the
extension of the inclusion $\Spec(\kappa(y))\inject X$, as explained
in~\S\ref{sec:reduction}.  Let $s$ be a nonzero meromorphic section of $L$, and
write $s = ft$ on an open neighborhood $\cU$ of $\red(x)$, where $t$ is a
nonvanishing section of $\cL$ on $\cU$ and $f$ is a meromorphic function on $U$.
We define
\[ \|s(x)\|_\cL \coloneq |f(x)| = \|\iota_y^*(f)\|_x. \]
Any unit $u$ in a neighborhood of $\red(x)$ pulls back via $\iota_y$ to a unit
in $\kappa(y)^\circ$, hence satisfies $\|\iota_y^*(u)\|_x = |u(x)| = 1$, so
$\|s(x)\|_\cL$ is well-defined.  We call $\|\scdot\|_\cL$ the
\defi{model metric} on $L$ associated to its integral model $\cL$.  By choosing
local sections, it follows from Lemma~\ref{lem:log.PL} that
$F\coloneq-\log\|s\|_\cL$ is a tropical meromorphic function.

\begin{rem*}
  Model metrics have the following intersection-theoretic interpretation over a
  discretely valued field $K$ (note that the definition of $\|\scdot\|_\cL$
  above does not use that $K$ is algebraically closed).  Suppose that $\Z$ is
  the value group of $K$.  For simplicity we restrict ourselves to a regular
  split semistable model $\cX$ of a smooth, proper, geometrically connected
  curve $X$.  A meromorphic section $s$ of $L$ can be regarded as a meromorphic
  section of $\cL$, hence has an order of vanishing $\ord_D(s)$ along any
  irreducible component $D$ of $\cX_k$.  If $\zeta\in X^\an$ is the point
  reducing to the generic point of $D$ then we have the equality
  \[ -\log\|s(\zeta)\|_\cL = \ord_D(s). \]
  This follows from the observation that $\ord_D\colon K(X)^\times\to\Z$ is also
  a valuation such that the induced map $\Spec(K(X)^\circ)\to\cX$ takes the
  closed point to the generic point of $D$.
\end{rem*}

We can now state the general slope formula, which can be derived from
Theorem~\ref{thm:poincare-lelong}.  See~\cite[Theorem~2.6]{KatzRZB-uniform-bounds}.

\begin{thmsub}[Slope Formula for line bundles]\label{thm:slope.formula}
  Let $X$ be a smooth, proper, connected $K$-curve, and let $\cX$ be a
  semistable $R$-model of $X$ with corresponding skeleton
  $\Gamma_\cX\subset X^\an$.  Assume that $\cX$ is not smooth, so that
  $\Gamma_\cX$ is not a point.  Let $\cL$ be a line bundle on $\cX$, let
  $L = \cL|_X$, let $s$ be a nonzero meromorphic section of $L$, and let
  $F = -\log\|s\|\colon\Gamma_\cX\to\R$.  Then
  \[ \tau_*(\div(s)) = \div(F) + \sum_\zeta\deg(\cL|_{D_\zeta})\,(\zeta), \]
  where the sum is taken over vertices $\zeta$ of $\Gamma_\cX$, and $D_\zeta$ is
  the irreducible component of $\cX_k$ with generic point $\red(\zeta)$.
\end{thmsub}

\subsection{The divisor of a regular differential}\label{sec:divis-regul-diff}
Recall that $\cX$ is a semistable model of $X$.  We take $L = \Omega^1_{X/K}$, the
cotangent sheaf.  This invertible sheaf has a canonical extension to $\cX$,
namely, the relative dualizing sheaf $\Omega^1_{\cX/R}$.  (Actually, the theory
of the relative dualizing sheaf is only well-developed for noetherian schemes,
of which $\cX$ is not an example.  We will ignore this technical difficulty
entirely; see~\cite[\S2.4]{KatzRZB-uniform-bounds} for details.)  We write $\|\scdot\| = \|\scdot\|_{\Omega^1_{\cX/R}}$ for the
corresponding model metric.  The adjunction formula implies that if
$D\subset\cX_k$ is an irreducible component, then
\begin{equation}\label{eq:adjunction}
\deg(\Omega^1_{\cX/R}|_D) = 2g(D) - 2 + r(D),
\end{equation}
where $g(D)$ is the geometric genus of $D$ and $r(D)$ is the number of points of
the normalization of $D$ which map to singular points of $\cX_k$, i.e., the
number of nodes lying on $D$, counting self-intersections twice.
Theorem~\ref{thm:slope.formula} and~\eqref{eq:adjunction} imply that if
$\omega\in H^0(X,\Omega^1_{X/K})$ is a nonzero regular differential and
$F = -\log\|\omega\|$, then
\begin{equation}\label{eq:effectivity.1}
  \div(F) + \sum_\zeta \big(2g(D_\zeta) - 2 + r(D_\zeta)\big)\,(\zeta) \geq 0,
\end{equation}
where the sum is defined as in Theorem~\ref{thm:slope.formula}.

We have the following tropical interpretation of~\eqref{eq:effectivity.1}.  Let
$\Gamma = \Gamma_\cX$ be the skeleton of $X$ associated to the model $\cX$.  For
a vertex $x\in\Gamma$, we let $\deg(x)$ denote the valency of $x$ in $\Gamma$,
and we let $g(x)$ denote the geometric genus of the irreducible component of
$\cX_k$ with generic point $\red(x)$.  If $x\in\Gamma$ is not a vertex then we
set $\deg(x) = 2$ and $g(x) = 0$.  The function $g\colon\Gamma\to\Z$ is a
\defi{weight function} on the vertices of $\Gamma$; hence $(\Gamma,g)$ is a
\defi{vertex-weighted metric graph}.  The \defi{canonical divisor} of the
vertex-weighted metric graph $(\Gamma,g)$ is
\begin{equation}\label{eq:canonical.Gamma}
  K_\Gamma \coloneq \sum_{x\in\Gamma} \big(2g(x) - 2 + \deg(x)\big)\,(x).
\end{equation}
Note that $K_\Gamma$ is supported on the vertices of $\Gamma$.  Since the edges
adjoining a vertex $x$ correspond to the nodal points lying on the irreducible
component of $\cX_k$ with generic point $\red(x)$ (again counting
self-intersections twice), the (purely combinatorial) definition of $K_\Gamma$
precisely encodes the multi-degree of the relative dualizing sheaf
$\Omega^1_{\cX/R}$ restricted to $\cX_k$.  Hence we obtain the following result.

\begin{corsub}\label{cor:effectivity}
  Let $X$ be a smooth, proper, connected $K$-curve, and let $\cX$ be a
  semistable $R$-model of $X$, with corresponding skeleton
  $\Gamma_\cX\subset X^\an$.  Assume that $\cX$ is not smooth, so that
  $\Gamma_\cX$ is not a point.  Let $\omega\in H^0(X,\Omega^1_{X/K})$ be a
  nonzero regular differential, and let
  $F = -\log\|\omega\|_{\Omega^1_{\cX/R}}\colon\Gamma_\cX\to\R$.  Then
  \[ \div(F) + K_{\Gamma_\cX} \geq 0. \]
\end{corsub}

If we decide that $F$ is the ``tropicalization'' of $\omega$, and that a
``section of the tropical canonical bundle'' on a weighted metric graph $\Gamma$
is a tropical meromorphic function $F\colon\Gamma\to\R$ such that
$\div(F) + K_{\Gamma}\geq 0$, then Corollary~\ref{cor:effectivity} asserts that
\begin{quotation}
  \emph{The tropicalization of a section of the canonical bundle is a section of
  the tropical canonical bundle.}
\end{quotation}

\begin{rem*}
  We should mention that the theory of divisors and linear equivalence on
  graphs, initiated primarily by Baker and Norine, has a rich and beautiful
  combinatorial structure that mirrors the analogous theory for algebraic
  curves.  For instance, there is a Riemann--Roch theorem in this context.
  See~\cite{BakerN:RR,Baker:specialization,BakerN:Harmonic}, for instance.
\end{rem*}

The \defi{genus} of a vertex-weighted metric graph $(\Gamma,g)$ is by definition
\[ g(\Gamma) \coloneq \sum_{x\in\Gamma} g(x) + h_1(\Gamma), \]
the sum of the weights of the vertices and the first Betti number of the graph
(as a simplicial complex).  If $\Gamma$ is the skeleton associated to a
semistable model $\cX$ of $X$ as above, then a standard calculation shows that
$g(\Gamma)$ is the arithmetic genus of $\cX_k$, and hence that $g(\Gamma)$ is
the genus of $X$.  See also~\cite[Theorem~4.6]{BoschL:stableReductionAnd}.

Consider now the following lemma, whose statement and proof are purely combinatorial.  

\begin{lemsub}[{\cite[Lemma~4.15]{KatzRZB-uniform-bounds}}]\label{lem:slope.bound}
  Let $(\Gamma,g)$ be a vertex-weighted metric graph of genus $g(\Gamma)$.  Let
  $F$ be a tropical meromorphic function on $\Gamma$ such that
  $\div(F)+ K_\Gamma\geq 0$.  Then for all $x\in \Gamma$ and all tangent
  directions $\bv$ at $x$, we have $|d_\bv F|\leq 2g(\Gamma)-1$.
\end{lemsub}

In other words, if $F$ is a section of the tropical canonical bundle, then all
slopes of $F$ are bounded by $2g(\Gamma)-1$.  Lemma~\ref{lem:slope.bound} and
Corollary~\ref{cor:effectivity} together give a bound on the slopes of the
tropicalization of a regular differential, which will be a key ingredient in our
application of the Chabauty--Coleman method in the sequel.  This also
demonstrates the utility of the non-Archimedean analytic language in reducing
algebro-geometric problems to well-studied combinatorial questions.

\section{Theories of $p$-adic Integration}
\label{sec:integration}

In this section, we fix a smooth, proper, connected $\C_p$-curve $X$, along with
a semistable model $\cX$, in the sense of~\S\ref{sec:reduction}.  Let $J$ be the
Jacobian of $X$.  It is known that $J$ extends to a smooth group scheme $\cJ$
over $R$, the ring of integers of $\C_p$.  This is due to the theory of N\'eron
models if $J$ is defined over a finite extension of $\Q_p$; otherwise, use~\cite[\S5]{BoschL:stableReductionAnd2}.

The discussion of $p$-adic integration in~\S\ref{sec:chabauty} was limited to
the following: on an open ball, the restriction of a regular $1$-form $\omega$
is exact, and thus has the form $df$ for some analytic function $f$. The
integral is then computed via the primitive $f$ as
\[
\int^P_Q \omega = \int^{P}_{Q} df \coloneq f(Q) - f(P).
\]
Restricting to such \defi{tiny integrals} -- i.e., those between points in the same
residue class -- suppresses a major technical difficulty: integrating between
residue classes.  Thankfully, tiny integrals are all that are needed for many
classical applications of the method, including both Chabauty and Coleman's
theorems.

If $P$ and $Q$ do not lie in the same tube, then there are \emph{multiple} ways
to $p$-adic analytically continue the integral $\int^P_Q \omega$.  We will
discuss two of them, namely, \defi{abelian integration} and
\defi{Berkovich--Coleman integration}.  For simplicity we restrict to
integrating regular $1$-forms between $\C_p$-points of $X$.

\subsection{Abelian integration}\label{sec:abelian-integration}
The group $J(\C_p)$ is a $\C_p$-Lie group, in the na\"ive sense that it locally
looks like an open neighborhood of $\C_p^g$, with ``smooth'' transition
functions given by convergent power series.  Such $p$-adic manifolds were
studied by Bourbaki in their treatise on Lie groups and Lie algebras \cite{Bourbaki-lie-groups-7-9}; using general considerations, one can prove that
there exists a unique homomorphism of $\C_p$-Lie groups
$\log\colon J(\C_p)\to\Lie(J)\cong\C_p^g$ whose linearization
\[ d\log\colon\Lie(J)\to\Lie(\Lie(J)) = \Lie(J) \]
is the identity map. See~\S\ref{sec:abel-integr-form} below for an
algebro-geometric construction in our situation, though.  Since $\C_p^g$ has no additive
torsion, the full torsion subgroup of $J(\C_p)$ is contained in $\ker(\log)$.

For $P\in J(\C_p)$ and $\omega\in\Omega_{J/\C_p}^1(J)$ we define
\[ \Abint_0^P\omega = \angles{\log(P),\,\omega} \]
where $\angles{\scdot,\scdot}$ is the pairing between $\Lie(J)$ and
$\Omega_{J/\C_p}^1(J)$.  For $P,Q\in J(\C_p)$ we set
\[ \Abint_P^Q\omega = \Abint_0^Q\omega -
\Abint_0^P\omega. \]
We call $\Abint$ the \defi{abelian integral} on $J$.

\subsubsection{The abelian integral and formal antidifferentiation}
\label{sec:abel-integr-form}
Let $\fJ$ be the completion of $\cJ$ along its identity section.  Choosing
coordinates, $\fJ\cong\Spf(R\ps{x_1,\ldots,x_g})$ is a commutative
$g$-dimensional formal group over $R$.  Let $\bF(\bx,\by)$ be the formal group
law, where $\bF = (F_1,\ldots,F_g)$, $\bx = (x_1,\ldots,x_g)$, etc.

Any cotangent vector $\bv = \sum_{i=1}^g a_i\,dx_i$ at the identity can be
extended uniquely (by translation) to give a translation-invariant $1$-form
$\omega$ on $\fJ$, and similarly for $q$-forms for any $q\geq 0$.  The power
series defining multiplication by $2$ has the form $[2](\bx) = 2\bx + $
higher-order terms, so $[2]^*\bv = 2\bv$, and hence $[2]^*\omega = 2\omega$ for
a translation-invariant $1$-form $\omega$.  Similarly, $[2]^*\eta = 2^q\eta$ for
a translation-invariant $q$-form $\eta$.  Taking $\eta = d\omega$ for a
translation-invariant $1$-form $\omega$, we have
\[ 2\,d\omega = d(2\omega) = d([2]^*\omega) = [2]^*d\omega = 4\,d\omega, \]
so any such $\omega$ is closed.  Writing $\omega = \sum_{i=1}^g f_i\,dx_i$ for
$f_i\in R\ps\bx$, since $\omega$ is closed, we can formally antidifferentiate
the $f_i$ so that $\omega = dh$ for $h\in\C_p\ps\bx$.  In other words,
translation-invariant $1$-forms are \emph{exact} on the generic fiber.

Let $\omega = dh$ be a translation-invariant $1$-form, always choosing the
antiderivative $h$ to have zero constant term.  We claim that $h$ defines a
homomorphism of formal groups from $\fJ_{\C_p}$ to $\hat\bG_{a,\C_p}$, i.e.,
that $h(F(\bx,\by)) = h(\bx) + h(\by)$.  Using translation-invariance, we have
\[ dh = \omega = T_{\by}^*\omega = d(T_{\by}^*h) = d(h(F(\bx,\by))), \]
where $T_{\by}$ is translation by $\by$, and differentiation is taken
with respect to $\bx$.  It follows that $h(F(\bx,\by)) = h(\bx) + c(\by)$, where
$c(\by)\in R\ps\by$ is the ``constant of integration.''  Substituting $\bx=\bb0$
gives $h(\by) = h(F(\bb0,\by)) = c(\by)$, which proves the claim.

The above association $\omega\mapsto h$ gives a homomorphism
$\fJ_{\C_p}\to\Lie(\fJ_{\C_p})$, in that a translation-invariant $1$-form gives
rise to a homomorphism $\fJ_{\C_p}\to\hat\bG_{a,\C_p}$; almost by definition, the
linearization of this homomorphism is the identity on $\Lie(\fJ_{\C_p})$.  In
coordinates, we have a basis $dx_1,\ldots,dx_g$ for the cotangent space of $\fJ$
(or of $\cJ$, or $J$) at the identity.  Let $\omega_i = dh_i$ be the
translation-invariant $1$-form extending $dx_i$.  Since $\omega_i = dx_i$ at the
identity, the linear term of $h_i$ is $x_i$, and therefore
$(h_1,\ldots,h_g)\colon\fJ_{\C_p}\to\hat\bG_a^g$ takes $\del/\del x_i$ to the
$\del/\del t_i$, where $t_i$ is the coordinate on the $i$th copy of $\hat\bG_a$.

\begin{egsub}\label{eg:Gm.1forms}
  So far our discussion has been intrinsic to the formal group $\fJ$, so we may
  take $\fJ=\hat\bG_m$, with group law $F(x,y) = xy + x + y$ and inverse
  $I(x) = (1+x)\inv-1 = -x + x^2 - x^3 + \cdots$.  If $\omega$ is the
  translation-invariant $1$-form associated to $dx$, then
  \[\begin{split}
    \omega(y) &= T_{I(y)}^*dx = \frac d{dx}(F(x,I(y)))\,dx \\
    &= \frac d{dx}(xI(y) + x + I(y))\,dx = (I(y) + 1)\,dx \\
    &= \big(1 - y + y^2 - y^3 + \cdots\big)\,dx.
  \end{split} \]
  (This is perhaps confusing because $dx$ denotes a cotangent vector field, so
  we need another variable $y$ for our power series.) Substituting $y=x$ again
  gives $\omega = dx/(1+x)$, the usual translation-invariant $1$-form.  Hence
  $\omega = d\log(1+x)$, where $\log(1+x) = x - x^2/2 + x^3/3 - \cdots$ is the
  Mercator series.  Of course $\log(1+x)$ defines a homomorphism from
  $\hat\bG_m$ to $\hat\bG_a$ over $\C_p$, which takes $\del/\del x$ to itself.
\end{egsub}

Returning to the discussion of Jacobians (really just abelian varieties), let
$J_+(\C_p)\subset J(\C_p)$ denote the set of points reducing to the identity in
the special fiber of $\cJ$.  Then $J_+(\C_p)$ is a \emph{subgroup} of $J(\C_p)$.
We have $J_+(\C_p)\cong\bB(1)^g_+(\C_p) = \fm_{\C_p}^g$ as sets, with the
coordinates $x_1,\ldots,x_g$ for the formal completion $\fJ$ defining this
bijection (via reduction modulo successively higher powers of $p$).  Any global
$1$-form $\omega\in H^0(J,\Omega^1_{J/\C_p})$ is translation-invariant due to
properness of $J$, hence exact when restricted to $\fJ_{\C_p}$.  The absolute values
of the coefficients of its formal antiderivative $h_\omega\in\C_p\ps{\bx}$ grow
at most logarithmically with the size of the exponent (i.e.\ $v(1/p^n) = -n$),
so the power series $h_\omega$ converges on $J_+(\C_p)$.  Since $h_\omega$ is
formally a homomorphism to $\hat\bG_a$, it also defines a homomorphism
$h_\omega\colon J_+(\C_p)\to\C_p$.  Therefore, the logarithm (in the sense
of~\S\ref{sec:abelian-integration}) is defined by
\[ \angles{\log(P),\,\omega} = h_\omega(P) \]
on $J_+(\C_p)$.  In summary, $\log$ is simply given by formal
antidifferentiation of a (translation-invariant) global $1$-form $\omega$ on
$J_+(\C_p)$.

One can show that $J(\C_p)/J_+(\C_p)$ is a torsion group. This is easiest to see
when $J$ is defined over a finite extension $F$ of $\Q_p$, in which case
$J(F)/J_+(F)$ injects into the group of closed points of the special fiber of
the N\'eron model of $J$.  Hence the logarithm can be defined on all points
$P\in J(\C_p)$ by multiplying by a suitable integer $n$ such that
$[n]P\in J_+(\C_p)$, then using the logarithm as defined above, then dividing by
$n$ in $\C_p$: that is, $\log(P) = n\inv\log([n]P)$.

\subsubsection{Abelian integration on a curve}\label{sec:abel-integr-curve}
Fix a base point $P_0\in X(\C_p)$, and let $\iota\colon X\inject J$ be the
Abel--Jacobi map with respect to $P_0$.  We use $\iota^*$ to identify
$H^0(J,\Omega^1_{J/\C_p})$ with $H^0(X,\Omega^1_{X/\C_p})$.  For
$P,Q\in X(\C_p)$ and $\omega\in H^0(J,\Omega^1_{J/\C_p})$ we define the
\defi{abelian integral} by
\[ \Abint_P^Q\iota^*\omega \coloneq \Abint_{\iota(P)}^{\iota(Q)}\omega. \]
The abelian integral is clearly independent of the choice of $P_0$.  Moreover,
it satisfies the following properties:
\begin{enumerate}
\item It is \emph{path-independent}, in that $\Abint_P^Q\omega$ makes no
  reference to a ``path'' from $P$ to $Q$.
\item For $P_1,P_2,P_3\in X(\C_p)$ and $\omega\in H^0(X,\Omega^1_{X/\C_p})$ we
  have
  \[ \Abint_{P_1}^{P_3}\omega = \Abint_{P_1}^{P_2}\omega +
  \Abint_{P_2}^{P_3}\omega. \]
\item For fixed $P,Q\in X(\C_p)$, the map $\omega\mapsto\Abint_P^Q\omega$ is
  $\C_p$-linear in $\omega$.
\item If $P,Q\in\red\inv(\td x)$ for $\td x\in\cX_k(k)$ a smooth point,
  then $\Abint_P^Q(\omega)$ is calculated by formally antidifferentiating
  $\omega$ with respect to a coordinate on $\red\inv(\td x)$.
\end{enumerate}
The only property which is not obvious from the definitions is the final one.
If $X$ is defined over a finite extension of $\Q_p$, then the N\'eron mapping
property implies that the Abel--Jacobi map extends to a map $\cX^\sm\to\cJ$; the
claim then follows from the discussion in~\S\ref{sec:abel-integr-form}.  In
general, it turns out that the claim is true for $P,Q$ contained in any open
subdomain of $X^\an$ which is isomorphic to $\bB(1)_+$;
see~\cite[Proposition~3.10]{KatzRZB-uniform-bounds}.  Of course, if
$U\subset X^\an$ is isomorphic to an open annulus, for instance, then there is
no reason for the abelian integral to be computed by formal antidifferentiation
on $U$, and in general \emph{it is not}.  This is a crucial point that Stoll
realized~\cite[Proposition~7.3]{stoll:uniform}.  See also~\S\ref{sec:berk-colem-integr}.

For us, the most important property of the abelian integral is:
\begin{quotation}
  Let $P,Q\in X(\C_p)$, and suppose that $[Q]-[P]$ represents a torsion point of
  $J(\C_p)$.  Then $\Abint_P^Q\omega = 0$ for all $\omega\in H^0(X,\C_p)$.
\end{quotation}

\subsection{Berkovich--Coleman integration}\label{sec:berk-colem-integr}
Historically, what we call the Berkovich--Coleman integration theory was first
developed for curves by Coleman and Coleman--de Shalit in \cite{ColemanDilogarithms,ColemanTorsion,CdS}.
Coleman's idea was to extend the tiny integrals on tubes by ``analytic
continuation by Frobenius'' -- essentially by insisting that the integral be
functorial with respect to pullbacks, and then considering pullbacks by
self-maps of affinoid subdomains of $X$ lifting the Frobenius on affine opens in
$\cX_k$.  This method is elegant and well-suited to computation; an algorithm for integration on hyperelliptic curves is given in \cite{balakrishnanBK:coleman}.

Berkovich \cite{Berkovich:integrationOfOne} then took the idea of analytic continuation by Frobenius
and radically extended it, constructing an integration theory that applies to
essentially any smooth analytic space.  It applies in particular to Jacobians,
so one can recover Coleman's integration of one-forms on curves (of arbitrary
reduction type) by understanding Berkovich integration on Jacobians and pulling
back by an Abel--Jacobi map.  We take this approach, as it is more suitable for
our purposes (as we would like to compare it to the abelian integral), and it is
in some sense just as explicit as Coleman's.

\subsubsection{Definition of the integral}
In order to define Berkovich--Coleman integration, one has to fix once and for
all a choice of a ``branch of the $p$-adic logarithm'', in the following sense.
The Mercator series $\log(1+x) = x - x^2/2 + x^3/3 - \cdots$ converges on
$1+\fm_{\C_p}$, the residue ball reducing to $1\in\bar\F_p$; this extends
uniquely to a homomorphism $\log\colon R^\times\to\C_p$, since every other
residue class contains a unique root of unity, which is killed by any
homomorphism to $\C_p$.  However, there are many ways to extend $\log$ to a
homomorphism $\Log\colon\C_p^\times\to\C_p$, since
$\C_p^\times\cong\Q\times R^\times$ as groups.  In effect, one has to (arbitrarily) fix
the value of $\Log(p)\in\C_p$; then for $x\in\Q$ and $u\in R^\times$, one has
$\Log(p^x u) = x\Log(p) + \log(u)$.  Note that this definition does not depend
on a choice of $x$-th power of $p$, since any other choice would differ by a
root of unity, which is killed by $\Log$.

For a smooth $\C_p$-analytic space $Y$, we let $\cP(Y)$ denote the set of
continuous paths between $\C_p$-points of $Y$, i.e., the set of continuous maps
$\gamma\colon[0,1]\to Y$ with $\gamma(0),\gamma(1)\in Y(\C_p)$, and we let
$Z^1_{\dR}(Y)$ denote the space of closed analytic $1$-forms on $Y$.  (In the
sequel, $Y$ will be the analytification of a smooth $\C_p$-variety, or an open
ball, or an open annulus.)  Again, the theory applies to much more general
differential forms; we restrict to regular one-forms for simplicity.

\begin{defnsub}\label{def:BC.int}
  The \defi{Berkovich--Coleman integration theory} is the unique pairing
  \[ \BCint\colon\cP(Y)\times Z^1_{\dR}(Y)\To\C_p \]
  for every smooth $\C_p$-analytic space $Y$, satisfying for all
  $\gamma\in\cP(Y)$ and $\omega\in Z^1_{\dR}(Y)$:
  \begin{enumerate}
  \item $\omega\mapsto\BCint_\gamma\omega$ is $\C_p$-linear in $\omega$ for
    fixed $\gamma$.
  \item $\BCint_\gamma\omega$ only depends on the fixed end-point homotopy class
    of $\gamma\in\cP(Y)$.
  \item If $\gamma,\gamma'\in\cP(Y)$ with $\gamma(1) = \gamma'(0)$, then
    \[ \BCint_{\gamma*\gamma'}\omega = \BCint_\gamma\omega +
    \BCint_{\gamma'}\omega, \]
    where $\gamma*\gamma'$ is the concatenation.
  \item If $f\colon Y\to Y'$ is a morphism and $\omega'\in Z^1_{\dR}(Y')$, then
    \[ \BCint_\gamma f^*\omega' = \BCint_{f(\gamma)}\omega'. \]
  \item If $\omega = df$ is exact, then
    $\BCint_\gamma\omega = f(\gamma(1)) - f(\gamma(0))$.
  \item If $Y = \bG_m^\an = \Spec(\C_p[t^{\pm1}])^\an$ and $\omega = dt/t$ is
    the invariant differential, then
    \[ \BCint_\gamma\frac{dt}t = \Log(\gamma(1)) - \Log(\gamma(0)). \]

  \end{enumerate}
\end{defnsub}

The existence and uniqueness of the Berkovich--Coleman integral is very deep,
and forms the content of Berkovich's book~\cite{Berkovich:integrationOfOne}.  We wish to state a
consequence of properties~(4)--(6) for emphasis, and for contrast with the
abelian integral:
\begin{quotation}
  The Berkovich--Coleman integral is local on $Y$, in that if $U\subset Y$ is an
  open subdomain and $\gamma\subset U$, then $\int_\gamma\omega$ can be computed
  on $U$ or on $Y$.  If $\omega$ happens to be exact on $U$, then
  $\int_\gamma\omega$ is computed by formal antidifferentiation.
\end{quotation}
We also point out that if $Y$ is simply-connected, then~(2) implies that
$\BCint$ is path-independent; in this case, we simply write
$\BCint_P^Q\coloneq\BCint_\gamma$ for any path $\gamma$ from $P$ to $Q$.

\subsubsection{Integration on totally degenerate Jacobians}
In complex analysis, a standard way to integrate holomorphic one-forms on an
abelian variety $A$ is first to pass to its universal cover, which is the vector
space $H^0(A,\Omega^1_{A/\C})^*$, where all closed holomorphic one-forms are
exact.  We proceed in essentially the same way for abelian varieties over
$\C_p$.  The analytification $Z^\an$ of a connected, smooth variety $Z$ is
locally contractible and locally path-connected, so it admits a universal cover
$\pi\colon\td Z\to Z^\an$ in the sense of point-set topology.  As in the complex
setting, the universal cover inherits a unique structure of analytic space
making $\pi$ into a local isomorphism.  However, the universal cover of an
abelian variety is no longer simply a vector space; in general it is much more
complicated.  We will discuss the uniformization theory of $J$, the Jacobian of
$X$, when $X$ is a Mumford curve, i.e., when $\cX_k$ has only rational
components.  In this case, the universal cover of $J^\an$ is an analytic torus,
and $J$ is said to be \defi{totally degenerate}.  Everything we will say
generalizes to the case of a general abelian variety, but this is much is more
technical, and the important ideas already appear in the totally degenerate
case.  For references, see~\cite{FdV} in the totally degenerate case,
and~\cite{BoschL:stableReductionAnd2,BoschL:degeneratingAbelianVarieties} in general.

\begin{egsub}\label{eg:unif.over.C}
  To motivate the non-Archimedean situation, first we recall one way to
  construct the Jacobian of a Riemann surface $X$.  There is a natural period
  mapping $\eta\colon H_1(X,\Z)\to H^0(X,\Omega^1)^*$ defined by
  \[ \gamma\mapsto \bigg(\omega\mapsto\int_\gamma\omega\bigg). \]
  One proves that the image of $\eta$ is a lattice.  The Jacobian of $X$ is
  \[ J = H^0(X,\Omega^1)^*/\eta H_1(X,\Z), \]
  which is a complex torus of dimension $g$.

  For simplicity, suppose now that $X$ is an elliptic curve.  Then
  $H^0(X,\Omega^1) \cong H^0(X,\sO) = \C$.  Fix a basis $e_1,e_2$ for
  $H_1(X,\Z)\cong\Z^2$, and choose an isomorphism
  $H^0(X,\Omega^1) \cong H^0(X,\sO)$ such that $e_1\mapsto\tau\in\bH$, the upper
  half plane, and $e_2\mapsto 1$.  Then $J = \C/\angles{1,\tau}$.  Composing
  with the exponential map $z\mapsto\exp(2\pi iz)\colon\C\to\C^\times$ kills
  $\Z$, so we have $J\cong\C^\times/q^\Z$, where $q = \exp(2\pi i\tau)$.

  In higher dimensions, in the complex situation one has a choice whether to think of the
  Jacobian $J$ as a quotient of a complex vector space, or of a complex torus.
  Only the latter viewpoint works in the non-Archimedean world.
\end{egsub}

A rough statement of the uniformization theorem for Jacobians of Mumford curves
is as follows.  See~\cite[\S7]{BoschL:stableReductionAnd2}
and~\cite[Chapter~6]{FdV}.

\begin{thmsub}\label{thm:unif.jac}
  Let $X$ be a genus-$g$ Mumford curve over $\C_p$ with Jacobian $J$.  Then
  there is a natural torus $T\cong\bG_m^g$ and a natural homomorphism
  $\eta\colon H_1(X^\an,\Z)\to T(\C_p)\cong(\C_p^\times)^g$ such that 
  $J^\an\cong T^\an/\eta H_1(X^\an,\Z)$.
\end{thmsub}

A few words about this theorem are in order.  First, $H_1(X^\an,\Z)$ is simply
the \emph{singular} homology of $X^\an$, in the sense of point-set topology.  It
is free of rank $g$.  The character lattice of $T$ is also canonically
isomorphic to $H_1(X^\an,\Z)$ -- this is due to the autoduality of $J$ -- hence
$J^\an$ is determined by a pairing
$H^1(X^\an,\Z)\times H^1(X^\an,\Z)\to\C_p^\times$.  We will use this fact later.
The action of $H_1(X^\an,\Z)$ on $T^\an$ is totally discontinuous, so that the
quotient makes sense, but in fact more is true.  Define
$\trop\colon(\C_p^\times)^g\to\Q^g$ by
\[ \trop(x_1,\ldots,x_g) = (v(x_1),\ldots,v(x_g)). \]
Choosing a basis for the character lattice of $T$, we can think of $\trop$ as a
homomorphism $T(\C_p)\to\Q^g$.  (If we did not want to choose a basis, then
$\trop$ would take values in $X_*(T)\tensor_\Z\Q$, where $X_*(T)$ is the
cocharacter lattice.)  Then it turns out that $\trop\circ\eta(H_1(X^\an,\Z))$ is
a \emph{lattice} in $\Q^g$, in that its (rank-$g$) image spans.
Total discontinuity of translation follows easily from this.

Given the homomorphism $\eta$ and the uniformization map
$\pi\colon T^\an\to J^\an$, it is now straightforward to compute the
Berkovich--Coleman integral on $J$.  Let $\omega\in H^0(J,\Omega^1_{J/\C_p})$
and let $\gamma$ be a path from the identity to a point $P\in J(\C_p)$.  Let
$\td\gamma$ be the unique lift of $\gamma$ to $T^\an$ starting at $1$, and let
$\td P = (x_1,\ldots,x_g)$ be the endpoint of $\td\gamma$.
Pulling back by $\pi$ gives an invariant one-form on $T^\an$, so
\[ \pi^*\omega = \sum_{i=1}^g a_i\frac{d t_i}{t_i}, \]
where $a_1,\ldots,a_g\in\C_p$ and $t_1,\ldots,t_g$ are coordinates on
$T\cong\bG_m^g$.  But by definition, $\BCint dt_i/t_i = \Log t_i$, so
\[ \BCint_\gamma\omega = \BCint_{\td\gamma} \pi^*\omega = \sum_{i=1}^g a_i
\Log(x_i). \]
We see then that the periods in the complex case have been replaced by the
$p$-adic logarithm of the image of the non-Archimedean uniformization lattice
$H_1(X^\an,\Z)\to(\C_p^\times)^g$.  Note though that the image of
$\Log\circ\eta$ cannot be a ``lattice'' in $\C_p^g$, since any nonzero subgroup
of $\C_p^g$ has an accumulation point at $0$ -- unlike the image of
$\trop\circ\eta$.  There is an odd interplay here between the $p$-adic logarithm
$\Log$ and the Archimedean logarithm $\trop$.

\subsubsection{Berkovich--Coleman integration on a curve}
Let $\iota\colon X\inject J$ be the Abel--Jacobi map with respect to a base
point $P_0$.  By functoriality of the Berkovich--Coleman integral, for
a path $\gamma$ and a $1$-form
$\omega\in H^0(J,\Omega^1_{J/\C_p})\cong H^0(X,\Omega^1_{X/\C_p})$, we have
\[ \BCint_\gamma\iota^*\omega = \BCint_{\iota(\gamma)}\omega. \]
This integral can also be computed in terms of the universal cover
$\pi\colon\td X\to X^\an$.  Indeed, the map $\iota$ lifts in a unique way to a
morphism of universal covers $\td\iota\colon \td X\to T^\an$ such that
$\td\iota$ sends a fixed lift $\td P_0$ of $P_0$ to the identity.  Choosing any
lift $\td\gamma$ of $\gamma$, then, we have
\[ \BCint_{\td\gamma}\pi^*\iota^*\omega =
\BCint_{\td\gamma}\td\iota^*\pi^*\omega =
\BCint_{\td\iota(\td\gamma)}\pi^*\omega, \]
the latter integral being calculated on $T^\an$ as above.  We will use this
viewpoint later.

Again we emphasize that the Berkovich--Coleman integral can be computed locally
on $X^\an$, and by antidifferentiation on any domain where the differential is
exact. 

\subsection{Comparing the integrals}
The following result essentially states that the abelian and Berkovich--Coleman
integrals on $J$ differ only by the existence of $p$-adic periods for the
latter.  As above we let $X$ be a Mumford curve with Jacobian $J$ and
uniformization $\pi\colon T^\an\to J^\an$.

\begin{propsub}[{\cite[Proposition~3.16]{KatzRZB-uniform-bounds}}]\label{prop:comparing-integrals-jac}
  Let $\omega\in H^0(J,\Omega^1_{J/\C_p})$.  Then the homomorphism
  \begin{equation}\label{eq:int.diff.1}
    P\mapsto \BCint_0^P\pi^*\omega - \Abint_0^{\pi(P)}\omega
    \colon T(\C_p)\To\C_p
  \end{equation}
  factors through $\trop\colon T(\C_p)\cong(\C_p^\times)^g\to\Q^g$.
\end{propsub}

In other words, $\BCint$ and $\Abint$ coincide on $\trop\inv(0)\subset T(\C_p)$,
which we identify (via $\pi$) with its image in $J(\C_p)$.  Allowing $\omega$ to
vary, we can regard~\eqref{eq:int.diff.1} as a homomorphism
$T(\C_p)\to\Hom(H^0(J,\Omega^1_{J/\C_p}),\C_p) = \Lie(J)$, which again factors
through $\trop$.  It follows formally that there exists a $\Q$-linear map
$L\colon\Q_p^g\to\Lie(J)$ such that, for all $P\in T(\C_p)$ and
$\omega\in H^0(J,\Omega^1_{J/\C_p})$, we have
\begin{equation}\label{eq:int.diff.2}
  \BCint_0^P\pi^*\omega - \Abint_0^{\pi(P)}\omega
  = \angles{L\circ\trop(P),\,\omega},
\end{equation}
where $\angles{\scdot,\scdot}$ here denotes the duality pairing between
$\Lie(J)$ and $H^0(J,\Omega^1_{J/\C_p})$.

\subsubsection{The tropical Abel--Jacobi map}\label{sec:tropical-abel-jacobi}
\begin{figure}[ht]
  \centering
  \includegraphics{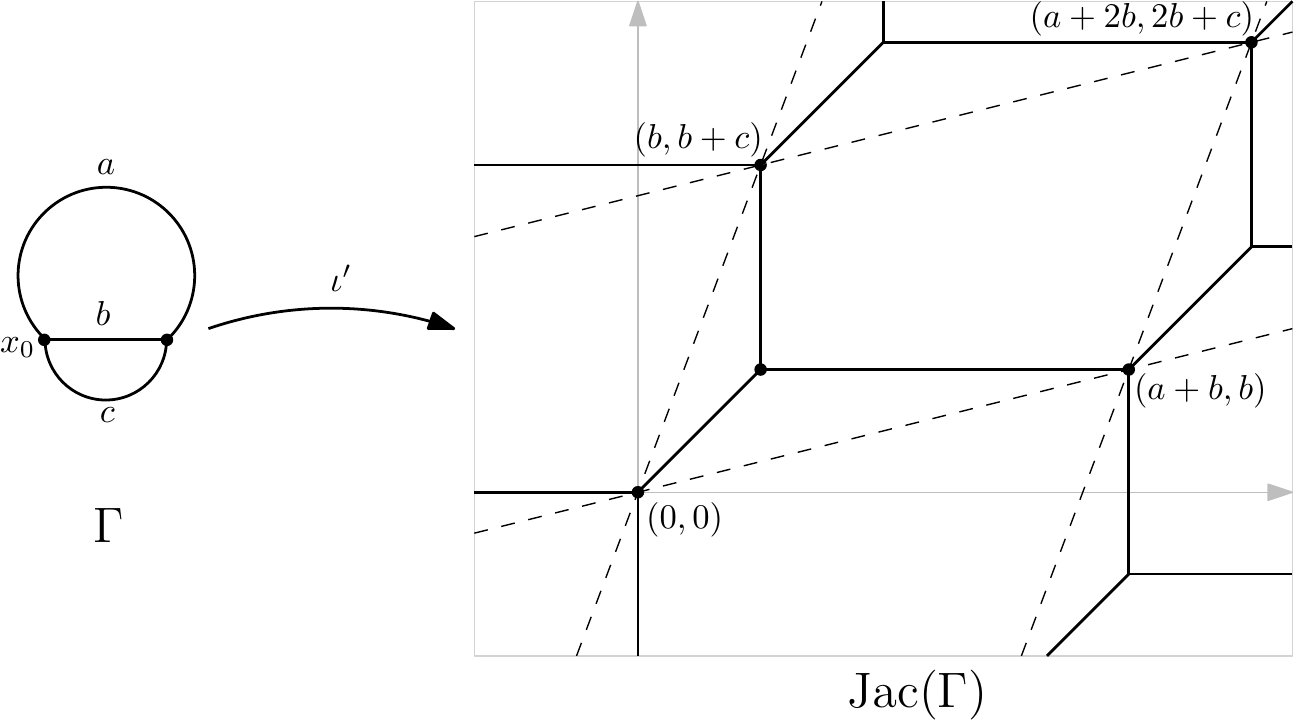}
  
  \caption{An illustration of the tropical Abel--Jacobi map $\iota'$
    of~\S\ref{sec:tropical-abel-jacobi}.  Here $a,b,c$ denote edge lengths, and
    $\Jac(\Gamma)$ is the quotient of $\R^2$ by the lattice generated by
    $(a+b,b)$ and $(b,b+c)$.  In this case, the tropical balancing condition of
    Theorem~\ref{thm:tropical.aj.is.nice} says that, at each vertex in the image
    of $\iota'$, the primitive integral vectors along the adjacent edges sum to zero.
    (Recreated from~\cite[Figure~7]{Mikhalkin2008}.)}
  \label{fig:tropical.aj}
\end{figure}

At this point we need to introduce the tropical Abel--Jacobi map, which controls
the difference between the integrals on $X^\an$.  For references,
see~\cite{Mikhalkin2008,BakerF:metric,bakerR:skeleton}.  Fix a skeleton $\Gamma$ of
$X$.  The \emph{edge length pairing} is the bilinear map
$[\scdot,\scdot]\colon C_1(\Gamma)\times C_1(\Gamma)\to\R$ on the group of
simplicial $1$-chains, defined on directed edges by
\[ [e,e] = \ell(e) \quad\text{and}\quad [e,e'] = 0, \]
where $\ell(e)$ is the length of $e$, and $e'\neq e$.  We restrict the pairing
to $H_1(\Gamma,\Z)\subset C_1(\Gamma)$.  As the edge length pairing is clearly
symmetric and non-degenerate, it induces a homomorphism
\[ \eta'\colon H_1(\Gamma,\Z)\To \Hom(H_1(\Gamma,\Z),\R) =
H^1(\Gamma,\R)\cong\R^g, \]
whose image is a lattice.  The \defi{Jacobian} of the metric graph $\Gamma$ is
by definition
\[ \Jac(\Gamma) \coloneq \R^g/\eta' H_1(\Gamma,\Z). \]
One proves as for algebraic curves that $\Jac(\Gamma)$ can be canonically
identified with $\Div^0(\Gamma)/\Prin(\Gamma)$, the group of degree-zero
divisors on $\Gamma$ modulo the divisors of tropical meromorphic functions
(see~\S\ref{sec:potent-theory-berk}).  The \defi{tropical Abel--Jacobi map} with
respect to a point $x_0\in\Gamma$ is the function
$\iota'\colon\Gamma\to\Jac(\Gamma)$ defined by $\iota'(x) = [x] - [x_0]$.  See Figure~\ref{fig:tropical.aj}.

Let $\td\iota'\colon\td\Gamma\to\R^g$ be a lift of $\iota'$ to universal covers.
This is a function from an infinite metric graph into a Euclidean space, whose
structure was studied by
Mikhalkin--Zharkov~\cite{Mikhalkin2008}.  Among many other
things, they prove:

\begin{thmsub}[Mikhalkin--Zharkov]\label{thm:tropical.aj.is.nice}
  Let $\td e\subset\td\Gamma$ be an edge, and let $e\subset\Gamma$ be its
  image.  
  \begin{enumerate}
  \item If $\Gamma\setminus e$ is disconnected, then $\td\iota'$ is constant on
    $\td e$.
  \item If $\Gamma\setminus e$ is connected, then $\td\iota'$ is affine-linear
    on $\td e$ with rational slope.
  \item $\td\iota'$ satisfies the tropical balancing condition at vertices.
  \end{enumerate}
\end{thmsub}

The balancing condition in the last part of
Theorem~\ref{thm:tropical.aj.is.nice} roughly says that at any vertex $x\in\td\Gamma$,
a weighted sum of the images of the tangent vectors at $x$ under $\td\iota'$ is
equal to zero.  This implies, for instance, that if $x$ has three adjacent edges
$\td e_1,\td e_2,\td e_3$, then their images under $\td\iota'$ are coplanar.
See the end of \S3
in~\cite{BakerF:metric} for details.  See also Figure~\ref{fig:tropical.aj}.

\subsubsection{Tropicalizing the Abel--Jacobi map}\label{sec:tropicalize-abel-jacobi}
The relationship between the algebraic and tropical Abel--Jacobi maps was
studied in~\cite{bakerR:skeleton}.  The relevant results are as follows.  Extend
$\trop\colon(\C_p^\times)^g\to\Q^g$ to a function
$\trop\colon T^\an=(\bG_m^\an)^g\to\R^g$ by the rule
\[ \trop(\|\scdot\|) = (-\log\|t_1\|,\ldots,-\log\|t_g\|), \]
where $t_i$ is a coordinate on the $i$th factor of $\bG_m$.  This map descends
to a function
\[ \bar\tau\colon J^\an = T^\an/\eta(H_1(X^\an,\Z))
\To\R^g/\trop\circ\eta (H_1(X^\an,\Z)). \]
The real torus on the right is called the \defi{skeleton} of $J$, and we denote it
by $\Sigma = \Sigma(J)$.  Berkovich~\cite[Chapter~6]{Berkovich:spectralTheory} showed that there is a natural embedding
$\Sigma\inject J^\an$, and that $J^\an$ deformation retracts onto $\Sigma$.

\begin{thmsub}[{\cite{bakerR:skeleton}}]
  Let $X$ be a Mumford curve with Jacobian $J$ and uniformization
  $\pi\colon T^\an\to J^\an$.
  \begin{enumerate}
  \item The lattice $\trop\circ\eta(H_1(X^\an,\Z))\subset\R^g$ coincides with
    the lattice $\eta'(H_1(\Gamma,\Z))$ induced by the edge length pairing.
    Hence $\Sigma = \Jac(\Gamma)$.
  \item Let $P_0\in X(\C_p)$ and let $x_0 = \tau(P_0)\in\Gamma$, where
    $\tau\colon X^\an\to\Gamma$ is the retraction map.  Let
    $\iota\colon X\inject J$ and $\iota'\colon\Gamma\to\Sigma=\Jac(\Gamma)$ be
    the corresponding Abel--Jacobi maps.  
    Then the following square is commutative:
    \begin{equation}\label{eq:aj.square}
      \xymatrix @=.2in{
        {X^\an} \ar[r]^\iota \ar[d]_\tau & {J^\an} \ar[d]^{\bar\tau} \\
        {\Gamma} \ar[r]_{\iota'} & {\Sigma}
      }
    \end{equation}

  \end{enumerate}
\end{thmsub}

The statement~(1) assumes that we have chosen compatible bases for
$T\cong\bG_m^g$ and $\R^g$.  In the basis-free version, the cocharacter lattice
of $T$ is $H^1(X^\an,\Z) = H^1(\Gamma,\Z)$, and both lattices
$\trop\circ\eta(H_1(X^\an,\Z))$ and $\eta'(H_1(\Gamma,\Z))$ live in
$H^1(X^\an,\R) = H^1(\Gamma,\R)$.

\subsubsection{Comparing the integrals on a curve}
Now we combine the results
of~\S\S\ref{sec:tropical-abel-jacobi}--~\ref{sec:tropicalize-abel-jacobi}.  We
mentioned in~\S\ref{sec:abel-integr-curve} that $\Abint_P^Q$ is computed by
formal antidifferentiation for $P,Q$ contained in the same open ball in $X^\an$.
As the same is true for $\BCint_P^Q$, it follows that $\BCint = \Abint$ on open
balls contained in $X^\an$.

By Theorem~\ref{thm:formal.fibers}, the $\C_p$-points of $X$ can be partitioned
into open balls and open annuli.  Comparing the integrals on annuli is
more subtle.  It is clear from the way that universal covers are constructed
that the universal cover $\td X$ of $X$ deformation retracts onto the universal
cover $\td\Gamma$ of its skeleton $\Gamma$.  Hence we can
lift~\eqref{eq:aj.square} to universal covers:
\begin{equation}\label{eq:aj.square.covers}
\xymatrix @=.2in{
  {\td X^\an} \ar[r]^{\td\iota} \ar[d]_\tau & {T^\an} \ar[d]^{\trop} \\
  {\td\Gamma} \ar[r]_{\td\iota'} & {\R^g} }
\end{equation}

Let $\td e\subset\td\Gamma$ be an open edge, with image $e\subset\Gamma$.  Let
$A=\tau\inv(e)$ and $\td A = \tau\inv(\td e)$, and note that $\td A\isom A$ are
open annuli.  Suppose that our base points $P_0$ and $x_0 = \tau(P_0)$ are
contained in $A$.
It follows from~\eqref{eq:int.diff.2} that there is a linear map
$L\colon\Q^g\to\Lie(J)$ such that, for all $\omega\in H^0(X,\Omega^1_{X/\C_p})$
and $P\in A(\C_p) = \td A(\C_p)$, we have
\[ \BCint_{P_0}^P\iota^*\omega - \Abint_{P_0}^{P}\iota^*\omega =
\angles{L\circ\trop\circ\td\iota(P),\,\omega},
\]
(Note that $A$ is simply connected, as it deformation retracts onto $e$, so that
$\BCint_{P_0}^P$ makes sense.)  Lifting to universal covers and using
commutativity of~\eqref{eq:aj.square.covers}, we have
$\trop\circ\td\iota = \td\iota'\circ\tau$.  But by
Theorem~\ref{thm:tropical.aj.is.nice}, $\td\iota'$ is \emph{affine-linear} on
$\td e$.  Recalling from~\S\ref{sec:skeleton-an-annulus} that
$\tau\colon A(\C_p)\to e$ is simply the valuation map $P\mapsto v(P)$ after
choosing an isomorphism $A\cong\bS(\rho)_+$, we have
derived the following important result of
Stoll~\cite[Proposition~7.3]{stoll:uniform}.

\begin{propsub}[Stoll]\label{prop:comp.int.annuli}
  With the above notation, there is a $\C_p$-linear map
  $a\colon H^0(X,\Omega^1_{X/\C_p})\to\C_p$ such that, 
  for all $P,Q\in A(\C_p)$, we have
  \[ \BCint_P^Q\omega - \Abint_P^Q\omega = a(\omega)\big(v(Q)-v(P)\big). \]
\end{propsub}

\begin{corsub}\label{cor:integrals.equal.annuli}
  Let $V$ be the subspace of $H^0(X,\Omega^1_{X/\C_p})$ consisting of all
  $\omega$ such that $\BCint_P^Q\omega = \Abint_P^Q\omega$ for all
  $P,Q\in A(\C_p)$.  Then $V$ has codimension at most one.
\end{corsub}

Corollary~\ref{cor:integrals.equal.annuli} is very important, because it
produces a single linear condition on $\omega$ for the Berkovich--Coleman
integral to coincide with the abelian integral on $A$.  As the former is
computed by formal antidifferentiation, and the latter can be chosen to vanish on rational or
torsion points, this is crucial to any application of the Chabauty--Coleman
method to annuli.

Using the tropical balancing condition in Theorem~\ref{thm:tropical.aj.is.nice},
one can extend Corollary~\ref{cor:integrals.equal.annuli} to the following more
general situation, which is important for applications 
to uniform Manin--Mumford.
Let $x\in\Gamma$ be a vertex, let $S_x$ be a star neighborhood, and let
$U_x = \tau\inv(S_x)$ be a basic wide open subdomain, as
in~\S\ref{sec:decomp-into-wide}.

\begin{corsub}\label{cor:integrals.equal.wideopens}
  Let $d$ be the valency of $x$ in $\Gamma$.  Let $V$ be the subspace of
  $H^0(X,\Omega^1_{X/\C_p})$ consisting of all $\omega$ such that
  $\BCint_P^Q\omega = \Abint_P^Q\omega$ for all $P,Q\in U_x(\C_p)$.  Then $V$
  has codimension at most $d-1$.
\end{corsub}

\begin{rem*}
  For simplicity of exposition, in~\S\ref{sec:uniformity-results} we have omitted the result that requires
  Corollary~\ref{cor:integrals.equal.wideopens}, a variant of
  Theorem~\ref{thm:geom.torsion.bound} below with the quadratic factor replaced
  by a linear factor in $g$ --
  see~\cite[Theorem~5.5(2)]{KatzRZB-uniform-bounds}.  We have included
  Corollary~\ref{cor:integrals.equal.wideopens} for its conceptual importance.
\end{rem*}

\begin{rem*}
  The conclusions of Corollaries~\ref{cor:integrals.equal.annuli}
  and~\ref{cor:integrals.equal.wideopens} still hold true when $J$ does not have
  totally degenerate reduction, i.e., when $X$ is not a Mumford curve.  The
  argument is more technical, however, as it involves the general uniformization
  theory of non-Archimedean abelian varieties.
  See~\cite[\S4]{KatzRZB-uniform-bounds} for details.
\end{rem*}

\subsection{The Stoll decomposition of a non-Archimedean curve}
\label{sec:stoll-decomposition}

  Stoll had the remarkable idea to cover $X(\Q)$ by sets bigger than tubes
  around $\F_p$-points of a regular model.  This offers advantages over the
  usual effective Chabauty arguments which concede that there may be a rational point 
  in each tube because there may be arbitrarily many residue points in the bad reduction
  case.  Being able to integrate on these larger sets necessitated the use of a more involved
  $p$-adic integration which we covered above.
  To obtain a uniform bound, one must pick an economical covering.  Here, we outline
  Stoll's choice of covering which comes from a minimal regular model.  We will employ this covering
  in the proofs of Theorem~\ref{T:KRZB} and Theorem~\ref{T:uniformity-K-torsion}.  We will use
  a different covering, one coming from a semistable model in the proof of 
  Theorem~\ref{T:uniformity-torsion}.  We will state the results in this section for $\Q_p$ with the
  understanding that analogous results hold for its finite extensions.

  The properties of the covering are summarized as follows:
  \begin{proposition}[Stoll, \cite{stoll:uniform}, Proposition~5.3 ] There exists
  $t\in\{0,1,2,\ldots,g\}$ such that $X(\Q_p)$ is covered by at most
  $(5q+2)(g-1)-3q(t-1)$ embedded open balls and at most embedded $2g-3+t$ open
  annuli, all defined over $\Q_p$. 
  \end{proposition}
  
  Let $\cX$ be a minimal regular model of $X$ over $\Z_p$.  Denote the components of the 
  special fiber $\cX_{\F_p}$ by $X_v$ and write their multiplicity as $m(X_v)$.  We observe that 
  a point of $X(\Q_p)$ must specialize to a
  component $X_v$ of multiplicity $1$.  Indeed, such a point extends to a section 
  $\sigma\colon\Spec(\Z_p)\to \cX$.  This section must intersect the special fiber with multiplicity 
  $1$.  Therefore, we only need to find a collection of subsets of $X^{\an}$ that contain all 
  points  specializing to smooth points of $\cX_{\F_p}$ on components of multiplicity $1$.
  
  Let $K_{\cX}$ denote the relative canonical bundle of $\cX$. The adjunction formula for
  surfaces states that
  \[K_{\cX}\cdot X_v=2p_a(X_v)-2-X_v^2\]
  where $p_a(X_v)$ is the arithmetic genus of $X_v$.
  In the case that $g(X)\geq 2$ and $\cX$ is a regular minimal model, then the $K$-degree
  of each component is nonnegative.  
   The total $K$-degree is
  \[2g(X)-2=\sum_v m(X_v) K_{\cX}\cdot X_v\]
  so there are at most $2g(X)-2$ components of positive $K_{\cX}$-degree.
 The intersection pairing on the components of the special fiber is negative semidefinite and $X_v^2=0$ 
  if and only $X_v$ is the only component of $\cX_{\F_p}$.  In the case that $\cX_v$ has more than
  one component,  then $K\cdot X_v=0$ if and only if $p_a(X_v)=0$ and $X_v^2=-2$.  Such curves
  of $K$-degree equal to $0$ are called $(-2)$-curves.  There are two possibilities for the 
  $(-2)$-curves:
  \begin{enumerate}
    \item They are part of a chain of $(-2)$-curves which meet distinct multiplicity $1$ components;
    \item They are so-called $\A^1$-components which meet a multiplicity $2$ component in 
    one point, meet several components in a single point, or meet a component in a multiplicity 
    $2$ intersection point.
  \end{enumerate}
  The number of $\A^1$-components and chains of $(-2)$-curves of multiplicity
  $1$ can be bounded by a combinatorial study of the arithmetic graph encoding
  the components of the special fiber and their intersections. Now, points
  specializing to smooth points of chains can be covered by open annuli.  This
  is essentially because blowing down chains of $\bP^1$'s yields a node, whose
  inverse image is an annulus, as in Theorem~\ref{thm:formal.fibers}(3).  Points
  specializing to the $\F_p$-points of $\A^1$- components and multiplicity $1$
  components of positive $K$-degree can be covered by tubes
  around $\F_p$-points.  The Hasse--Weil bound gives an upper bound on the number
  of such tubes.
  
\section{Uniformity results}
\label{sec:uniformity-results}

Let $X$ be a smooth, proper, geometrically connected curve of genus $g$ over
$\Q$.  All arguments hold equally well over any number field; we restrict
ourselves to the rationals for concreteness.  Let $J$ be the Jacobian of $X$.
Fix a prime number $p$.  Choose an Abel--Jacobi map
$\iota\colon X_{\C_p}\inject J_{\C_p}$.  

We would like to give bounds on $\#X(\Q)$ and on the size of a torsion packet
$\#\iota\inv(J(\C_p)_{\tors})$ that depend only on $p$ and $g$.  Let $S$ be either
one of these sets.  Our strategy is as follows:
\begin{enumerate}
\item Find nonzero $\omega\in H^0(X_{\C_p},\Omega^1_{X_{\C_p}/\C_p})$ such that
  $\Abint_P^Q\omega = 0$ for all $P,Q\in S$.
  
  \smallskip
  For uniform Mordell when the rank condition is satisfied, we find such $\omega$ using the classical
  Chabauty--Coleman argument, i.e., by taking the closure of $J(\Q)$ in
  $J(\Q_p)$.  For uniform Manin--Mumford, any nonzero $\omega$ works, since
  $\Abint$ vanishes on torsion points of $J$.

  \medskip

\item Decompose $X^\an$ into basic wide open domains $U$ such that for
  each $U$, there exists $\omega$ as in~(1) such that
  \newcommand\condE{\textbf{(E)}}
  \newcommand\condI{\textbf{(I)}}
  \begin{itemize}
    \setlength{\itemindent}{.3in}
  \item[\condE] $\omega = df$ is exact on $U$.
  \item[\condI] $\BCint_P^Q\omega = \Abint_P^Q\omega$ for all $P,Q\in U(\C_p)$.
  \end{itemize}

  \smallskip These conditions guarantee that an antiderivative $f$ vanishes on
  $S\cap U(\C_p)$.  They each impose some number of linear conditions on the
  space of suitable $\omega\in H^0(X_{\C_p},\Omega^1_{X_{\C_p}/\C_p})$.
  Indeed,~\condI\ imposes $d-1$ linear conditions by
  Corollary~\ref{cor:integrals.equal.wideopens}, where $d$ is the valency of the
  vertex in $U$.  The number of conditions imposed by~\condE\ is simply the
  dimension of the de Rham cohomology group $H^1_{\dR}(U)$, which was computed
  by Coleman and is essentially the same as for a complex ``pair of pants'' with $d$
  holes.

  \medskip

\item Use Lemma~\ref{lem:slope.bound} to bound the slopes of
  $G = -\log\|\omega\|$ on the skeleton of each $U$.

  \medskip

\item Use the bound in~(3) to bound the slopes of the antiderivative, i.e., of
  $F = -\log|f|$.

  \smallskip
  This is a replacement for the ``$p$-adic Rolle's theorem'' part of the usual
  Chabauty--Coleman argument, i.e., bounding the number of zeros of the
  antiderivative of $\omega$ in terms of the number of zeros of $\omega$.

\medskip

\item Use Corollary~\ref{cor:zeros.from.slopes} to bound the number of zeros of
  $f$ on $U$.
\end{enumerate}
By construction, all points of $S\cap U(\C_p)$ are zeros of $f$, so as long as
we can bound the number of such $U$ uniformly in terms of $g$, then we will be
done.  For applications to uniform Mordell, this was already done by Stoll, as
explained in~\S\ref{sec:stoll-decomposition}.  For applications to uniform
Manin--Mumford, this is an elementary argument in the combinatorics of stable
genus-$g$ graphs, as we will explain below.  

We may modify step~(2) to pick a number of functions on each $U$ such that each 
rational or torsion point is a zero of one of the functions.  This circumvents some of
the conditions that~\textbf{(I)} imposes. 

Essentially all of the above steps depend only on machinery which we have
already explained, except for step~(4).

\subsection{$p$-adic Rolle's theorem for annuli and basic wide opens}
In this subsection, we let $X$ denote a smooth, proper, connected curve over $\C_p$.
Let $\omega\in H^0(X,\Omega^1_{X/\C_p})$ be a nonzero $1$-form and let
$G = -\log\|\omega\|$, the valuation with respect to the model metric
$\|\scdot\|$ defined in~\S\ref{sec:divis-regul-diff}.  Let $x$ be a vertex of a
skeleton $\Gamma$ of $X$, and let $U = U_x$ be a basic wide open 
subdomain centered a $x$, as in~\S\ref{sec:decomp-into-wide}.  Suppose that
$\omega = df$ is exact on $U$, where $f$ is an analytic function on $U$, and let
$F = \log|f|\colon U\to\R$.  We wish to bound the number of zeros of $f$ on $U$,
which, thanks to Corollary~\ref{cor:zeros.from.slopes}, is equivalent to
bounding the slopes of $F$ on the edges adjacent to $x$.  If $e$ is an open edge
adjacent to $x$ contained in $U$, then its inverse image under retraction
$A = \tau\inv(e)$ is an open annulus, isomorphic to $\bS(\rho)_+$ for some
$\rho$.  Since $A$ is also a basic wide open subdomain (centered at any point on
its skeleton), we may as well assume $U = A = \bS(\rho)_+$ is an annulus.

\begin{remsub}
  The analogous, simpler problem of bounding the number of zeros of $f$ in terms
  of those of $\omega$ on an open ball was discussed
  in~\S\ref{ss:local-bounds}.  (This result is much more similar to the
  classical Rolle's theorem.)
\end{remsub}

Let $t$ be a parameter on $\bS(\rho)_+$.  That is, choose a parameter $t$ on
$\bA^1 = \Spec(\C_p[t])$, so that
$\bS(\rho)_+ = \{x\mid|t(x)|\in(\rho,1)\}$.  Let $r = -\log\rho$, so the
skeleton of $\bS(\rho)_+$ is $(0,r)$.  The canonical line bundle on
$\bS(\rho)_+$ is trivial; a nonvanishing section is the invariant differential
$dt/t$.  Hence we can write $\omega = h(t)dt/t$, where $h$ is an analytic
function on $\bS(\rho)_+$.  As explained in Example~\ref{eg:annulus}, such an
analytic function has the form
\[ h = \sum_{i=-\infty}^\infty a_it^i, \]
where $|a_i|\tau^i\to 0$ as $t\to\pm\infty$ for all $\tau\in(\rho,1)$.
Integrating term by term, we have
\[ f = \sum_{i\neq 0} \frac{a_i}it^i + C. \]

Unfortunately, $v(a_i/i) < v(a_i)$ when $p\mid i$.  This means that the Newton
polygon of $f$ looks different from the Newton polygon of $h$.  In fact, it is
not hard to see that, although $h$ has finitely many zeros, and hence finitely
many slopes in its Newton polygon, in general $f$ may have infinitely many such
slopes, hence \emph{infinitely many zeros} in $\bS(\rho)_+$.  However, since $f$
is analytic on any closed sub-annulus $A'\subset\bS(\rho)_+$, it can have only
finitely many zeros on $A'$.  Translating back into slopes, it turns out that
one can indeed bound the slopes of $F$ in terms of those of $G$, but only on a
smaller sub-annulus.  In other words, some \emph{overconvergence} condition is
necessary.

To make the bound precise, we introduce the following quantity.

\begin{defnsub}\label{def:Np}
  Let $r$ be a positive real number, let $N_0$ be an integer, and let $p$ be a
  prime.  Define $N_p(r,N_0)$ to be the smallest positive integer $N$ such that
  for all $n\geq N$, one has
  \begin{equation}\label{eq:one.annulus.bound}
    r(n - N_0) > \lfloor\log_p(n)\rfloor.
  \end{equation}
\end{defnsub}

\begin{remsub}
\label{R:bounds-on-NpAB}
The integer $N_p(r,N_0)$ gets larger as $N_0$ increases and as $r$ decreases,
and gets smaller as $p$ increases.  If $N_0\geq 0$ and $p \geq N_0 + 2$ then
$N_p(r,N_0) = N_0 + 1$ because $\lfloor\log_p(N_0+1)\rfloor = 0 < r$.  One
should think of $N_p(r,N_0)-N_0$ as the correction to the $p$-adic Rolle's theorem
coming from the fact that $1/p$ has negative valuation.  See
also~\cite[\S6]{Stoll:2006-independence-of-rational-points} for a more sophisticated approach to
the same problem. (Stoll's correction factor $\delta(\scdot,\scdot)$ is slightly
better, but ours is easier to define.)
\end{remsub}

\begin{figure}[ht]
  \centering
  \begin{tikzpicture}[every circle/.style={fill=black, radius=.7mm},
    every node/.style={font=\tiny}]
    \coordinate (ZM) at (0,0);
    \coordinate (ZP) at (4,0);
    \coordinate (R) at (2.5, 0);
    \draw (ZM) node[left] {$y$} -- (ZP) node[right] {$x$};
    \draw[->, very thick] (ZP) -- ++(-.5cm,0)
          node [pos=0.8, below] {$\bv_x$};
    \fill (ZM) circle[] (ZP) circle[] (R) circle;
    \path (R) ++(1mm,0) node [below] {$z$};
    \draw[->, very thick] (R) -- ++(-.5cm,0) node[pos=.8, below] {$\bv$};
    \draw[|<->|, thin] (0,5mm) -- (2.5,5mm) node[pos=.5, above] {$a$};
    \draw[|<->|, thin] (2.5,5mm) -- (4,5mm) node[pos=.5, above] {$r-a$};
  \end{tikzpicture}
  
  \caption{Illustration of the notation used in
    Proposition~\ref{p:annularbounds}.  The interval represents the edge $e$,
    which has length $r$ and endpoints $x,y$.}
  \label{fig:annulus.setup}
\end{figure}
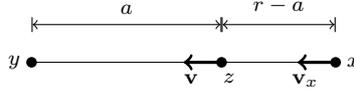

The rest of the notation in Proposition~\ref{p:annularbounds} is defined in Figure~\ref{fig:annulus.setup}.

\begin{propsub}\label{p:annularbounds} 
  Let $\omega\in H^0(X,\Omega^1_{X/\C_p})$ be a nonzero global differential, and
  suppose that $\omega$ is~\emph{exact} on $\bS(\rho)_+$, so $\omega = df$ for
  an analytic function $f$ on $\bS(\rho)_+$.  Let $F = -\log|f|$ and
  $G = -\log\|\omega\|$, and let $N_0 = d_{\bv_x}F_0(x)$.  Choose
  $a\in(0,r)$, and let $\bv$ be the tangent direction at $z$ in the
  direction of $y$.  Then $d_{\bv} F(z)\leq N_p(r-a,N_0)$.
\end{propsub}

The proof is a straightforward but tedious argument involving Newton polygons,
and can be found in \cite[Proposition~4.7]{KatzRZB-uniform-bounds}.  Combining
Proposition~\ref{p:annularbounds} with the slope formula in the form of
Corollary~\ref{cor:zeros.from.slopes}, the statement about the combinatorics of
sections of the tropical canonical bundle in Lemma~\ref{lem:slope.bound}, and
the fact that $G=-\log\|\omega\|$ is in fact a section of the tropical
canonical bundle as in Corollary~\ref{cor:effectivity}, we can bound the number
of zeros of the antiderivative of a global $1$-form which is exact on a wide open
subdomain $U$ in terms of the valency of the central vertex, the genus $g$, and our
fixed prime $p$ -- after shrinking $U$ a little bit.

To make this precise, let $U$ be a basic wide open subdomain in $X^\an$, so $U$
is the inverse image under retraction of a neighborhood of a point $x$ in a
skeleton $\Gamma$.  Let $e_1,\ldots,e_d$ be the open edges adjacent to $x$ in
$U$, where $d$ is the valency of $x$.  For
$a < \min\{\ell(e_1),\ldots,\ell(e_d)\}$ let $U_a$ denote the basic wide open
subdomain of $U$ obtained by deleting an open annulus of modulus $\exp(-a)$, the
inverse image under retraction of an interval of length $a$ at the end of each
$e_i$.  In other words, $U_a$ is obtained by shortening each edge of its
skeleton by $a$ units.

\begin{thmsub}\label{thm:wideopen.bounds}
  With the above notation, let $\omega\in H^0(X,\Omega^1_{X/\C_p})$ be a nonzero
  global differential, and suppose that $\omega$ is exact on $U$, so
  $\omega = df$ for an analytic function $f$ on $U$.  Then $f$ has at most
  $d\,N_p(a,\,2g-1)$ geometric zeros, counted with multiplicity, on
  $U_a$.  If $U$ is defined with respect to a star neighborhood in a
  skeleton with no genus-zero leaves, then we may replace $2g-1$ by $2g-2$.

  In particular, if $U$ is an open annulus, then $f$ has at most
  $2\,N_p(a,2g-1)$ zeros on $U_a$.
\end{thmsub}

\subsection{Uniform bounds on rational points}\label{sec:unif-bounds-rati}
We are now in a position to put all of the ingredients together to prove our
uniform Mordell theorem, subject to a restriction on the rank of the Jacobian.
As in the beginning of this section, let $X$ be a smooth, proper, geometrically
connected curve of genus $g$ over $\Q$, let $J$ be the Jacobian of $X$, and fix
a prime number $p\geq 3$.  Suppose that $X$ has a rational point, so that we can
choose an Abel--Jacobi map $\iota\colon X\inject J$ defined over $\Q$.

\begin{thmsub}\label{thm:with.ss.model}
  Let $G\subset J(\Q)$ be a subgroup of rank at most $g-3$.  Then
  \begin{equation}\label{eq:rational.bound}
    \#\iota^{-1}(G) \leq \big(5pg + 6g - 2p - 8\big)\,(4g-2).
  \end{equation}
\end{thmsub}

Before sketching the proof, we state the two main consequences we have in mind.
The first is an application of Theorem~\ref{thm:with.ss.model} with $G = J(\Q)$.

\begin{corsub}\label{cor:rational.points}
  If $J(\Q)$ has rank at most $g-3$, then
  \[ \#X(\Q)\leq \big(5pg + 6g - 2p - 8\big)\,(4g-2). \]
\end{corsub}

In particular, taking $p=3$ yields $\#X(\Q)\leq 84g^2 - 98g + 28$.  This extends
the main result of~\cite{stoll:uniform} to the case of arbitrary curves with
small Mordell--Weil rank.

In the next corollary, we take $G = J(\Q)_{\tors}$, which has rank zero.

\begin{corsub}\label{cor:rational.torsion.packets}
  Any \emph{rational} torsion packet has size at most 
  \[ \#\iota\inv(J(\Q)_{\tors}) \leq\big(5pg + 6g - 2p - 8\big)\,(4g-2). \]
\end{corsub}

Although the conclusion of Corollary~\ref{cor:rational.torsion.packets} is
much weaker than a full uniform Manin--Mumford theorem, it has the advantage that it
is \emph{completely unconditional}: it applies to every $X$.

\begin{proof}[Sketch of the proof of Theorem~\ref{thm:with.ss.model}]
  First we use the results of~\S\ref{sec:stoll-decomposition} to decompose
  $X(\Q_p)$ into some number of open annuli and open balls, all defined over
  $\Q_p$.  The number of annuli and balls is uniformly bounded in $g$, so it
  suffices to uniformly bound the number of $\Q$-points contained in each.

  By the usual Chabauty argument, the space $V$ of $1$-forms
  $\omega\in H^0(J_{\C_p},\Omega^1_{J_{\C_p}/\C_p})$ vanishing on $G$ has
  codimension at most equal to the rank of $G$, hence dimension at least $3$.
  For each open ball $B$, suppose that $B$ contains a $\Q$-rational point.  Then
  we apply the standard Chabauty--Coleman argument, as strengthened e.g., by
  Stoll~\cite{Stoll:2006-independence-of-rational-points} when $p$ is small
  relative to $g$, to any nonzero $\omega\in V$ to bound the size of
  $B\cap X(\Q)$.

  Now let $A\cong\bS(\rho)_+$ be an open annulus in the Stoll decomposition.
  Since $A$ is defined over $\Q_p$, its modulus must be an element of
  $|\Q_p^\times|$, so $A\cong\bS(|p|^b)$ for some $b\geq 1$.  Hence
  \[ A(\Q_p) = \big\{x\in\Q_p^\times\mid v(x)\in(0,b)\cap\Z\big\}
  = \big\{x\in\Q_p^\times\mid v(x)\in[1,b-1]\cap\Z\big\}.  \]
  Suppose that
  $A\cap X(\Q)\neq\emptyset$.  This implies $b\geq 2$.  By
  Corollary~\ref{cor:integrals.equal.wideopens}, the dimension of the space of
  $\omega\in V$ such that $\BCint\omega = \Abint\omega$ on $A$ is at least $2$.
  A $1$-form $\omega = \sum a_it^i\,dt/t$ is exact if and only if $a_0 = 0$;
  hence the dimension of the space of $\omega\in V$ such that
  $\BCint\omega=\Abint\omega$ on $A$ and such that $\omega = df$ is exact on $A$
  is at least $1$.  In particular, there exists a nonzero $\omega$ satisfying
  all of these properties.  Choosing the antiderivative $f$ to vanish on one
  point of $A\cap X(\Q)$, it will also vanish on all other points of
  $A\cap X(\Q)$.  Now we apply Theorem~\ref{thm:wideopen.bounds} to $f$ with
  $a\in(0,1)$ and $U=A$.  Because each point of $A\cap X(\Q)$ is contained in 
  $U_r\cap X(\Q)$, we conclude that $\#A\cap X(\Q)\leq 2N_p(a,2g-1)$.  Taking the
  limit as $a\to 1$ yields $\#A\cap X(\Q)\leq 2N_p(1,2g-1)$.  Finally, one shows
  that for $p\geq 3$, the inequality $N_p(1,N_0)\leq 2N_0$ always holds; this
  allows us to avoid having a term of $N_p(\scdot,\scdot)$
  in~\eqref{eq:rational.bound}.
\end{proof}

\subsection{Uniform bounds on torsion packets}\label{sec:unif-bounds-tors}
Now we explain our progress towards a uniform Manin--Mumford theorem.
Suppose that the genus $g$ is at least $2$.  Then $X_{\C_p}$ admits a stable
model; it follows from this that $X_{\C_p}$ admits a skeleton $\Gamma$ which is a
stable graph, in the sense that all vertices of positive weight have valency at
least $3$.  For a vertex $x\in\Gamma$ we let $\deg(x)$ denote its valency and
$g(x)$ its weight.  We set
\[ E(g,p) \coloneq
\begin{cases}
  \#\GSp_{2g}(\F_5) &\text{if } p\neq 5 \\
  \#\GSp_{2g}(\F_7) &\text{if } p = 5.
\end{cases}\]

\begin{thmsub}\label{thm:geom.torsion.bound}
  Let $\Gamma$ be the minimal skeleton of $X_{\C_p}$, considered as a
  vertex-weighted metric graph.  If $g > 2g(x) + \deg(x)$ for all vertices $x$
  of $\Gamma$, then
  \begin{equation}\label{eq:torsion.bound}
    \#\iota^{-1}(J(\C_p)_{\tors}) \leq 
    (16g^2-12g)\,N_p\big((4E(g,p))^{-1},\,2g-2\big).
  \end{equation}
\end{thmsub}

\begin{remsub}
  The hypotheses of Theorem~\ref{thm:geom.torsion.bound} are satisfied if $X$
  is a Mumford curve of genus $g$ and all vertices of $\Gamma$ have valency at
  most $g-1$, e.g., if $g\geq 4$ and $\Gamma$ is trivalent.  
\end{remsub}

\begin{proof}[Sketch of the proof of Theorem~\ref{thm:geom.torsion.bound}]
  We will apply Theorem~\ref{thm:wideopen.bounds} to bound the number of zeros
  of a function $f$ on a subdomain $U_a$ contained in a basic wide open $U$.  Recall 
  that $U_a$ is obtained from $U$ by deleting open annular ends.  We must choose the
  moduli of the annular ends so that $U_a$ will be a member of a collection of subsets of 
  $X^{\an}$ whose union contains all the geometric torsion points.  We cannot choose $a$
  close to $1$ as in the proof of Theorem~\ref{thm:with.ss.model} because we
  may not cover every geometric torsion point.  Instead, we will cover $X^{\an}$ with
  basic wide opens indexed by flags of vertices of edges $(v,\epsilon)$ in $\Gamma$.  
  The basic wide opens $U_{x,\epsilon,r}$ will be of the form $(U_x)_r$ where $U_x$ is as 
  in subsection~\ref{sec:decomp-into-wide} and $r$ is less than half the length of the edges
  adjacent to $v$.  Such a collection will cover $X^{\an}$.  The reason for having 
  elements of the collection indexed by edges $\epsilon$ is to make a choice of $1$-form 
  $\omega_{(v,\epsilon)}=df_{(x,\epsilon)}$ exact on $U_{x,\epsilon,r}$ such that the 
  Berkovich--Coleman and abelian integrals agree on the annular end of $U_x$ 
  corresponding to $\epsilon$.  Therefore, every geometric torsion point of $(U_x)_r$ is
  guaranteed to be a geometric zero of $f_{(x,\epsilon)}$ for some $\epsilon$ adjacent to $v$.
  
  We begin by bounding the lengths of the edges of $\Gamma$. Unlike in the proof of 
  Theorem~\ref{thm:with.ss.model}, we cannot assume that an annulus obtained as the 
  inverse image of an edge in $\Gamma$ is defined over $\Q_p$.  These annuli are 
  defined over a finite extension $K$ of $\Q_p$ over which the stable model of $X$ is 
  defined.  Using results of Deligne--Mumford~\cite[Theorem~2.4]{DeligneMumford:irreducibility} 
  and Silverberg--Zarhin~\cite[Corollary~6.3]{silverbergZarhin:semistableReduction},
  one can bound $[K:\Q_p]\leq E(g,p)$.  See the proof of Theorem~5.5
  in~\cite{KatzRZB-uniform-bounds} for details.  Since the modulus of an annulus
  defined over $K$ is an element of $|K^\times|$, and since the ramification
  degree of $K/\Q_p$ is at most $E(g,p)$, it follows that each edge has length at least
  $1/E(g,p)$.  This explains the presence of $(4E(g,p))^{-1}$
  in~\eqref{eq:torsion.bound}.

  We discuss the choice of $\omega_{(v,\epsilon)}$ on $U_{v,\epsilon,r}$.  We have a 
  $g$-dimensional space $\Omega^1(X)$ of regular $1$-forms on $X$.  By 
  \cite{ColemanReciprocity}, $U_x$ is a basic wide open in the analytification of a 
  genus $g(x)$ curve $Y$ of good reduction.  Because it has  $\deg(x)$ ends, the de 
  Rham cohomology of $U_x$ is described by the following 
  exact sequence by a result of Coleman \cite{ColemanReciprocity}:
    \[ 0 \To H^1_{\dR}(Y)^{\alg} \To H^1_{\dR}(U) \xrightarrow{\Dsum\Res}
  \Dsum_{i=1}^{\deg(x)} \C_p \overset{\sum}\To \C_p \To 0 \]
  where $\Res$ is the residue map at the ends.  It follows that $H^1_{\dR}(U)$ is a 
  $(2g(x)+\deg(x)-1)$-dimensional vector space and the composition 
  \[\Omega^1(X) \hookrightarrow H^1_{\dR}(X)^{\alg} \To H^1_{\dR}(U)\]
  has a kernel of dimension at least $2$.  Consequently, by 
  Corollary~\ref{cor:integrals.equal.wideopens}, we may pick a $1$-form $\omega_{v,\epsilon}$
  on the kernel such that
  \[\BCint_P^Q\omega = \Abint_P^Q\omega\]
  for $P,Q\in X(\C_p)$ reducing to points of $v\cup \epsilon$. 

   By a straightforward argument using the combinatorics of semistable graphs, we
   see that each vertex of $\Gamma$ has valency at most $2g$, so that the number of 
   zeros of $f_{(v,\epsilon)}$ on $U_{v,\epsilon,r}$ is at most $2gN_p\big((4E(g,p))^{-1},\,2g-2\big)$.  
   The number of flags $(v,\epsilon)$ is equal to twice the number of edges which is at most $4g-3$,
   giving the stated bound.   
\end{proof}

\section{Other directions}
\label{sec:other-directions}

In this section, we will discuss other methods for producing bounds on the number of rational and torsion points.  The general theme here will be to find a set that is guaranteed to contain the rational and torsion points and bounding that set.  We will only discuss bounds on the total number of points rather than on their heights.  

\subsection{Buium's jet space method}

The method of Buium~\cite{buium:p-adic-jets} gives uniform bounds on the torsion points of curves in the case of good reduction at a prime $p$.  Here we will follow the exposition of Poonen~\cite{poonen:computing-torsion}.  The method we describe is a reduction mod $p$ of the $p$-adic analogue of a method used by Buium to address Lang's conjecture in characteristic $0$.  His method is inspired by an analogy between $p$-adic fields and function fields.  Specifically, in the characteristic $0$ setting,  a $k[\varepsilon]/\varepsilon^{n+1}$-point of a variety $X$ can be viewed as an $n$th order jet.  Consequently, one may use Weil restriction of scalars to produce a jet scheme $X^n$ satisfying $X^n(k)=X(k[\varepsilon]/\varepsilon^{n+1})$.  Similarly, given a variety over $\Z_p$, we may view  a $\Z/p^{n+1}$-point as a sort of $n$th order jet.  Here, we will use the more classical language of Greenberg transforms instead of Buium's language of $p$-jets.

Buium's result is the following:
\begin{theorem} Let $\iota\colon X\to J$ be an Abel--Jacobi map from a smooth curve of genus $g\geq 2$ defined over a number field $K$.  Let $\fp$  be a prime of $K$ with $p=\operatorname{char} \fp>2g$.  Assume that $K/\Q$ is unramified at $\fp$ and that $X$ has good reduction at $\fp$.  Then
\[\#\iota^{-1}\big(J(\overline{\Q})_{\tors}\big) \leq p^{4g}3^g(p(2g-2)+6g)g!.\]
\end{theorem}

Note that in contrast to our bound on geometric torsion, this bound works at primes of good reduction.  This does not give a uniform bound because the smallest prime of good reduction can be arbitrarily high, for example, a hyperelliptic curve of genus $2$ whose ramification points over $\P^1$ come together at all primes less than a given large number.

For ease of notation, we explain Buium's argument in the case where
$K=\Q$.  Let $p$ be a prime of good reduction for $X$.  For a perfect field $k$ of characteristic $p$ and a scheme $X$ locally of finite type over $W_{n+1}(k)$, the ring of the length $n+1$ Witt vectors, there is a scheme $X_0^n$ such that for any $k$-algebra $L$, there is a functorial identification $X(W_{n+1}(L))=X_0^n(L)$.   The scheme $X_0^n$ is called the Greenberg transform of $X$.  In the case under consideration, we will let $k=\F_p$ and $n=1$, so $X_0^1$ satisfies, in particular $X(\Z/p^2)=X_0^1(\F_p)$.   There is a natural reduction map $\pi\colon X(\Q_p^{\unr})\to X_0^n(\overline{\F}_p)$ where $\Q_p^{\unr}$ is the maximal unramified extension of $\Q_p$.

By general considerations involving algebraic groups, the Greenberg transform of
the Jacobian $J$ fits in an exact sequence
\[\xymatrix{
0\ar[r]&L\ar[r]&J_0^1\ar[r]&J_0\ar[r]&0}
\]
where $L\cong \G_a^g$ is a vector group over $\F_p$.  Because $X$ has good reduction at $p$, $J_0$ is proper.  It follows that $pJ_0^1$, the subgroup of $J_0^1$ given by the image of the multiplication-by-$p$ map, is proper and the quotient $J_0^1/pJ_0^1$ is isomorphic to $\G_a^g$.

Now we make use of the following lemma:
\begin{lemma} If $X$ is a curve of genus $g$ over $\Q_p$ with good reduction and $p>2g$ then $X\cap J_{\tors}\subset J(\Q_p^{\unr})_{\tors}$.
\end{lemma}

This lemma is due to Coleman \cite{coleman:ramified_torsion_curves}.  We can combine it with the following observation which follows from basic facts about formal groups (see, e.g.,~\cite{katz:galois-properties-of-torsion})

\begin{lemma} The reduction map $\pi\colon J(\Q_p^{\unr})_{\tors}\to J^1_0(\overline{\F}_p)$ is injective.
\end{lemma}

We summarize with the diagram:
\[
\xymatrix{
J(\Q_p^{\unr})_{\tors}\ar[r]&J_0^1(\overline{\F}_p)\ar[d]\ar[r]&J_0^0(\overline{\F}_p).\\
& J_0^1/pJ_0^1\cong \G_a^g&
}
\]

Now, by bounding the rank of the $p$-power torsion, we can conclude that that the image of $J(\Q_p^{\unr})_{\tors}$ in $J_0^1/pJ_0^1$ has cardinality at most $p^{2g}$.  Consequently if $b_1,\dots,b_N$ are elements of $J(\Q_p^{\unr})_{\tors}$ that surject onto the image, we have the following inclusion of sets in $J_0^1$:
\[\pi(\iota(X)\cap J(\overline{\Q}_p)_{\tors})\subset \bigcup_i \left(X_0^1\cap (b_i+pJ_0^1\right)).\]
Buium shows using a very clever argument that $X_0^1$ is affine, and because $pJ_0^1$ is proper, the intersections  $X_0^1\cap (b_i+pJ_0^1)$ are all finite.  The cardinality of the intersection can be bounded by an intersection theory computation \cite[Proof of Theorem 1.11]{buium:p-adic-jets}.  

It is an open question to combine our work with Buium's to get an
unconditional uniform bound on geometric torsion.  In forthcoming work
with Taylor Dupuy \cite{DKRZB}, we will extend Buium's argument and combine it with our annular bound to obtain a geometric torsion bound for curves of compact-type reduction.  For more general reduction types, the delicate balance between properness and affineness is upset, and some new ideas are needed.

\subsection{Kim's non-Abelian Chabauty program}

The hypothesis that $r < g$ is essential to the analysis of rational points via Chabauty's method and Coleman's theory of (abelian) $p$-adic integration. An especially inventive new direction to circumvent this limitation is M. Kim's ``non-abelian Chabauty'' \cite{Kim:siegel}, which exploits mildly non-abelian enrichments of the Jacobian of a curve and \emph{iterated}, rather than abelian integration, to derive Diophantine results. The initial successes of \cite{Kim:siegel} include a new proof of Siegel's theorem, and, assuming a few standard conjectures, a new proof of the Mordell conjecture.

There is great interest (and difficulty) in, and voluminous collective activity (e.g., a 2007 Banff workshop) dedicated to, making Kim's methods explicit and adapting his techniques to derive other Diophantine results. One recent major advance is \cite{balakrishnanBS:quadratic-arxiv} which extracts from Kim's work \emph{explicit formulas} for $p$-adic analytic functions vanishing on $X(\Z)$ in the very special case that $X$ is hyperelliptic and $\rank_{\Z} J(\Q) = g$. 

We give a brief exposition of Kim's method following \cite{Kim:fundamentalgroups}.  Let $X$ be an algebraic variety defined over $\Q$ and set $\overline{X}=X\times \overline{\Q}$.  Set $\Gamma=\Gal(\overline{\Q}/\Q)$.  Let $\pi_1^{\et}(X,b)$ be the \'etale fundamental group based at $b\in X(\Q)$.  For any $x\in X(\Q)$, the fundamental torsor $\pi_1^{\et}(\overline{X},b,x)$ gives an element of $H^1(\Gamma,\pi_1^{\et}(\overline{X},b))$ yielding a {\em period map}
\[X(\Q)\to H^1(\Gamma,\pi_1^{\et}(\overline{X},b)).\]
For technical reasons, one may replace the \'{e}tale fundamental group by the Tannakian fundamental group associated to locally constant unipotent $\Q_p$-vector bundles.  Here, a unipotent $\Q_p$-vector bundle is a vector bundle on $X$ equipped with a filtration
\[L=L_0\supset L_1 \supset \dots \supset L_n \supset L_{n+1}=0\]
whose associated graded bundles are trivial $\Q_p$-vector bundles.  Taking the fiber at a point gives a fiber functor, and so one may define a Tannakian fundamental group $\pi_1^{u,\Q_p}(\overline{X},b)$ which is equal to the $\Q_p$-pro-unipotent completion of the \'{e}tale fundamental group.  Again, one has a period map 
 \[X(\Q)\to H^1(\Gamma,\pi_1^{u,\Q_p}(\overline{X},b)).\]
The $\Q_p$-pro-unipotent fundamental group has quotients $[\pi_1^{u,\Q_p}(\overline{X},b))]_n$ corresponding to taking vector bundles with filtrations of length $n$.  By taking a model for $X$, we can refine the period map to
\[X(\Z_S)\to  H_f^1(\Gamma,\pi_1^{u,\Q_p}(\overline{X},b)) \subset H^1(\Gamma,\pi_1^{u,\Q_p}(\overline{X},b))\]
where $S$ is a finite set of primes containing the primes of bad reduction of $X$ and the subscript $f$ denotes local conditions at cofinitely many primes.  One may construct a $p$-adic period morphism using $\Gamma_p=\Gal(\overline{\Q_p}/\Q_p)$ and obtain a commutative diagram
\[
\xymatrix{
X(\Z_S)\ar[d]\ar[r]&X(\Z_p)\ar[d]\\
H_f^1(\Gamma,[\pi_1^{u,\Q_p}(\overline{X},b)]_n)\ar[r]&H_f^1(\Gamma_p,[\pi_1^{u,\Q_p}(\overline{X},b)]_n)}.
\]
Here, the lower row should be thought of as analogous to the inclusion $J(\Q)\hookrightarrow J(\Q_p)$ (technically its image under $p$-adic logarithm).
Now, one has the isomorphism $H_f^1(\Gamma_p,[\pi_1^{u,\Q_p}(\overline{X},b)]_n)\cong [\pi^{DR}_1(X_{\Q_p},b)]_n/F^0$ where $F^0$ is the $0$th piece of the Hodge filtration on the truncated unipotent de Rham fundamental group.  Subject to deep conjectures on Galois representations of geometric origins, in the case of curves one expects the map $H_f^1(\Gamma,[\pi_1^{u,\Q_p}(\overline{X},b)]_n)\to H_f^1(\Gamma_p,[\pi_1^{u,\Q_p}(\overline{X},b)]_n)\cong [\pi^{DR}_1(X_{\Q_p},b)]_n/F^0$ to have  image equal to a proper subset.   Now, the target space has coordinates given by $p$-adic iterated integrals
\[\int_b^x \omega_1\omega_2\cdots\omega_m.\]
Consequently, one would obtain a locally analytic function vanishing on the $S$-integral points.

The closely related quadratic Chabauty method introduced by Balakrishnan, Besser, and M\"{u}ller \cite{balakrishnanBS:quadratic-arxiv} can be used to study $p$-integral points on hyperelliptic curves whose Mordell--Weil rank is equal to their genus.  The authors produce a function that takes on finitely many values on the $p$-integral points.  Using this function, they are able to give an explicit method for computing the $p$-integral points.  The function arises from iterated $p$-adic integrals and the Coleman--Gross height pairing.

It is an interesting and difficult question to give explicit uniform bounds for the non-abelian or quadratic Chabauty methods.  Leaving aside the conjectures on Galois representations, one needs to bound the zeros of a $p$-adic iterated integral.  Locally, there are considerable difficulties in maintaining control over the Newton polygon because one is multiplying functions (which have an uncontrolled constant terms coming from constants of integration).  Here, having alternative interpretations of the iterated integrals would be helpful.
More globally, there is the same issue of the divergence between the analogues
of abelian and Berkovich--Coleman integrals in the bad reduction case.  Work in
progress by Besser--Zerbes addresses this question.  The situation is complicated
because iterated integrals are not additive along paths but instead obey a more
involved concatenation formula.

\subsection{The Poonen--Stoll sieve}

 The Poonen--Stoll sieve \cite{PoonenS:mostOddDegree} is a technique that allowed Poonen and Stoll to prove that most odd degree hyperelliptic curves have only one rational point.  Here, ``most'' is defined with respect to the family $\cF_g$ of hyperelliptic curves 
 \[y^2=x^{2g+1}+a_1x^{2g}+a_2x^{2g-1}+\dots+a_{2g}x+a_{2g+1}.\]
  A height function is given by
 \[H(C)=\max\left\{|a_1|,|a_2|^{1/2},\dots,|a_{2g+1}|^{1/(2g+1)} \right\} \]
 and for $X\in \R$,
 \[\cF_{g,X}\coloneq\left\{C\in \cF_g \mid H(C)<X\right\}.\]
 The density of a subset $S\subseteq \cF_g$ is 
 \[\mu(S)=\lim_{X\to\infty} \frac{\#(S\cap \cF_{g,X})}{\#\cF_{g,X}}.\]
 The main result is that the lower density of the set of curves in $\cF_g$ with exactly one rational point is a positive number that tends to $1$ as $g$ goes to infinity.  The proof makes use of the Bhargava--Gross theorem on the equidistribution of $2$-Selmer elements \cite{bhargava2012average}.

We will discuss only the step in the proof analogous to the Chabauty method.  Recall that the classical Chabauty method uses the $p$-adic closure $\overline{J(\Q)}\subset J(\Q_p)$ as a stand-in for the rational points of the Jacobian.  The Poonen--Stoll method has a stand-in coming from the Selmer group.  For a field $K$, we apply Galois cohomology of $\overline{K}/K$ to the short exact sequence
\[\xymatrix
{0\ar[r]&J[p](\overline{K})\ar[r]&J(\overline{K})\ar[r]^{\cdot p}& J(\overline{K})\ar[r]&0}
\]
to obtain a connecting homomorphism $\frac{J(K)}{pJ(K)}\hookrightarrow H^1(K,J[p]).$
We make use of the inclusion of $\Q$ into its $v$-adic completions to obtain a commutative diagram
\[
\xymatrix{
\frac{J(\Q)}{pJ(\Q)}\ar[d]\ar@{^{(}->}[r]^>>>>>>>>\delta &H^1(K,J[p])\ar[d]^{\res}\\
\prod_v \frac{J(\Q_v)}{pJ(\Q)_v}\ar@{^{(}->}[r]^>>>>{\delta'} & \prod_v H^1(\Q_v,J[p]).
}\]
The $p$-Selmer group is defined to be
\[\Sel_p J=\{\xi\in H^1(\Q,J[p]) \mid \res({\xi})\in\im(\delta')\}.\]
By construction, there is an inclusion $J(\Q)/pJ(\Q)\hookrightarrow \Sel_p J$, and there is also a natural homomorphism $\Sel_p J\to J(\Q_p)/pJ(\Q_p)$ coming from $\res$.  
The image of $\Sel_p J$ in $J(\Q_p)/pJ(\Q_p)$ is a stand-in for the image of $J(\Q)/pJ(\Q)$.  
The logarithm map $\log\colon J(\Q_p)\to \Z_p^g$ induces a homomorphism $\log\otimes \F_p\colon J(\Q_p)/pJ(\Q_p)\to \F_p^g$.  We will projectivize and consider the intersection of the image of $X(\Q_p)$ with the image of $\Sel_p J$ in $\P^{g-1}(\F_p)$.  We have a commutative diagram
\[\xymatrix{
X(\Q)\ar@{^{(}->}[rr]\ar@{^{(}->}[d]&&X(\Q_p)\ar@{^{(}->}[d]&&\\
J(\Q)\ar[d]\ar@{^{(}->}[rr]&&J(\Q_p)\ar[d]_{\pi_p}\ar[r]^\log&\Z_p^g\ar[d]\ar@{-->}[dr]^\rho&\\
\frac{J(\Q)}{pJ(\Q)}\ar@{^{(}->}[r]&\Sel_p J \ar[r] \ar@/_1.3pc/[rr]_\sigma&\frac{J(\Q_p)}{pJ(\Q_p)}\ar[r]^{\log\otimes \F_p}&\F_p^g\ar@{-->}[r]&\P^{g-1}(\F_p).
}\]
where $\sigma$ is defined by composition.

The analogue of the Chabauty method is the following statement:
\begin{proposition} Suppose that $\sigma$ is injective.  Then we have $\overline{J(\Q)}[p^\infty]=0$ and $\rho\log(\overline{J(\Q)})\subseteq \P\sigma(\Sel_p J)$.  If, in addition, $\rho\log(C(\Q_p))$ and $\P\sigma(\Sel_p J)$ are disjoint, then $C(\Q_p)\cap \overline{J(\Q)}\subseteq J(\Q_p)[p']$.
\end{proposition}

Here, $J(Q_p)[p']$ denotes the prime-to-$p$ torsion.  Because $\rho$ is not defined on all of the logarithm of torsion, we cannot conclude that $C(\Q)=\emptyset$ as one might initially expect.   The main result follows from restricting to the case $p=2$ and making use of the Bhargava--Gross result to show that on a set of curves of positive density $\sigma$ is injective and $\rho\log(C(\Q_p))\cap\P\sigma(\Sel_p J)=\emptyset$ .  Therefore, the only rational points are torsion.  By cutting away small sets from $\cF_g$, one may ensure that the only torsion points are Weierstrass points.  But the only prime-to-$2$ Weierstrass torsion point is the point at infinity.

The Poonen--Stoll sieve has been exploited by Stoll \cite{Stoll:chabautywoMordellWeil} to give an algorithm for finding rational points on curves.  The algorithm, which requires conditions on the maps $\sigma$ and $\pi_p$ to succeed, does not require finding generators of a finite-index subgroup of the Mordell--Weil group.  The difficulty of that step is the major obstacle to executing the usual explicit Chabauty method in practice.  

\bibliographystyle{amsalpha}
\bibliography{master}

\def\cprime{$'$}
\providecommand{\bysame}{\leavevmode\hbox to3em{\hrulefill}\thinspace}
\providecommand{\MR}{\relax\ifhmode\unskip\space\fi MR }
\providecommand{\MRhref}[2]{%
  \href{http://www.ams.org/mathscinet-getitem?mr=#1}{#2}
}
\providecommand{\href}[2]{#2}
\begin{thebibliography}{KRZB16}

\bibitem[AB15]{Amini2012}
Omid Amini and Matthew Baker, \emph{Linear series on metrized complexes of
  algebraic curves}, Math. Ann. \textbf{362} (2015), no.~1-2, 55--106.
  \MR{3343870}

\bibitem[Bak08]{Baker:specialization}
Matthew Baker, \emph{Specialization of linear systems from curves to graphs},
  Algebra Number Theory \textbf{2} (2008), no.~6, 613--653, With an appendix by
  Brian Conrad. \MR{2448666 (2010a:14012)}

\bibitem[BBK10]{balakrishnanBK:coleman}
Jennifer~S. Balakrishnan, Robert~W. Bradshaw, and Kiran~S. Kedlaya,
  \emph{Explicit {C}oleman integration for hyperelliptic curves}, Algorithmic
  number theory, Lecture Notes in Comput. Sci., vol. 6197, Springer, Berlin,
  2010, pp.~16--31. \MR{2721410 (2012b:14048)}

\bibitem[BBM16]{balakrishnanBS:quadratic-arxiv}
Jennifer~S. Balakrishnan, Amnon Besser, and J.~Steffen Müller, \emph{Quadratic
  {C}habauty: {$p$}-adic heights and integral points on hyperelliptic curves},
  J. Reine Angew. Math. \textbf{720} (2016), 51--79. \MR{3565969}

\bibitem[Ber90]{Berkovich:spectralTheory}
Vladimir~G. Berkovich, \emph{Spectral theory and analytic geometry over
  non-{A}rchimedean fields}, Mathematical Surveys and Monographs, vol.~33,
  American Mathematical Society, Providence, RI, 1990. \MR{1070709 (91k:32038)}

\bibitem[Ber07]{Berkovich:integrationOfOne}
\bysame, \emph{Integration of one-forms on {$p$}-adic analytic spaces}, Annals
  of Mathematics Studies, vol. 162, Princeton University Press, Princeton, NJ,
  2007. \MR{2263704}

\bibitem[BF11]{BakerF:metric}
Matthew Baker and Xander Faber, \emph{Metric properties of the tropical
  {A}bel-{J}acobi map}, J. Algebraic Combin. \textbf{33} (2011), no.~3,
  349--381. \MR{2772537 (2012c:14124)}

\bibitem[BG13]{bhargava2012average}
Manjul Bhargava and Benedict~H. Gross, \emph{The average size of the 2-{S}elmer
  group of {J}acobians of hyperelliptic curves having a rational {W}eierstrass
  point}, Automorphic representations and {$L$}-functions, Tata Inst. Fundam.
  Res. Stud. Math., vol.~22, Tata Inst. Fund. Res., Mumbai, 2013, pp.~23--91.
  \MR{3156850}

\bibitem[BGW13]{bhargava2013pencils}
Manjul Bhargava, Benedict~H Gross, and Xiaoheng Wang, \emph{Pencils of quadrics
  and the arithmetic of hyperelliptic curves}, Preprint available at
  arXiv:1310.7692 (2013).

\bibitem[Bha13]{bhargava2013most}
Manjul Bhargava, \emph{Most hyperelliptic curves over $\mathbb{Q}$ have no
  rational points}, Preprint available at arXiv:1308.0395 (2013).

\bibitem[BJ16]{bakerJ:degeneration-survey}
Matthew Baker and David Jensen, \emph{Degeneration of linear series from the
  tropical point of view and applications}, Nonarchimedean and Tropical
  Geometry, Simons Symp., Springer, 2016.

\bibitem[BL84]{BoschL:stableReductionAnd2}
Siegfried Bosch and Werner L{\"u}tkebohmert, \emph{Stable reduction and
  uniformization of abelian varieties. {II}}, Invent. Math. \textbf{78} (1984),
  no.~2, 257--297. \MR{0767194 (86j:14040b)}

\bibitem[BL85]{BoschL:stableReductionAnd}
\bysame, \emph{Stable reduction and uniformization of abelian varieties. {I}},
  Math. Ann. \textbf{270} (1985), no.~3, 349--379. \MR{0774362 (86j:14040a)}

\bibitem[BL91]{BoschL:degeneratingAbelianVarieties}
\bysame, \emph{Degenerating abelian varieties}, Topology \textbf{30} (1991),
  no.~4, 653--698. \MR{1133878 (92i:14043)}

\bibitem[BN07]{BakerN:RR}
Matthew Baker and Serguei Norine, \emph{Riemann-{R}och and {A}bel-{J}acobi
  theory on a finite graph}, Adv. Math. \textbf{215} (2007), no.~2, 766--788.
  \MR{2355607 (2008m:05167)}

\bibitem[BN09]{BakerN:Harmonic}
\bysame, \emph{Harmonic morphisms and hyperelliptic graphs}, Int. Math. Res.
  Not. IMRN (2009), no.~15, 2914--2955. \MR{2525845 (2010e:14031)}

\bibitem[Bom90]{Bombieri:mordell-conjecture-revisited}
Enrico Bombieri, \emph{The {M}ordell conjecture revisited}, Ann. Scuola Norm.
  Sup. Pisa Cl. Sci. (4) \textbf{17} (1990), no.~4, 615--640. \MR{1093712
  (92a:11072)}

\bibitem[Bou98]{Bourbaki-lie-groups-1-3}
Nicolas Bourbaki, \emph{Lie groups and {L}ie algebras. {C}hapters 1--3},
  Elements of Mathematics (Berlin), Springer-Verlag, Berlin, 1998, Translated
  from the French, Reprint of the 1989 English translation. \MR{1728312
  (2001g:17006)}

\bibitem[Bou05]{Bourbaki-lie-groups-7-9}
\bysame, \emph{Lie groups and {L}ie algebras. {C}hapters 7--9}, Elements of
  Mathematics (Berlin), Springer-Verlag, Berlin, 2005, Translated from the 1975
  and 1982 French originals by Andrew Pressley. \MR{2109105 (2005h:17001)}

\bibitem[BPR13]{Baker2013}
Matthew Baker, Sam Payne, and Joseph Rabinoff, \emph{On the structure of
  nonarchimedean analytic curves}, Tropical and {N}on-{A}rchimedean {G}eometry,
  Contemp. Math., vol. 605, Amer. Math. Soc., Providence, RI, 2013,
  pp.~93--121. \MR{3204269}

\bibitem[BR10]{Baker2010}
Matthew Baker and Robert Rumely, \emph{Potential theory and dynamics on the
  {B}erkovich projective line}, Mathematical Surveys and Monographs, vol. 159,
  American Mathematical Society, Providence, RI, 2010. \MR{2599526}

\bibitem[BR15]{bakerR:skeleton}
Matthew Baker and Joseph Rabinoff, \emph{The skeleton of the {J}acobian, the
  {J}acobian of the skeleton, and lifting meromorphic functions from tropical
  to algebraic curves}, Int. Math. Res. Not. IMRN (2015), no.~16, 7436--7472.
  \MR{3428970}

\bibitem[BS10]{BruinS:2010-mordell-weil-sieve}
Nils Bruin and Michael Stoll, \emph{The {M}ordell-{W}eil sieve: proving
  non-existence of rational points on curves}, LMS J. Comput. Math. \textbf{13}
  (2010), 272--306. \MR{2685127}

\bibitem[BS13]{bhargava2013average}
Manjul Bhargava and Arul Shankar, \emph{The average size of the 5-{S}elmer
  group of elliptic curves is 6, and the average rank is less than 1}, Preprint
  available at arXiv:1312.7859 (2013).

\bibitem[Bui96]{buium:p-adic-jets}
Alexandru Buium, \emph{Geometry of {$p$}-jets}, Duke Math. J. \textbf{82}
  (1996), no.~2, 349--367. \MR{1387233}

\bibitem[CdS88]{CdS}
Robert Coleman and Ehud de~Shalit, \emph{{$p$}-adic regulators on curves and
  special values of {$p$}-adic {$L$}-functions}, Invent. Math. \textbf{93}
  (1988), no.~2, 239--266. \MR{948100 (89k:11041)}

\bibitem[Cha41]{Chabauty1941}
Claude Chabauty, \emph{Sur les points rationnels des courbes alg\'ebriques de
  genre sup\'erieur \`a l'unit\'e}, C. R. Acad. Sci. Paris \textbf{212} (1941),
  882--885.

\bibitem[CHM97]{CaporasoHM:uniformity}
Lucia Caporaso, Joe Harris, and Barry Mazur, \emph{Uniformity of rational
  points}, J. Amer. Math. Soc. \textbf{10} (1997), no.~1, 1--35. \MR{1325796
  (97d:14033)}

\bibitem[Col82]{ColemanDilogarithms}
Robert~F. Coleman, \emph{Dilogarithms, regulators and {$p$}-adic
  {$L$}-functions}, Invent. Math. \textbf{69} (1982), no.~2, 171--208.
  \MR{674400 (84a:12021)}

\bibitem[Col85a]{Coleman:effectiveChabauty}
\bysame, \emph{Effective {C}habauty}, Duke Math. J. \textbf{52} (1985), no.~3,
  765--770. \MR{808103 (87f:11043)}

\bibitem[Col85b]{ColemanTorsion}
\bysame, \emph{Torsion points on curves and {$p$}-adic abelian integrals}, Ann.
  of Math. (2) \textbf{121} (1985), no.~1, 111--168. \MR{782557 (86j:14014)}

\bibitem[Col87]{coleman:ramified_torsion_curves}
\bysame, \emph{Ramified torsion points on curves}, Duke Math. J. \textbf{54}
  (1987), no.~2, 615--640. \MR{899407}

\bibitem[Col89]{ColemanReciprocity}
\bysame, \emph{Reciprocity laws on curves}, Compositio Math. \textbf{72}
  (1989), no.~2, 205--235. \MR{1030142 (91c:14028)}

\bibitem[DKRZB]{DKRZB}
Taylor Dupuy, Eric Katz, Joseph Rabinoff, and David Zureick-Brown, \emph{Total
  jet spaces and uniform bounds for torsion points on curves with compact type
  reduction}, in preparation.

\bibitem[DM69]{DeligneMumford:irreducibility}
P.~Deligne and D.~Mumford, \emph{The irreducibility of the space of curves of
  given genus}, Inst. Hautes \'Etudes Sci. Publ. Math. (1969), no.~36, 75--109.
  \MR{0262240 (41 \#6850)}

\bibitem[Fal86]{Faltings:bookArithmeticGeometry}
Gerd Faltings, \emph{Finiteness theorems for abelian varieties over number
  fields}, Arithmetic geometry ({S}torrs, {C}onn., 1984), Springer, New York,
  1986, Translated from the German original [Invent. Math. {{\bf{7}}3} (1983),
  no. 3, 349--366; ibid. {{\bf{7}}5} (1984), no. 2, 381; MR 85g:11026ab] by
  Edward Shipz, pp.~9--27. \MR{861971}

\bibitem[FvdP04]{FdV}
Jean Fresnel and Marius van~der Put, \emph{Rigid analytic geometry and its
  applications}, Progress in Mathematics, vol. 218, Birkh\"auser Boston, Inc.,
  Boston, MA, 2004. \MR{2014891 (2004i:14023)}

\bibitem[GG93]{GordonG:1993-computing-the-mordell-weil}
Daniel~M. Gordon and David Grant, \emph{Computing the {M}ordell-{W}eil rank of
  {J}acobians of curves of genus two}, Trans. Amer. Math. Soc. \textbf{337}
  (1993), no.~2, 807--824. \MR{1094558}

\bibitem[Gou97]{gouvea:pAdicNumbers-an-introduction}
Fernando~Q. Gouv{\^e}a, \emph{{$p$}-adic numbers}, second ed., Universitext,
  Springer-Verlag, Berlin, 1997, An introduction. \MR{1488696}

\bibitem[Hin88]{Hindry1988}
Marc Hindry, \emph{Autour d'une conjecture de {S}erge {L}ang}, Invent. Math.
  \textbf{94} (1988), no.~3, 575--603. \MR{969244}

\bibitem[Ho14]{ho:howManySurvey}
Wei Ho, \emph{How many rational points does a random curve have?}, Bull. Amer.
  Math. Soc. (N.S.) \textbf{51} (2014), no.~1, 27--52. \MR{3119821}

\bibitem[Kat81]{katz:galois-properties-of-torsion}
Nicholas~M. Katz, \emph{Galois properties of torsion points on abelian
  varieties}, Invent. Math. \textbf{62} (1981), no.~3, 481--502. \MR{604840
  (82d:14025)}

\bibitem[Kim05]{Kim:siegel}
Minhyong Kim, \emph{The motivic fundamental group of {$\bold P^1 -
  \{0,1,\infty\}$} and the theorem of {S}iegel}, Invent. Math. \textbf{161}
  (2005), no.~3, 629--656. \MR{2181717 (2006k:11119)}

\bibitem[Kim10]{Kim:fundamentalgroups}
\bysame, \emph{Fundamental groups and {D}iophantine geometry}, Cent. Eur. J.
  Math. \textbf{8} (2010), no.~4, 633--645. \MR{2671216}

\bibitem[Kla93]{Klassen:1993-algebraic-points-low-degree}
Matthew~James Klassen, \emph{Algebraic points of low degree on curves of low
  rank}, ProQuest LLC, Ann Arbor, MI, 1993, Thesis (Ph.D.)--The University of
  Arizona. \MR{2690239}

\bibitem[KRZB16]{KatzRZB-uniform-bounds}
Eric Katz, Joseph Rabinoff, and David Zureick-Brown, \emph{Uniform bounds for
  the number of rational points on curves of small {M}ordell--{W}eil rank},
  Duke Math. J. \textbf{165} (2016), no.~16, 3189--3240. \MR{3566201}

\bibitem[KZB13]{katzZB:tropicalChabauty}
Eric Katz and David Zureick-Brown, \emph{The {C}habauty-{C}oleman bound at a
  prime of bad reduction and {C}lifford bounds for geometric rank functions},
  Compos. Math. \textbf{149} (2013), no.~11, 1818--1838. \MR{3133294}

\bibitem[LT02]{LoreniniT:thue}
Dino Lorenzini and Thomas~J. Tucker, \emph{Thue equations and the method of
  {C}habauty-{C}oleman}, Invent. Math. \textbf{148} (2002), no.~1, 47--77.
  \MR{1892843 (2003d:11088)}

\bibitem[MP12]{McCallumP:chabautySurvey}
William McCallum and Bjorn Poonen, \emph{The method of {C}habauty and
  {C}oleman}, Explicit methods in number theory, Panor. Synth\`eses, vol.~36,
  Soc. Math. France, Paris, 2012, pp.~99--117. \MR{3098132}

\bibitem[MZ08]{Mikhalkin2008}
Grigory Mikhalkin and Ilia Zharkov, \emph{Tropical curves, their {J}acobians
  and theta functions}, Curves and abelian varieties, Contemp. Math., vol. 465,
  Amer. Math. Soc., Providence, RI, 2008, pp.~203--230. \MR{2457739}

\bibitem[Par16]{Park-symmetric-chabauty-arxiv}
Jennifer Park, \emph{Effective {C}habauty for symmetric powers of curves},
  arXiv:1606.05195 (2016).

\bibitem[Poo01]{poonen:computing-torsion}
Bjorn Poonen, \emph{Computing torsion points on curves}, Experiment. Math.
  \textbf{10} (2001), no.~3, 449--465. \MR{1917430 (2003k:11104)}

\bibitem[PS14]{PoonenS:mostOddDegree}
Bjorn Poonen and Michael Stoll, \emph{Most odd degree hyperelliptic curves have
  only one rational point}, Ann. of Math. (2) \textbf{180} (2014), no.~3,
  1137--1166. \MR{3245014}

\bibitem[PZ08]{Pila2008}
Jonathan Pila and Umberto Zannier, \emph{Rational points in periodic analytic
  sets and the {M}anin-{M}umford conjecture}, Atti Accad. Naz. Lincei Cl. Sci.
  Fis. Mat. Natur. Rend. Lincei (9) Mat. Appl. \textbf{19} (2008), no.~2,
  149--162. \MR{2411018 (2009d:11110)}

\bibitem[Ray83]{raynaud:maninMumford}
M.~Raynaud, \emph{Courbes sur une vari\'et\'e ab\'elienne et points de
  torsion}, Invent. Math. \textbf{71} (1983), no.~1, 207--233. \MR{688265
  (84c:14021)}

\bibitem[Sik09]{Siksek:2009-symmetric-power-chabauty}
Samir Siksek, \emph{Chabauty for symmetric powers of curves}, Algebra Number
  Theory \textbf{3} (2009), no.~2, 209--236. \MR{MR2491943 (2010b:11069)}

\bibitem[Sko34]{skolem-1934-verfahren}
Thoralf Skolem, \emph{Ein verfahren zur behandlung gewisser exponentialer
  gleichungen und diophantischer gleichungen}, C. r \textbf{8} (1934),
  163--188.

\bibitem[Sto]{stoll:uniform}
Michael Stoll, \emph{Uniform bounds for the number of rational points on
  hyperelliptic curves of small {M}ordell--{W}eil rank}, to appear in the
  Journal of the European Math Society.

\bibitem[Sto06]{Stoll:2006-independence-of-rational-points}
\bysame, \emph{Independence of rational points on twists of a given curve},
  Compos. Math. \textbf{142} (2006), no.~5, 1201--1214. \MR{2264661
  (2007m:14025)}

\bibitem[Sto15]{Stoll:chabautywoMordellWeil}
\bysame, \emph{Chabauty without the {M}ordell--{W}eil group}, preprint
  arXiv:1506.04286 (2015).

\bibitem[SW13]{shankar2013average}
Arul Shankar and Xiaoheng Wang, \emph{Average size of the 2-{S}elmer group of
  {J}acobians of monic even hyperelliptic curves}, Preprint available at
  arXiv:1307.3531 (2013).

\bibitem[SZ95]{silverbergZarhin:semistableReduction}
A.~Silverberg and Yu.~G. Zarhin, \emph{Semistable reduction and torsion
  subgroups of abelian varieties}, Ann. Inst. Fourier (Grenoble) \textbf{45}
  (1995), no.~2, 403--420. \MR{1343556 (96h:11057)}

\bibitem[Szp85]{szpiro:peu}
Lucien Szpiro, \emph{Un peu d'effectivité}, Astérisque (1985), no.~127,
  275--287, Seminar on arithmetic bundles: the Mordell conjecture (Paris,
  1983/84). \MR{801928}

\bibitem[Tho15]{Thorne:E6-arithmetic}
Jack~A. Thorne, \emph{{$E_6$} and the arithmetic of a family of
  non-hyperelliptic curves of genus 3}, Forum Math. Pi \textbf{3} (2015), e1,
  41. \MR{3298319}

\bibitem[Thu05]{Thuillier2005}
A.~Thuillier, \emph{Th{\'e}orie du potentiel sur les courbes en
  g{\'e}om{\'e}trie analytique non archim{\'e}dienne. {A}pplications {\`a} la
  th{\'e}orie d'{A}rakelov}, Ph.D. thesis, University of Rennes, 2005, Preprint
  available at
  \url{http://tel.archives-ouvertes.fr/docs/00/04/87/50/PDF/tel-00010990.pdf}.

\bibitem[Ull98]{Ullmo1998}
Emmanuel Ullmo, \emph{Positivit\'e et discr\'etion des points alg\'ebriques des
  courbes}, Ann. of Math. (2) \textbf{147} (1998), no.~1, 167--179. \MR{1609514
  (99e:14031)}

\bibitem[Voj91]{Vojta:siegels-theorem-in-the-compact-case}
Paul Vojta, \emph{Siegel's theorem in the compact case}, Ann. of Math. (2)
  \textbf{133} (1991), no.~3, 509--548. \MR{1109352 (93d:11065)}

\end{thebibliography}

\vspace{.2in}

\end{document}